\documentclass{article}
\usepackage{latexsym}
\usepackage{epsfig}  
\begin{document}
\fontsize{13}{13}\selectfont
\title{\huge The Family BlowUp Formula of the Family
Seiberg-Witten Invariants}
\author{Ai-Ko Liu\footnote
{email address: akliu@math.berkeley.edu, Current Address: 
Mathematics Department of U.C. Berkeley}}

\maketitle


\newtheorem{conj}{Conjecture}[section]
\newtheorem{main}{Main Theorem}[section]
\newtheorem{theo}{Theorem}[section]
\newtheorem{lemm}{Lemma}[section]
\newtheorem{prop}{Proposition}[section]
\newtheorem{rem}{Remark}[section]
\newtheorem{cor}{Corollary}[section]
\newtheorem{mem}{Example}[section]
\newtheorem{defin}{Definition}[section]
\newtheorem{axiom}{Axiom}[section]
\newtheorem{obs}{Observation}
\newtheorem{assump}{Assumption}
\newtheorem{summ}{Summary}
\newtheorem{theo-corr}{Theorem-Corollary}[section]
\newtheorem{axioms}{Axiom:}[section]
\newtheorem{assum}{Assumption} 
\newtheorem{warn}{Warning}

In this paper one studies the family blowup formula of the family 
 Seiberg-Witten invariants[LL]. The paper is one among a series
 of papers aiming at studying the family Seiberg-Witten invariants.
  The family blowup formula has several interesting applications in
 symplectic geometry and in enumerative geometry.  The applications of
 the formula will be presented in the other papers.  Recently, the author
 has applied the technique of the family blowup formula to resolve some
 conjecture [Liu1] by G$\ddot o$ttsche [Got] and G$\ddot o$ttsche-Yau-Zaslow [YZ], 
[Got]. Furthermore, a proof of the Harvey-Moore conjecture [Liu2] is given 
along the
 line of [Liu1], in which the existence of family blowup formula has played 
an essential role. 
The current paper
 contains the material which provides the foundation of these
applications. It is also interesting to compare the results in this
paper with the one in [LL1] about the wall crossing formula of the
 family Seiberg-Witten invariants.

 The derivation of blowup formula has a long history in Donaldson theory.
 After being conjectured by the various experts about its existence, it was
 calculated first by R. Stern and R. Fintushal[FS2] the universal
 formula. Soon after the Seiberg-Witten theory had been developed in [W], the
 much more simplified blowup formula were derived by the various
 experts immediately.
  Despite of its simplicity, it has played a very crucial role in
 understanding the four-manifold topology.  It was in the two
 papers[LL1],[LL2], T.J. Li and the author found out the link between blowup
 formula and certain discrepancy of Taubes' ``SW=Gr'' in the $b^+_2=1$
 category.  Later it was D. Mcduff who modified the definition of
 Gromov-Taubes invariants in this special case and proposed a modified
 definition of the Gromov-Taubes invariant
 that was believed to be identifiable to the Seiberg-Witten
 invariants.

 It was in [LL1] that the family Seiberg-Witten invariants were defined
 and studied by the current author and T.J. Li.  Soon after our study,
 it was found that the family invariants shared the same discrepancy
 as the ordinary Seiberg-Witten invariant of $b^+_2=1$ four-manifold. 
   It was proposed by the author in a discussion with T.J. Li to use
 the blowup formula in studying this phenomena.  It turns out that the
 formula has some rather interesting applications in enumerative
 geometry, too. 
  The details will be presented elsewhere[Liu1].

 The author wants to thank
 Prof. C. H. Taubes and Prof. S.T. Yau for their encouragement. 
 The author also likes to thank T.J. Li with
 whom the theory of family Seiberg-Witten invariants were jointly
 developed [LL1].

 The organization of the current paper is as following. 
 In the first section, we set up the family blowup construction and
 set up the notations that will be frequently used in the following sections.
 In section \ref{section; proof}, one derives the family blowup formula in the
 ${\cal C}^{\infty}$ category using the language of $spin^c$ spinors and 
connections.
  The readers with an algebraic background can skip over the derivation in 
section \ref{section; proof} and
 jump to section \ref{section; define}. 

 After deriving the blowup formula, we outline a few applications 
of the blowup formula in the various
sub-sections of section \ref{section; proof}. In sub-section \ref{subsection; uf}, 
the application to enumeration of
 singular curves with prescribed singular multiplicities [Liu1] is
 addressed in certain detail.

  In section \ref{section; define}, one develops a version of 
algebraic Seiberg-Witten invariants ${\cal ASW}$ for algebraic surfaces. 
The definitions are 
separated into cases and are discussed in different sub-sections. 
 In sub-section \ref{subsection; fsw} we discuss the relationship between
 ${\cal ASW}$ and the usual Seiberg-Witten invariants.
 In section \ref{section; ap}, we prove the family blowup formula for the 
algebraic Seiberg-Witten invariants.
 Finally, in the subsection \ref{subsection; universal}, we construct the 
universal obstruction bundles for the
 universal families. 

\medskip
   
 As a preliminary, let us start by stating the basic facts about family 
Seiberg-Witten invariants and
 $spin^c$ structures.

 A $spin^c$ structure on a four-manifold $M$ determines a $U(2)$ bundle over $M$.
 Following the usual convention, we use its determinant line bundle ${\cal L}$ to 
parametrize the $spin^c$ structures on $M$. Thus, $spin^c$ structures on $M$ 
can be 
 identified non-canonically with (up to torsions) $H^2(M, {\bf Z})$.

 The Seiberg-Witten invariant on $M$ are defined using the Seiberg-Witten moduli 
spaces
  with the expected 
dimension formula $d_{SW}({\cal L})={c_1({\cal L})^2-2\chi(M)-3\sigma(M)\over 4}$
 on
 the moduli spaces.  For $b^+_2>1$ manifolds, the invariants defined are 
diffeomorphism
 invariants of the four-manifolds. For $b^+_2=1$ manifolds, the invariants 
defined 
depends on additional chamber structures of the SW equations.

 The family Seiberg-Witten invariants are a natural generalization of 
Seiberg-Witten invariants
 to fiber bundles ${\cal X}\mapsto B$ of smooth four-manifolds.

 Given a monodromy invariant fiberwise $spin^c$ structure and a fiberwise 
homotopic class of
 section $[B, {\cal X}]_{fiber}$, one may define $FSW$ (see [LL1]).
  
 The expected dimension formula of a relative $spin^c$ structure ${\cal L}$ is
 $$dim_{\bf R}B+{c_1({\cal L})^2-2\chi({\cal X}/B)-3\sigma({\cal X}/B)\over 4}.$$

 If $dim_{\bf R}B < b^+_2-1$, the family invariant defined are independent to
 the relative
 smooth Riemannian metrics and the choices of families of relative self-dual 
two forms used
 to define the family $SW$ equations.

 If $dim_{\bf R}B\geq b^+_2-1$, the family invariant defined may depend on the 
additional 
chamber structures.

\bigskip

 Let us state the main theorems in this paper. The detailed discussion on the 
notations will be
 discussed in section \ref{section; notation}.

 Let ${\cal X}\mapsto B$ be a smooth fiber bundle over
 a smooth oriented even dimensional base $B$ of oriented four-manifolds 
with $b^+_2\geq 1$.
 Let $s:B\mapsto {\cal X}$ be a smooth cross section such that
 the normal bundle ${\bf N}_{s(B)}{\cal X}$ is identified with
 a complex rank two bundle ${\bf N}_s$. Through tubular neighborhood theorem
 it induces  fiberwise almost complex structures in a neighborhood of 
$s(B)\subset {\cal X}$.

Let ${\cal X}'$ be the relative
 almost complex blowing up of ${\cal X}$ along $s(B)$, let $E$ denote
the exceptional line bundle associated to the exceptional locus $\cong
 {\bf P}_B({\bf N}_s)$.  

Let ${\cal L}$ denote a relative $spin^c$ structure of ${\cal X}\mapsto B$
 and let ${\cal L}_0$ denote the pull-back of  ${\cal L}$ by 
$s:B\mapsto {\cal X}$. 
In the additive notation, ${\cal L}+mE$, $m$ odd, represents the
 $spin^c$ structure associated with the tensor product ${\cal L}\otimes
 E^{\otimes m}$.
  
 Fix a fiberwise 
homotopic class of ${\cal C}^{\infty}$ sections $[B, {\cal X}']$, it 
induces a fiberwise homotopic class of sections $[B, {\cal X}]$ through the
 blowing down map ${\cal X}'\mapsto {\cal X}$.  The pure and 
 mixed family invariants
of ${\cal L}$, ${\cal L}+mE$ are defined as in [LL1].

\medskip

\begin{main}(Blowup formula for pure Family invariants)  \label{main; bu}  
Let $m$ be an odd integer, then
 we have the following family blowup formulae relating the pure invariant of
 ${\cal L}+mE$ with the mixed invariants of ${\cal L}$.
   Suppose that $m\geq 3$ and both of the $spin^c$ structures have non-negative 
family
 Seiberg-Witten dimensions, i.e. 
 $$dim_{\bf R}B+{c_1({\cal L})^2-2\chi({\cal X}/B)-3
\sigma({\cal X}/B)\over 4}-{m^2-1\over 4}\geq 0,$$
 then
$$FSW_B(1, {\cal L}+mE)=\sum_{i\geq 0} FSW_B(c_i(\sqrt{{\cal L}_0\otimes
 det({\bf N}_s)^{-1}}\otimes{\bf S}^{m-3\over 2}({\bf N}_s \oplus {\bf C}_B)), 
{\cal L}).$$

 If $m$ is a negative odd integer, then

 $$FSW_B(1, {\cal L}+mE)=\sum_{i\geq 0}FSW_B( c_i(\sqrt{{\cal L}_0\otimes
 det({\bf N}_s)^{-1}}\otimes 
{\bf S}^{-m-3\over 2}({\bf N}_s^{\ast} \oplus {\bf C}_B)), 
{\cal L}).$$

\medskip

Let $m=1$, then $FSW_B(1, {\cal L}+mE)=FSW_B(1, {\cal L})$.
 
\end{main}

 The similar blowup formula of the mixed invariants will be stated and proved
 in the section \ref{section; proof}.

\medskip
 
 In section \ref{section; define}, we define a version of algebraic 
Seiberg-Witten invariants on 
 algebraic surfaces $M$ for
 cohomology classes 
$C={c_1({\cal L}\otimes K_M)\over 2}\in H^{1, 1}(M, {\bf C})\cap H^2(M, {\bf Z})$.

\medskip

\begin{main} \label{main; ASW}
 Let $M$ be an algebraic surface and let $C$ be a integral $(1, 1)$ 
cohomology class, 
then there exists an ${\cal ASW}(C)\in {\bf Z}$ defined in terms of the moduli 
space of
 algebraic curves dual to $C$.

 For $p_g=0$ surfaces, 
the ${\cal ASW}(C)$ can be identified (up to signs) with the usual $SW$ 
invariant of the class 
$2C-c_1(K_M)\in H^2(M, {\bf Z})$ in
the chamber deformed by large multiples of Kahler forms.
\end{main}

\medskip

 The major distinction of algebraic Seiberg-Witten invariants from the usual 
$SW$ on
 symplectic four-manifolds is that
 ${\cal ASW}$ are not of simple type, I.e. ${\cal ASW}(C)$ may be non-zero for
 classes with $d_{GT}(C)={C^2-C\dot c_1(K_M)\over 2}>0$.

\medskip

 The construction for algebraic Seiberg-Witten invariants (based on Kuranishi 
models) 
can be extended to algebraic families
 ${\cal X}\mapsto B$ or its relative blowing up ${\cal X}'\mapsto B$. 
Then we have the
 corresponding family blowing up formulae relating the ${\cal AFSW}$ 
of different classes,

\begin{main}(Blowup Formula for Algebraic Family Seiberg-Witten Invariants). 
\label{main; abu}

 Let $c$ be a class of algebraic cycle in $B$. Then the algebraic mixed invariants
 of $C+mE$ and $C$ are related by 
 
$${\cal AFSW}_{{\cal X}'\mapsto B}(c, C+mE)=\sum_{i\geq 0}
{\cal ASW}_{{\cal X}\mapsto B}(
c\cap c_i({\bf E}\otimes ({\bf S}^{m-2}({\bf C}_B\oplus 
{\bf N}_{s(B)}{\cal X}))), C) $$

 for $m\geq 2$ and 

$${\cal AFSW}_{{\cal X}'\mapsto B}(c, C+mE)=\sum_{i\geq 0}{\cal ASW}_{{\cal X}
\mapsto B}(
c\cap c_i({\bf E}\otimes ({\bf S}^{-m-1}({\bf C}_B\oplus 
{\bf N}^{\ast}_{s(B)}{\cal X}))), C) $$

 for $m\leq -1$.

 For $m=0, 1$, we have $${\cal AFSW}_{{\cal X}'\mapsto B}(c, C+mE)
={\cal AFSW}_{{\cal X}\mapsto B}(c, C).$$

\end{main}

\medskip
 For the assumptions imposed on ${\cal X}', {\cal X}$ and 
the detailed assumptions of
theorem \ref{main; abu}, please consult section 
\ref{section; ap} and theorem \ref{theo; al-bl}.
\medskip

\section{ The family blowing up construction} \label{section; notation}

\medskip

  Suppose that one is given a fiber bundle ${\cal X}$ over a compact oriented 
manifold $B$ 
whose typical
fibers are diffeomorphic to an oriented four-manifold $M$ with $b^+_2>0$.
Let $s:B\mapsto {\cal X}$ be a smooth cross section and  a tubular
neighborhood in ${\cal X}$ is denoted by 
  ${\cal N}$. The tubular neighborhood theorem allows us to
 identify ${\cal N}$ with a real four dimensional vector (ball) bundle 
over $s(B)$. One imposes the extra condition on the section 
$s$ requiring that the real four-plane bundle carries complex structures
and fixes one complex structure in our discussion. Therefore, the
 normal bundle is viewed as a rank two complex vector bundle, denoted by 
${\bf N}_s$.

Let ${\bf C}_B$ be the trivial complex line bundle over $B$
and let $\overline{\bf P}=\overline{{\bf P}_B({\bf N}_s\oplus {\bf C}_B)}$ be the
 projectification of the bundle ${\bf N}_s\oplus {\bf C}_B$, with the fiber-wise
 orientation reversed. The trivial factor 
${\bf C}_B$ in ${\bf N}_s\oplus {\bf C}_B$ 
defines a smooth  
section to $\bar{\bf P}$ whose tubular neighborhood is diffeomorphic to 
$\bar{\bf N}_s$, the 
 total space of the bundle ${\bf N}_s$, with bundle orientation reversed. 
With the preceding convention understood, we can perform the fiberwise
connected sum by deleting the two tubular neighborhoods of ${\cal X}$ and 
 $\bar{\bf P}$ and
gluing their complements via a fiber-wise orientation-reversing
diffeomorphism. The new fiber bundle,
denoted by ${\cal X}'_s$,  ${\cal X}'_s\cong {\cal X}\sharp_f \bar{\bf P}$, 
is called the family blowing up
of ${\cal X}$ along $s$. Unlike the case when $B$ is a point,
the existence of $s$ is no longer a completely trivial matter.
However, if the fiber bundle ${\cal X}\mapsto B$
  has a fiber-wise almost complex structure,
then any section $s$ can be blown up topologically.
Another new feature is that
different choices of fiberwise homotopy classes
 of these cross sections may result in non-diffeomorphic fiber bundles 
${\cal X}'_s$. 

 When the cross section $s$ has been fixed, we may drop the subscript $s$ in 
${\cal X}'_s$
 and denote the resulting fiber bundle by ${\cal X}'$.

\medskip

\begin{rem}\label{rem; total}
If the fiber bundle 
 carries a family of fiberwise symplectic forms (parameterized by the base),
then it induces a fiberwise almost complex structure on the fibers and we can
blow up any smooth cross section. Moreover,
 one can mimic the symplectic blowing up construction of 
Guillemin-Sternberg-Mcduff [Mc] to 
 construct a new family of fiberwise symplectic forms on the ``blown up'' fiber
 bundle.
Even if we have considered the
 family blowing ups in the symplectic category, nevertheless we do not require the
 additional condition that the 
 total space of the fiber bundle to be symplectic. Sometimes the additional
 condition
 is met and the family blowing up construction is really the blowing-up
 of a real codimension four symplectic cross section in a symplectic total space.
  We definitely want to relax the condition here as some natural families (e.g. 
twistor 
 families of $K3$ or $T^4$ or its induced families on the universal spaces [Liu1])
simply do not satisfy this additional condition. However fiberwise blowing ups
 along the cross 
sections of these non-symplectic fiber bundles is crucial in counting holomorphic
 curves in 
the twistor families[Liu1].     
\end{rem}

  It is a well known fact that
 the cohomology of a ${\bf CP}^2$ bundle is, 
 as a module 
 of the cohomology over the base $H^{\ast}(B;{\bf Z})$,
generated by the various powers of ``hyper-plane class''. The 
 choices of the hyper-plane classes are not unique. Different choice of the
hyper-lane class give different generators of the same module. The same
 assertion applies to $\bar{\bf P}$ as well once we flip the orientation. 
 However if we view $\bar{\bf P}$ as the projectification of ${\bf N}_s\oplus
 {\bf C}_B$ with a 
reversed orientation, it does give us a canonical choice of exceptional class $E$.
 From now on
 let us fix this choice implicitly. Then 
 $E^2\cdot \pi^{\ast}_{\bar{\bf P}}[B]=-1$ and $E^3$ lies in  
$\oplus_{i=0}^2 H^{\ast}(B)E^i$.

\medskip

\section{ The Blowup Formula and the Proof of the Family Blowup Formula}
\label{section; proof}

\medskip

 Recall that the usual blowup formula of the Seiberg-Witten invariants
 relates the Seiberg-Witten invariant of a given Spin$^c$ structure on 
$M$ with the invariant of a corresponding Spin$^c$ structure 
 on $M \sharp \overline{ CP^2}$. Before we give a proof of the family
 blowup formula, let us
 review the idea behind the proof of the original formula. 
 The proof of the original formula is based on the usage of the long neck
 metrics. Suppose there is a smooth metric on $M\sharp \overline{ CP^2}$
 such that the metric is isometric to the product metric on the cylindrical
 collar $S^3\times [-L,L]$ between $M-B_{p_1}(\epsilon_1)$ and 
$\overline{ CP^2}-B_{p_2}(\epsilon_2)$. 
 We let $L$ go to $\infty$ and discuss the structure of the moduli 
spaces.  As $S^3$ has positive scalar curvature metrics, the spinor part of the
 Seiberg-Witten solutions tend to 
vanish on the long neck. This implies that any
 solution can be decomposed into one over $M-B_{p_1}(\epsilon_1)$ and the other one
 is over $\overline{ CP^2}-B_{p_2}(\epsilon)$.  
 Conversely we can glue solutions on $M-B_{p_1}(\epsilon_1)$ and 
$\overline{ CP^2}-B_{p_2}(\epsilon_2)$ 
back to  solutions on the connected sum $M\sharp \overline{ CP^2}$.
As $\overline{ CP^2}$ is negative definite, we can choose the metric such that
 the solutions on $\overline {CP^2}$ are reducible. Namely,
 $\psi \equiv 0$ and the connection is anti-self-dual. The reader can consult 
[FS], page 226, theorem
 8.5. where the authors considered long necks of lens spaces $L(p^2, 1-p)$ 
instead of $S^3$. 

The $spin^c$ structure reduces to some odd
multiple of $E$ on $\overline{ CP^2}$. We denote the $spin^c$ determinant 
line bundle on $M$ by ${\cal L}$.
The ${\cal L}$ induces a $spin^c$ determinant line bundle on 
$M\sharp \overline{\bf CP^2}$, also denoted by the
same symbol.

If the absolute value of the multiplicity is equal to one, then
 $dim({\cal M}_{\cal L})=dim({\cal M}_{{\cal L}\pm E})$, 
and 
 $SW({\cal L})=SW({\cal L}\pm E)$, where the left hand side
 is a Seiberg-Witten invariant of $M$ and the right hand side is that of
 $M\sharp \overline{ CP^2}$. 

On the other hand if the multiplicity is bigger than
 one in absolute value, the glued moduli space is smooth 
${\cal M}_{{\cal L}+mE}\cong {\cal M}_{\cal L}$ 
but it is not of the expected dimension.
 Thus the obstruction bundle must be inserted in the calculation of the
 Seiberg-Witten invariants. 

 As usual let $e$ denote the Euler class of the $S^1$ bundle ${\bf e}$ 
over ${\cal M}_{{\cal L}+mE}$ constructed
 from the quotient of global $U(1)$ gauge transformations. Let 
${\bf Obs}_m$ denote the obstruction
 bundle on ${\cal M}_{L+mE}$.

Then the Seiberg-Witten invariant $SW({\cal L}+mE)$
 is calculated by the following expression
$$\int_{{\cal M}_{\cal L}}e^{dim {\cal M}_{{\cal L}+mE}/2} \cup e({\bf Obs}_m)$$
$$=\int_{{\cal M}_{\cal L}}e^{dim {\cal M}_{{\cal L}+mE}/2} 
\cup c_{top}({\bf Obs}_m)$$

where $e({\bf Obs}_m)$ denotes the Euler class of the vector bundle ${\bf Obs}_m$.

 To prove that $SW({\cal L})=SW({\cal L}+mE)$, it suffices to prove that 
$e({\bf Obs}_m)$ is a ${m^2-1\over 8}$ power of $e$.  
 
 We sketch a proof of the following
 simple lemma which is well known to the experts.

\medskip

\begin{lemm}\label{lemm; bun}  
Let $\tilde{\bf e}$ be the complex line bundle over ${\cal M}_{{\cal L}+mE}$ 
induced by the
 principal $S^1$ bundle ${\bf e}$.
The obstruction bundle ${\bf Obs}_m$ over ${\cal M}_{{\cal L}+mE}$ is
 isomorphic to the
 bundle $\tilde{\bf e}\otimes_{\bf C} {\bf C}^{m^2-1\over 8}$.
\end{lemm} 

\noindent Proof of lemma \ref{lemm; bun}:
\medskip

 In fact the formal dimension of the moduli space ${\cal M}_{mE}$ on 
$\overline{\bf CP^2}$
 is ${-m^2-3\over 4}$. On the other hand the expression $2\chi+3\sigma$ 
decreases by one after a single blowing up. Thus the real formal dimensions of the
 moduli spaces $dim {\cal M}_{\cal L}$ and $dim {\cal M}_{{\cal L}+mE}$  
differ by $-m^2+1\over 4$.

Because the unique solution of the $spin^c$ structure $mE$ on 
$\overline{\bf CP^2}$ is reducible, the
 solution is fixed 
 by the $S^1$ action.  Thus the obstruction bundle over ${\cal M}_{mE}$ is 
nothing but a 
${m^2-1\over 8}$ dimensional complex vector space over a single point (the 
complex structure of the
 bundle is inherited from the complex spinors). Under
 the gluing construction the reducible solution on $\overline {CP^2}$ is glued
 to a solution on $M$. As a result, the final $S^1$ action on the 
glued based moduli space comes from the diagonal embedding of the $S^1$ actions
 on both sides of configuration spaces. Under this action, the obstruction
 bundle is identified with
 $\tilde{\bf e}\otimes_{\bf C} {\bf C}^{m^2-1\over 8}$. $\Box$

  After reviewing the idea to derive the 
 blowup formula, we may generalize it
 to the family invariants. Recall that in defining the family invariant, the
 tautological 
class $e$ is not canonically defined. It depends on the choices of a homotopic
 classes of 
 cross sections of the fiber bundles
${\cal X}\mapsto B$, ${\cal X}'\mapsto B$.
 The choice was made implicitly in the definition
 of the family invariants.

Given the data as in section \ref{section; notation}, we relate the family
 Seiberg-Witten
 invariants of ${\cal L}+mE$ on ${\cal X}\sharp_f \bar{\bf P}$ and ${\cal L}$
on ${\cal X}$.
As before the family dimensions of the moduli spaces of these two 
$spin^c$ structures differ by $-m^2+1\over 4$.
 Using the same type of analysis as before, we would like to stretch the
 long neck such that the metric becomes the product metric along the long 
neck. 

\medskip

\begin{lemm} \label{lemm; psc}
The unit sphere bundle of ${\bf N}_s$ carries a fiberwise positive scalar
 curvature metric.
\end{lemm}

\medskip
\noindent Proof:  Given the hermitian metric on ${\bf N}_s$, it induces 
Euclidean metrics 
on the fibers of ${\bf N}_s$. Then the unit spheres in the fibers are given 
the induced 
Riemannian metric over which $SO(4)$ acts transitively. Then the unit spheres 
$\cong S^3$ carries
 the round metric, which is well known to be of positive scalar curvature. $\Box$

Because the $S^3$ bundle carries fiberwise positive scalar curvature metric, one 
can mimic the previous discussion on long neck metric (see also [FS]) 
to decompose or glue the family moduli spaces.
 
 By using the special kind of fiberwise metric, the family moduli spaces for
 ${\cal L}+mE$ is the fiber product of the corresponding  family moduli spaces
 for ${\cal L}$ over ${\cal X}\mapsto B$ and for 
${\cal M}_{\pi^{\ast}_{\bar{\bf P}}
{\cal L}_0+mE}$ over $\bar{\bf P}\mapsto B$, 
where ${\cal M}_{{\cal L}_0+mE}\cong {\cal M}_{mE}$ is diffeomorphic to $B$
 under the natural projection map. Here the line bundle 
${\cal L}_0$ over $B$ denotes the pull-back of  ${\cal L}$ by 
$s:B\mapsto {\cal X}$.
Because ${\cal M}_{mE}$ is not of the expected family dimension,
 it carries an obstruction bundle of real rank ${m^2+3\over 4}$.

 If $m=\pm 1$, we may still conclude that these two family Seiberg-Witten
invariants coincide. $FSW_B({\cal L}\pm E)=FSW_B({\cal L})$. 

When $m\not=\pm 1$, the analogue of lemma \ref{lemm; bun} is 

\begin{prop}  \label{prop; bun2}
 Let ${\cal X}'={\cal X}\sharp_f \bar{\bf P}$ be the fiberwise connected sum of
 ${\cal X}$ with
 $\bar{\bf P}$ with the long neck metric.
 Then there is a real ${m^2-1\over 4}$ dimensional obstruction bundle 
 ${\bf Obs}_m$ over ${\cal M}_{{\cal L}+mE}\cong {\cal M}_{\cal L}$. 

 Let $A_0$ denote the unique fiberwise anti-self-dual connections over
 $\bar{\bf P}\mapsto B$ for
 the $spin^c$ line bundle 
$\pi^{\ast}_{\bar{\bf P}}{\cal L}_0+mE$ with respect to the positive scalar
 curvature 
fiberwise metrics on $\bar{\bf P}$.
 The obstruction bundle ${\bf Obs}_m$ over ${\cal M}_{{\cal L}+mE}$ can be 
identified
 with $\tilde{\bf e}\otimes_{\bf C} \pi^{\ast}_{{\cal M}_{\cal L}}{\bf W}_m$, 
where ${\bf W}_m$ is the 
 ${m^2-1\over 8}$ dimensional complex vector bundle over $B$, the cokernel 
bundle of
 $D_{A_0}$.  
\end{prop}

 The derivation of the proposition is parallel to lemma \ref{lemm; bun} 
except that 
a vector space ${\bf C}^{m^2-1\over 8}$ should be replaced 
by a vector bundle of the same rank.

\noindent Proof of prop \ref{prop; bun2}: 
 The line bundle $\pi^{\ast}_{\bar{\bf P}}{\cal L}_0$ are pull-back from the
 base $B$ and it has a trivial
 first Chern class on the fibers. Its presence does not affect the dimension count.

 Given the unique reducible solution $(A_0, 0)$ on $\bar{\bf P}/B$ of 
$\pi^{\ast}_{\bar{\bf P}}{\cal L}_0+mE$,
 consider the deformation complex of the family Seiberg-Witten equations 
at $(A_0, 0)$,

 $$d^{\ast}\oplus d^+\oplus D_{A_0}:
\Omega^1_{\bar{\bf P}/B}\oplus \Gamma(\bar{\bf P}/B, 
{\cal S}^+_{\pi^{\ast}_{\bar{\bf P}}{\cal L}_0+mE})\longrightarrow
 \Omega^0_{ \bar{\bf P}/B}\oplus \Omega^2_{+, \bar{\bf P}/B}\oplus 
\Gamma(\bar{\bf P}/B, {\cal S}^-_{\pi^{\ast}_{\bar{\bf P}}{\cal L}_0+mE}).$$

 Because ${\bar{\bf P}}$ is simply connected with the 
positive scalar curvature metric, 
 the kernel of the above deformation complex is trivial. 
Because $\bar{\bf P}$ is negative 
definite, $Coker(d^+)=0$. Thus, the cokernel of the deformation complex is
 equal to ${\bf R}_B\oplus Coker(D_{A_0})$, 
where ${\bf R}_B$ denotes the trivial real rank one 
line bundle of constant functions on $B$. This direct factor 
can be identified with the lie algebra
 of the global $S^1$ gauge action which fixes $(A_0, 0)$. 
By index calculation, 
$Coker(D_{A_0})\subset \Gamma(\bar{\bf P}/B, 
{\cal S}^-_{\pi^{\ast}_{\bar{\bf P}}{\cal L}_0+mE})$ 
is of complex rank ${m^2-1\over 8}$.

  After we graft the reducible solution on $\bar{\bf P}$ to 
${\cal X}'$, the non-reducible 
solution on ${\cal X}'$ is not fixed by the $S^1$ action. 
Thus, only the ${\bf W}_m=Coker(D_{A_0})$ factor
 of the obstruction bundle is grafted to an obstruction bundle 
on ${\cal M}_{{\cal L}+mE}$.

 By the same argument as in lemma \ref{lemm; bun}, the complex 
line bundle $\tilde{\bf e}$ (which
depends on the choice of a cross section $B\mapsto {\cal X}'$) 
is tensored with 
$\pi^{\ast}_{{\cal M}_{\cal L}}{\bf W}_m$.

 Thus, we have $${\bf Obs}_m \cong \tilde{\bf e}\otimes 
\pi^{\ast}_{{\cal M}_{\cal L}} {\bf W}_m,
 {\bf W}_m=Coker(D_{A_0}).$$  $\Box$

 Because ${\cal M}_{\pi^{\ast}_{\bar{\bf P}}{\cal L}_0+mE}
\cong {\cal M}_{mE}\cong B$ through the projection map,
we can view ${\bf W}_m$ as a bundle over $B$.
  The following proposition identifies the bundle ${\bf W}_m$ 
in the $K$ group.

\begin{prop} \label{prop; find}
  Given a positive and odd $m$, the vector bundle ${\bf W}_m$ can be
 identified with
 $$\sqrt{{\cal L}_0\otimes det({\bf N}_s)^{-1}}\otimes 
{\bf S}^{m-3\over 2}({\bf N}_s\oplus {\bf C}_B),$$

 in the rational K group $K(B)\otimes_{\bf Z} {\bf Q}$.

 For $m$ odd and negative, the vector bundle ${\bf W}_m$ can be 
identified with
 $$\sqrt{{\cal L}_0\otimes det({\bf N}_s)^{-1}}\otimes 
{\bf S}^{-m-3\over 2}({\bf N}_s^{\ast}\oplus {\bf C}_B),$$
in the rational K group $K(B)\otimes_{\bf Z} {\bf Q}$.
\end{prop}

 The proposition is the key to prove the family blowup formula.

\noindent Proof of proposition \ref{prop; find}: It suffices to 
calculate the Chern character of 
${\bf W}_m$ and show that it is equal to the Chern character of the 
right hand side bundle.
 In the proof of the previous proposition prop. \ref{prop; bun2}, we know that
 $[{\bf W}_m]=[-IND(D_{A_0})]$ in $K(B)$, the negative of family index of 
$D_{A_0}$.

 In general the Chern character of the index bundle is calculated by the
 family index theorem [BGV]. As the calculation is purely topological, it 
does not depend on the explicit choice of the connection $A_0$. From now 
on we ignore the
dependence of the family index virtual bundle on $A_0$.
 The Chern character 
$ch(IND(D_{A_0}))$ is calculated by the following expression,

$$\int_{\bar{\bf P}\mapsto B}\hat{A}_{\bar{\bf P}/B}ch(\sqrt{{\cal L}_0
\otimes \cdot E^{\otimes m}}).$$

 The symbol $\int_{\bar{\bf P}\mapsto B}$ denotes the push-forward map 
from
 $H^{\cdot}(\bar{\bf P})$ to $H^{\cdot}(B)$.

  We notice that ${\bf P}={\bf P}({\bf N}_s\oplus {\bf C}_B)$ and 
$\bar{\bf P}$ have the opposite orientations.
 If we flip the orientation of $\bar{\bf P}$ back to ${\bf P}$ and 
switch the positive and negative spinors as well,

 $$ch(-IND(D_{A_0}))=-\int_{\bar{\bf P}\mapsto B}\hat{A}_{\bar{\bf P}/B}
ch(\sqrt{{\cal L}_0\otimes \cdot E^{\otimes m}})$$
$$=\int_{{\bf P}\mapsto B}\hat{A}_{{\bf P}/B}ch(\sqrt{{\cal L}_0
\otimes \cdot {\bf H}^{\otimes -m}}).$$

 In the above formula, the line bundle $E$ has been replaced by 
${\bf H}^{\ast}$, the tautological line bundle of 
the projective bundle ${\bf P}({\bf N}_s\oplus {\bf C}_B)$.

 To continue the calculation, we apply the following simple lemma,

\begin{lemm}\label{lemm; c1} 
 Let ${\bf V}$ be a complex rank $n$ vector bundle over a smooth manifold $B$. Let
 us denote ${\bf P}({\bf V})$ to be the projective space bundle over $B$ 
formed by projectifying
 ${\bf V}\mapsto B$. Then the first Chern class of the relative tangent bundle 
 along the fibers, $c_1({\bf T}_{{\bf P}({\bf V})/B})$ is given
 by $n{\bf H}+\pi^{\ast}_{{\bf P}({\bf V})}c_1({\bf V})$.
\end{lemm}

\bigskip
\noindent Proof of lemma \ref{lemm; c1}:
 
Recall the following well known short exact sequence of
 ${\bf T}_{{\bf CP}^{n-1}}$,

$$0\rightarrow {\bf C}\rightarrow  
{\bf H} \otimes \pi^{\ast}_{{\bf P}^{n-1}}{\bf C}^n 
\rightarrow {\bf T}_{{\bf CP}^{n-1}} \rightarrow 0.$$

 The relative version of the sequence on ${\bf V}$ gives 

$$0\rightarrow {\bf C}_{{\bf P}({\bf V})} \rightarrow {\bf H}\otimes 
\pi^{\ast}_{{\bf P}({\bf V})} {\bf V} 
 \rightarrow {\bf T}_{{\bf P}({\bf V})/B}\mapsto 0.$$

 Then 

$$c_1({\bf T}_{{\bf P}({\bf V})/B})=c_1({\bf H}\otimes 
\pi^{\ast}_{{\bf P}({\bf V})} {\bf V})
=n{\bf H}+\pi^{\ast}_{{\bf P}({\bf V})}c_1({\bf V}).$$ $\Box$
 
 In our discussion, we take $n=3$ and ${\bf V}={\bf N}_s\oplus {\bf  C}_B$.

 Having identified 
 the relative first Chern class of ${\bf P}$, we are ready to continue the
calculation. 
  Rewrite this family index push-forward as
$$\int_{{\bf P}/B}\hat A_{{\bf P}/B}ch(\sqrt{{\bf H}^3\otimes 
det({\bf N}_s)})
ch(\sqrt{{\cal L}_0\otimes {\bf H}^{\otimes (-m-3)}\otimes 
det({\bf N}_s)^{-1}}).$$
 The relative $\hat A_{{\bf P}/B}$ and $ch(\sqrt{{\bf H}^3\otimes
 det({\bf N}_s)})$ combine into the relative Todd
 class. And the calculation is reduced to a Grothendieck-Riemann-Roch type
 calculation on the projective space bundle ${\bf P}\mapsto B$, 

$$\int_{{\bf P}/B} Todd_{{\bf P}/B}ch(\sqrt{{\cal L}_0\otimes 
{\bf H}^{\otimes (-m-3)}\otimes det({\bf N}_s)^{-1}})$$

$$=ch(IND(\bar\partial+\bar\partial^{\ast}))\cdot 
ch(\sqrt{{\cal L}_0\otimes det({\bf N}_s)^{-1}}).$$

$$\bar\partial+\bar\partial^{\ast}:\Omega^{0,0}_{{\bf P}/B}\otimes 
{\bf H}^{-m-3\over 2}
\oplus \Omega^{0,2}_{{\bf P}/B}\otimes {\bf H}^{-m-3\over 2}\longrightarrow 
 \Omega^{0,1}_{{\bf P}/B}\otimes {\bf H}^{-m-3\over 2}$$

 If $m\leq -3$, ${-m-3\over 2}\in {\bf N}\cup \{0\}$, 
then the above Family Riemann-Roch formula can be re-interpreted as 
(up to tensoring with
 $\sqrt{{\cal L}_0\otimes det({\bf N}_s)^{-1}}$) the Chern character 
of the push-forward of
 ${\bf H}^{-m-3\over 2}$ along the projective bundle ${\bf P}\mapsto B$. 

\medskip
 
 Let us cite the following well known fact on projective spaces,

\begin{lemm}\label{lemm; symmetric}
 Let ${\bf P}(V)$ be the projective space formed by projectifying a 
complex vector space $V$ and ${\bf H}$ denote the
 holomorphic hyperplane bundle on ${\bf P}(V)$ with the standard 
$\bar\partial$ operator, then for $p\geq 0$,
 $H^0_{\bar\partial}({\bf P}(V), {\bf H}^p)$ is naturally isomorphic 
to ${\bf S}^{p}(V^{\ast})$, the $p$th-symmetric power
 of linear functionals on $V$. 
\end{lemm}

\medskip
We adopt the convention that ${\bf S}^0(V)={\bf C}$.

 By applying the family version of lemma \ref{lemm; symmetric}, the 
Chern character of ${\bf W}_m$
  is equivalent to the Chern character of the following
 bundle, 

 $$ \sqrt {{\cal L}_0\otimes det({\bf N}_s)^{-1}} \otimes 
{\bf S}^{-m-3\over 2}
({\bf N}_s^{\ast}\oplus {\bf C}_B).$$

  For $m\geq 3$ we apply relative Serre duality.
 Notice that the usual surface Riemann-Roch theorem [GH]
 for ${\bf CP}^2$ has the following structure,

$$h^0({\bf CP}^2, {\bf H}^p)-h^1({\bf CP}^2, {\bf H}^p)+h^2({\bf CP}^2, 
{\bf H}^p)=1+{(p^2+3p)\over 2}
={(p+1)(p+2)\over 2}.$$

 It is easy to see that Serre duality ${\bf H}^p\mapsto {\bf H}^{-p-3}$ 
induces a 
symmetry on the formula. If $m-3\geq 0$, we replace ${-m-3\over 2}$ by
$-{(-m-3)\over 2}-3={(m-3)\over 2}$. Notice that 

$$H^2({\bf CP}^2, {\bf H}^{-(m+3)\over 2})\cong
 H^0({\bf CP}^2, {\bf H}^{m-3\over 2})^{\ast}.$$

 And in this case the chern character is equivalent to that of 

$$ \sqrt{{\cal L}_0\otimes det({\bf N}_s)^{-1}}\otimes 
{\bf S}^{(m-3)\over 2}({\bf N}_s
\oplus {\bf C}_B).$$

 This ends the proof of proposition \ref{prop; find}. $\Box$ 

We are ready to prove the following the blowup formula of the
 family invariants.

\medskip

\begin{theo}(Blowup formula for the Pure Invariants) \label{theo; blup}
  Let ${\cal X}$, ${\cal X}'$, ${\bf N}_s$ be as defined in 
section \ref{section; notation} and let 
${\cal L}, {\cal L}_0$ be the $spin^c$ determinant line bundle 
over ${\cal X}\mapsto B$ and its 
pull-back by $s:B\mapsto {\cal X}$. 

Let $m>1$ be an odd integer bigger
 such that the $spin^c$ structure ${\cal L}+mE$ has non-negative 
family Seiberg-Witten dimension, then
 we have the following blowup formula relating the pure invariant of
 ${\cal L}+mE$ over ${\cal X}'\mapsto B$ with the mixed invariants of 
${\cal L}$ over ${\cal X}\mapsto B$,

$$FSW_B(1, {\cal L}+mE)=\sum_{i\geq 0} FSW_B(c_i(\sqrt{{\cal L}_0\otimes
 det({\bf N}_s)^{-1}}\otimes {\bf S}^{m-3\over 2}({\bf N}_s\oplus 
{\bf C}_B)),{\cal L}).$$

 Let $m<-1$ be an odd integer, then we have

$$FSW_B(1, {\cal L}+mE)=\sum_{i\geq 0} FSW_B(c_i(\sqrt{{\cal L}_0\otimes
 det({\bf N}_s)^{-1}}\otimes {\bf S}^{-m-3\over 2}({\bf N}_s^{\ast}
\oplus {\bf C}_B)),{\cal L}).$$

  Let $m=\pm 1$, then we have $$FSW_B({\cal L}+mE)=FSW_B({\cal L}).$$
\end{theo}

 \noindent Proof:  The theorem is a consequence of  proposition
 \ref{prop; find}.

  Fix a fiberwise homotopic class of fiberwise metrics and self-dual 
two forms pair on ${\cal X}'$.
  Also fix a fiberwise homotopic class of cross sections $B\mapsto 
{\cal X}'$. Following [LL1], one may
 define the family Seiberg-Witten invariant of ${\cal L}+mE$ in a 
specific chamber.

 For $m=\pm 1$, $$FSW_B(1, {\cal L}\pm E)=\int_{{\cal M}_{\cal L}}
e^{{dim_{\bf R}B\over 2}+
{c_1({\cal L})^2-2\chi(M)-3\sigma(M)\over 8}}=FSW_B(1, {\cal L}).$$

 For $m\not=1$, the pure invariant is equal to

 $$FSW_B(1, {\cal L}+mE)=\int_{{\cal M}_{{\cal L}}}e^{{dim_{\bf R}B\over 2}+
{c_1({\cal L})^2-2\chi(M)-3\sigma(M)-m^2+1\over 8}}
 c_{top}(\tilde{\bf e}\otimes {\bf W}_m).$$

 By proposition \ref{prop; find}, we can replace ${\bf W}_m$ by 
 either 

$$\sqrt{{\cal L}_0\otimes
 det({\bf N}_s)^{-1}}\otimes {\bf S}^{m-3\over 2}({\bf N}_s\oplus {\bf C}_B)$$

 or 

$$ \sqrt{{\cal L}_0\otimes
 det({\bf N}_s)^{-1}}\otimes {\bf S}^{-m-3\over 2}({\bf N}_s^{\ast}
\oplus {\bf C}_B)$$

 depending on $m\geq 3$ or $m\leq -3$. We get the blowup formula by expanding the
 top Chern class in terms of the powers of $c_1(\tilde{\bf e})=e$. $\Box$

 Similarly, by inserting an even dimensional cohomology class $\eta \in
 H^{\ast}(B, {\bf Z})$ into the
 definition of the family
 invariant, we can get the similar formula relating the mixed invariants
 before and after the blowup process.

\begin{theo}(Blowup formula for the Family Mixed Invariants) \label{theo; mixed}
 Suppose that $m\geq 3$, then it follows that
$$FSW_B(\eta,{\cal L}+mE)=\sum_{i\geq 0} FSW_B(\eta\cup c_i({\sqrt 
{\cal L}_0\otimes
 det({\bf N}_s)^{-1}}\otimes {\bf S}^{m-3\over 2}({\bf N}_s\oplus {\bf C}_B)), 
{\cal L}).$$

 Suppose that $m\leq -3$, then 
$$FSW_B(1, {\cal L}+mE)=\sum_{i\geq 0} FSW_B(c_i(\sqrt{{\cal L}_0\otimes
 det({\bf N}_s)^{-1}}\otimes {\bf S}^{-m-3\over 2}({\bf N}_s^{\ast}\oplus 
{\bf C}_B)),{\cal L}).$$

\end{theo}

\medskip
 
 The proof is almost identical to theorem \ref{theo; blup}. We omit it. $\Box$ 

The family blowup formula has a very interesting dependence on
 ${\cal L}$ through ${\cal L}_0$, unlike
 the usual blowup formula of the Seiberg-Witten theory. 
Only in the special cases where the ${\cal L}_0$ become trivial, e.g.
\ we have a trivial
 constant cross section in the trivial product family, does the dependence go 
away. The formula also depends on the complex rank two normal bundle 
${\bf N}_s$ explicitly. 
The gluing technique localizes the effect of the blowing up to the family
 invariant into a vicinity of the cross section $s:B\mapsto {\cal X}$.  In a sense 
the blowup formula should be viewed as a localization theorem of the
 Family Seiberg-Witten invariants.

We will spend some time in the next subsection \ref{subsection; apply}
 to find out the geometric meaning of the dependence.  
The formula also have a nontrivial
 dependence on the odd integer $m$, while the original blowup formula
 has no explicit dependence in $m$ at all. The reader may notice that there is a  
 duality symmetry between $m\mapsto -m$. In the proof we see that 
 the ${\bf Z}_2$ symmetry roots at the Serre duality of ${\bf CP}^2$.

\medskip

\subsection{\bf Applications to Counting Singular Curves} \label{subsection; apply}

\medskip

In this subsection, we discuss the relation between the
family blowup formula and the counting singular curves with prescribed
 multiplicities. 
The purpose is to link up the cohomological
 information of the family blowup formula derived in the previous section 
with the algebraic-geometric data which also appears in the discussion upon 
ideal sheaves of points.

We would like to achieve a few goals in this subsection:

\medskip

\noindent I. Understand the algebraic structure appearing in the family
 blowup formula. 
 Relate it with a pseudo-holomorphic or algebraic geometric question of
 counting curves with
 singularities of prescribed multiplicities.

\medskip

\noindent II. Motivate the definition of algebraic Seiberg-Witten invariant
 (see section 
\ref{section; define}) and
 an algebraic proof of family blowup formula. (see section \ref{section; ap})

\medskip

\noindent III. In proposition \ref{prop; vanishing}, use an example to
 illustrate why the 
$SW$ simple type condition
 for $b^+_2>1$ symplectic four-manifolds implies the vanishing of all the 
family invariants used to
count singular curves. Again, motivate to define algebraic Seiberg-Witten 
invariant without the
 simple type property.

\medskip
 
\noindent IV. Give some simple example of applying family blowup formula 
to count singular curves, as a
 warm up to the main application in [Liu1].

\medskip

The blowing up construction 
does not require the fiber bundle ${\cal X}\mapsto B$ to carry fiber-wise 
almost complex structures. If it does, the tubular neighborhood of 
$s:B\mapsto {\cal X}$ inherits
 the almost complex structures from the ambient space. As there are abundant 
examples of fiber bundles with fiberwise almost complex structures,
 there are countless examples for which the family blowup formula can 
be applied to. 
  
 We also recall that when a four-manifold $M$ carries an almost complex
 structure,
 any $spin^c$ determinant line bundle ${\cal L}$ can be re-written as 
$2C+K_M^{-1}$ in the
additive notation. The class $C$ is the cohomology class appearing in 
the
 Gromov-Taubes theory such that $SW({\cal L})=Gr(C)$.
 Considers a family which carries
 fiber-wise almost complex structures, then $det({\bf N}_s)$
 is isomorphic to ${\bf K}^{-1}_{{\cal X}/B}|_{s(B)}$. 
 Then the expression $\sqrt{{\cal L}_0\otimes det({\bf N}_s)^{-1}}$ is
 nothing but
 the restriction of $\sqrt{ {\cal L}\otimes {\bf K}_{{\cal X}/B}}$ to 
the section $s:B\mapsto {\cal X}$.
  This is the first hint that the pure topological
 derivation of the family blowup formula may have something to do with 
the Gromov-Taubes
  theory. 
 
 The dependence on ${\bf N}_s$ is manifestly present
 at the term ${\bf S}^{m-3\over 2}({\bf N}_s\oplus {\bf C}_B)=
\oplus_{i\leq {m-3\over 2}}{\bf S}^i({\bf N}_s)$. Given the complex normal bundle
 ${\bf N}_s$, it is exactly the bundle of polynomial algebra in ${\bf N}_s$
 of degree less than or equal to ${m-3\over 2}$.
 
To understand
 the structure, let us consider a very special case of the family blowup formula.
 
\begin{mem} \label{mem; com-sur}

 Let $M$ be a complex surface. Consider ${\cal X}=M\times M\mapsto M$
 to be the product
fiber bundle. Instead of using the trivial constant cross section, we consider
 the diagonal cross section $\Delta: M \mapsto M\times M$. 
 It is well known that the normal bundle of $\Delta(M)\subset M\times M$
 is isomorphic to the tangent bundle ${\bf T}_M$ itself. 
 Blowing up $\Delta(M)\subset M\times M$ in the complex category produces 
${\cal X}'=Blowup_{\Delta(M)}M\times M$.
 ${\cal X}'\mapsto M$ is the blown up fibration from $M\times M\mapsto M$.

Let ${\bf E}_C$ be a complex line bundle
 on $M$ with $c_1({\bf E})=C\in H^2(M, {\bf Z})$. The pull-back of 
${\bf E}_C$ to ${\cal X}=M\times M$
 induces a line bundle on the fiber bundle $M\times M\mapsto M$. 
 The restriction of any complex line bundle ${\bf E}_C$ to the diagonal 
section $\Delta
: M\mapsto M\times M$, is isomorphic to ${\bf E}_C$.  In the example take 
${\cal L}={\bf E}^2_C\otimes {\bf K}^{-1}_M$. As usual let $E$ denote the 
exceptional class of the blowing up.

 According to family blowup formula for $m=-2p-1<-1$, 

 $$FSW_M(1, {\cal L}-(2p+1)E)=
\sum_i FSW_B(c_i({\bf E}_C\otimes \bigl(\oplus_{0\leq k\leq p-1}
{\bf S}^{k}({\bf T}^{\ast}_M)\bigr)), {\cal L}).$$

 The obstruction
 bundle whose Chern
 classes are inserted in the blowup formula is given by 
 $${\bf E}_C\otimes(\oplus_{k\leq p-1}{\bf S}^i({\bf T}^{\ast}_M)).$$

 Taking a closer look at the bundle, one realizes that it is nothing
 but the $p-1$jets bundle of the line bundle ${\bf E}_C$. 
Take an arbitrary point
 $x \in M$. Consider an arbitrary germ of smooth section ${\bf s}$ of ${\bf E}_C$. 
 Recall that given a connection on a smooth line bundle ${\bf E}_C$, the covariant
 derivative defines a map $$\nabla:\Gamma(U, {\bf E}_C)\otimes {\bf T}_M
\mapsto \Gamma(U, {\bf E}_C),$$ for open sets $U$ containing $x$.

 We consider $({\bf s}(x),(\nabla {\bf s})(x), 
{\bf S}^2(\nabla\nabla){\bf s} (x),\cdots,
{\bf S}^{p-1}(\nabla\cdots\nabla){\bf s} (x))$, where ${\bf S}^k(\nabla\cdots \nabla)$ denotes 
the symmetrized $k-$th order covariant derivative operator. 
When we let the point $x$ move
 along the base manifold $M$, we find that the datum exactly have the right 
dependence as a section of 
${\bf E}_C\otimes(\oplus_{i\leq p-1}{\bf S}^i({\bf T}^{\ast}_M))$.

\end{mem}

 \medskip 

 Now we are ready to perform 
the calculation.
 From the family blowup formula (see page \pageref{theo; blup}), 
the pure family invariant on the blown up fiber bundle
 is expressed as the sum
 three terms of mixed family invariants over $M\times M\mapsto M$. 
The first term is the pure invariant.
 The second term is the mixed invariant with 
$c_1({\bf E}_C\otimes {\bf S}^{p-1}({\bf T}^{\ast}_M \oplus {\bf C}_M))$ inserted. 
The third term is the mixed family
 invariant with 
$c_2({\bf E}_C\otimes
{\bf S}^{p-1}({\bf T}^{\ast}_M\oplus {\bf C}_M))$ 
inserted, calculated with respect to the $spin^c$ structure ${\cal L}=
{\bf E}_C^2\otimes {\bf K}^{-1}_M$
 on $M\times M\mapsto M$.

Because the fiber bundle ${\cal X}=M\times M\mapsto M=B$
 is a trivial product bundle,
it is easy to see that the first two terms 
 always vanish.  It is because the family moduli space of 
${\cal L}={\bf E}_C^2\otimes {\bf K}_M^{_1}$
 over $B=M$ will be a trivial
 product. It implies the vanishing of the pure and mixed invariants unless the 
base class cohomological insertion has a degree equal to the base dimension. 
 Thus the pure invariant of 
the blown up fiber bundle $FSW_M(1, {\cal L}-(2p+1)E)$ is equal 
to the product of usual $SW$ invariant of
${\cal L}$ (over $B=$a point) and 
$\int_{M}c_2({\bf E}_C\otimes {\bf S}^{p-1}({\bf T}^{\ast}_M \oplus {\bf C}_M))$.

 Suppose we choose $p>0$, it is not hard to see that 
$${c_1({\cal  L}+mE)^2-2\chi(M)-3\sigma(M)+1\over 4}+4<{
c_1({\cal L})^2-2\chi-3\sigma
\over 4}.$$
 In order that $FSW_M(1, {\cal L}-(2p+1)E)$ to be nonzero,
 the moduli space dimension of ${\cal L}$ over $B=pt$, must be strictly
 positive.

Suppose that $M$ is a symplectic manifold with $b^+_2>1$, then
according to Taubes' result [T2], the
 manifold is of Seiberg-Witten
 simple type. That is to say that all the moduli spaces of basic 
classes (whose $SW\not=0$) are of 
 expected dimension zero (see [FM] for a derivation for Kahler surfaces). 
This observation gives us the following vanishing 
result,

\begin{prop}\label{prop; vanishing}
 Let $m$ be an odd integer whose absolute value is 
bigger than one. Let ${\cal L}$ be an arbitrary $spin^c$ structure on a 
Kahler surface with $p_g>0$ (or  more generally any
 symplectic four manifold with $b^+_2>1$), then the pure Seiberg-Witten
 invariants of ${\cal L}+mE$, $FSW_B(1, {\cal L}+mE)$
 on $M_2=Blowup_{\Delta(M)}(M\times M)\mapsto M$ vanish.
\end{prop}

 The space $M_2$ is the $l=2$ version of the universal space $M_l$ (see [Liu1]).

\noindent Proof of the proposition: From 

 $$FSW_B(1, {\cal L}-(2p+1)E)=\int_M c_2( {\bf E}_C\otimes 
{\bf S}^{p-1}({\bf T}^{\ast}_M \oplus {\bf C}_M))
\cdot SW({\cal L}),$$

 in order that $FSW_B(1, {\cal L}-(2p+1)E)\not=0$, $SW({\cal L})$ 
must be non-zero. But this implies that
 ${\cal L}$ is a basic class. If the surface has $p_g>0$, 
$c_1({\cal L})^2-2\chi-3\sigma=0$. This violates the
 bound on the expected dimension we get above.  The argument for 
$m>1$ case is almost identical to the $m=-(2p+1)<-1$ case.
 If $M$ is a $b^+_2>1$ symplectic four-manifold, one may replace 
complex blowing up by almost complex blowing up and
 the same argument works, too. $\Box$

Let us explain the relevance of this piece of calculation with 
symplectic
 geometry. If we consider ${\cal L}={\bf E}^2_C\otimes {\bf K}_M^{-1}$, 
then the Seiberg-Witten
 invariant of ${\cal L}$ is equal to the enumeration of pseudo holomorphic curves
Poincare dual to the cohomology class $C=c_1({\bf E}_C)$. When we consider the 
 singular curves dual to $C$ with a singularity of the prescribed multiplicity $p$,
 they can be resolved into smooth curves on the blown up manifold dual to 
$C-pE$ . Or in term of $spin^c$ structure, $2({\bf E}-pE)+{\bf K}_M^{-1}-E
={\cal L}-(2p+1)E$. When the singular point is allowed to move on the whole 
$B=M$, the calculation of
 $FSW_M(1, {\cal L}-(2p+1)E)$ should be directly related to the counting of 
curves dual to $C$ with a multiplicity
 $p$ singularity in $M$.
 Thus the vanishing result prop. \ref{prop; vanishing} on $FSW_M(1, 
{\cal L}+mE)$ indicates that such
 counting of the pseudo-holomorphic singular curves always gives the
answer $0$ when the corresponding symplectic four-manifold is of 
Seiberg-Witten simple type.   

 On the other hand, if the symplectic four-manifold $M$ is of $b^+_2=1$, then
 it is not of non-simple type in the Taubes chamber (determined by large 
deformations of its symplectic forms)
 and the usual $SW$ invariants of ${\cal L}$ for ${c_1^2({\cal L})-
2\chi-3\sigma\over 4}>0$ is calculated by 
 the so-called wall crossing formula [KM], [LL2] and is $\pm 1$ for $b_1=0$
 manifolds. 
 Then the family blowup formula predicts that the pure invariants on $M_2\mapsto
 M$ are given by 
$\pm \int_M c_2({\bf E}_C\otimes {\bf S}^{p-1}(
{\bf T}^{\ast}_M\oplus {\bf C}_M))$, the same
 as gotten by counting of singular curves by algebraic geometers [V].

 In the case when $b^+_2=3$ 
and the manifolds carry hyper-winding families of symplectic forms, 
one thickens the base by multiplying with a copy of $S^{b^+_2-1}$.

\medskip

\subsection{\bf The Family Blowup Formula and the Universal Family}
\label{subsection; uf}

\medskip

 One important application of the family blowup 
formula is in the long paper [Liu1], concerning
 the enumeration of singular curves with nodal or
 other singularities. Let us 
 give a slightly slow paced discussion about its 
relationship with family blowup formula. 
 This subsection is an extension of the example \ref{mem; com-sur}.

 Let $M$ be a symplectic four manifold with a compatible 
almost complex structure $J:TM\mapsto TM$.
 We fix such an almost complex structure and view $M$ as 
an almost complex manifold.
 Recall that the sequence of universal spaces $M_k, k\in {\bf N}$ and 
 $f_k:M_{k+1}\mapsto M_k$ can be constructed inductively
 by the recipe in [Liu1]. Set $M_0=pt$ and $M_1=M$. Then 
define $f_0:M_1\mapsto M_0$ to be the constant map.
 Suppose that $M_l, l\leq k, l\in {\bf N}$ and $f_{l-1}:M_l\mapsto 
M_{l-1}$ have been defined for all $l\leq k$ such that
 $f_i$ are smooth pseudo-holomorphic submersions.
 We define $M_{k+1}$ and $f_k: M_{k+1}\mapsto M_k$ by the following recipe.

  Take the fiber product of $M_k\times_{M_{k-1}} M_k$ 
through $f_{k-1}:M_k\mapsto M_{k-1}$ and 
  $f_{k-1}:M_k\mapsto M_{k-1}$.  Then $f_{k-1}:M_k\mapsto M_{k-1}$ 
maps relatively into the diagonal of
 $M_k\times_{M_{k-1}}M_k$ as an almost complex manifold. 
Consider the almost complex blowing up of 
 $\Delta_{M_{k-1}}:M_k\mapsto M_k\times_{M_{k-1}}M_k$ as a 
complex codimension two sub-manifold.
 Define $M_{k+1}$ to be the blown up manifold and it maps 
naturally to either copy of $M_k$ surjectively.
 
  By almost complex blowing up, we mean the following: The normal 
bundle of the relative diagonal is well known to 
 be isomorphic to the relative tangent bundle of $M_k\mapsto M_{k-1}$.
 Because $f_{k-1}$ is pseudo-holomorphic,
  the relative tangent bundle ${\bf T}_{M_k/M_{k-1}}$, as the kernel 
of $df_{k-1}$ is stable under the almost complex
 structure of $M_{k}$ and becomes a complex rank two vector bundle 
over $M_k$. Through the isomorphism, the normal
 bundle of the relative diagonal ${\bf N}_{\Delta(M_k)}M_k
\times_{M_{k-1}}M_k$ inherits the structure of a complex
 vector bundle. To construct the almost complex blowing up of 
$\Delta(M_k)\subset M_k\times_{M_{k-1}}M_k$ replaces 
 $\Delta(M_k)$ by ${\bf P}({\bf T}_{M_k/M_{k-1}})$. The almost 
complex structure outside of the blown-up locus 
 $ M_k\times_{M_{k-1}}M_k-\Delta_{M_{k-1}}(M_k)$ is unchanged. 
The almost complex structure on 
 ${\bf T}_{M_{k+1}}|_{{\bf P}({\bf T}_{M_k/M_{k-1}})}$ is induced
 from the natural almost complex structure from
 ${\bf T}_{{\bf P}({\bf T}_{M_k/M_{k-1}})}$ and 
${\bf N}_{\Delta(M_k)}M_k\times_{M_{k-1}}M_k\cong 
{\bf T}_{M_k/M_{k-1}}$. It is easy to check that the almost 
complex
structure induced in this way is ${\cal C}^{\infty}$ on the whole 
$M_{k+1}$. Moreover the natural map
 $f_k: M_{k+1}\mapsto M_k$, defined as the composition of 
$M_{k+1}\mapsto M_k\times_{M_{k-1}}M_k$ and
 $M_k\times_{M_{k-1}}M_k\mapsto M_k$ (either copy) 
is a composition of pseudo-holomorphic maps and is therefore
 pseudo-holomorphic.

\medskip

\begin{lemm}\label{lemm; consecutive}
 The fiber bundle $f_k:{\cal X}_{k+1}=M_{k+1}\mapsto M_k=B$ 
can be constructed from the product fiber bundle 
${\cal X}_0=M\times M_k\mapsto M_k$ by $k$ consecutive blowing
 ups ${\cal X}_i, 1\leq i\leq k$.

\noindent (i). ${\cal X}_i\mapsto M_k=B$ is constructed from 
${\cal X}_{i-1}\mapsto M_k=B$
 by blowing up a cross section $M_k\mapsto {\cal X}_{i-1}$.

\noindent (ii). Consider $f_{k-1, i}:
M_k\mapsto M_i, i\leq k-1$ to be $f_{k-1, i}=
f_i\circ f_{i+1}\circ \cdots f_{k-1}$ and
 $f_i: M_{i+1}\mapsto M_i$,  
then  ${\cal X}_i$ can be identified with the fiber product 
$M_k\times_{M_i}M_{i+1}$ of $f_{k-1, i}$ with $f_i$, the
 pull-back of $f_i: M_{i+1}\mapsto M_i$ by $f_{k-1, i}:M_k\mapsto M_i$.
\end{lemm}

\medskip

\noindent Proof of the lemma: The proof is essentially the same as the 
proof of  lemma 3.1. in [Liu1]. $\Box$

 Knowing that fiber bundle
 ${\cal X}_k$ and the original fiber bundle ${\cal X}_0=M\times M_k$
 are related by $k$ different blowing ups, the cohomology class 
$C-\sum_{i\leq k} m_iE_i$ on ${\cal X}_k$ 
 is related to $C$ on ${\cal X}_0$ through $C-m_1E_i$, $C-m_1E_1-m_2E_2$, 
$\cdots$. To simplify the notation, 
 we have identified $C$ on ${\cal X}_0$ with its 
pull-back on the blown-up manifolds ${\cal X}_i, i\leq k$ and denote them
 by the same symbol $C$.

\begin{prop}\label{prop; blowdown}
  For $C^2-C\cdot c_1(K_M)-\sum_{i\leq k}(m_i^2-m_i)\geq 0$,
the pure family invariant on ${\cal X}_k\mapsto M_k$, $FSW_{M_k}(1, 
{\bf E}_{2C-\sum_{i\leq k}(2m_i+1)E_i}\otimes 
{\bf K}^{-1}_M)$ is equal to

 $$\hskip -.9in SW({\bf E}_{2C}\otimes {\bf K}^{-1}_M)\cdot \int_{M_k}c_{2k}\bigl(
 {\bf E}_C\otimes {\bf S}^{m_1-1}(
f^{\ast}_{k-1, 1}{\bf T}^{\ast}_M\oplus {\bf C}_{M_k})\oplus 
{\bf E}_{C-m_1E_1}\otimes 
{\bf S}^{m_2-1}( f_{k-1, 2}^{\ast}
{\bf T}^{\ast}_{M_2/f_1^{\ast}M_1}\oplus {\bf C}_{M_k})\oplus $$
$$\cdots
 {\bf E}_{C-\sum_{i\leq k-1}m_iE_i}\otimes {\bf S}^{m_k-1}( f_{k-1, k}^{\ast}
{\bf T}^{\ast}_{M_k/f_{k-1}^{\ast}M_{k-1}}
 \oplus {\bf C}_{M_k})\bigr).$$

\end{prop}

\medskip

\noindent Proof of the proposition:  By using family blowup
 formula consecutively, one may relate the
 pure family invariant 
$FSW_{M_k}(1, {\bf E}_{2C-\sum_{i\leq k}(2m_i+1)E_i}\otimes 
{\bf K}^{-1}_M)$ on ${\cal X}_k\mapsto M_k$
 to a combination of 
the mixed family invariants $FSW_{M_k}(\eta, {\bf E}_{2C}\otimes 
{\bf K}^{-1}_M)$ over $M\times M_k\mapsto M_k$, 
 $\eta\in H^{\ast}(M_k, {\bf Z})$.  Because the map 
$f_{k-1, i}$ factors through $f_i$, then
 ${\cal X}_i$ can be viewed as the pull-back of 
$M_{i+1}\times_{M_i}M_{i+1}$ by $f_{k-1, i+1}: M_k\mapsto M_{i+1}$.

 Thus, the relative diagonal of $M_{i+1}\times_{M_i}M_{i+1}$ 
pulls back to a cross section ${\bf s}_i$ of ${\cal X}_i\mapsto M_k$.
On the other hand, pull-back by $M_k\mapsto M_{i+1}$ of the blowing
 up of the $\Delta_{M_i}M_{i+1}\subset M_{i+1}\times_{M_i}M_{i+1}$
 is nothing but ${\cal X}_{i+1}$. Thus, ${\cal X}_{i+1}$ can be
 thought as constructed from ${\cal X}_i$ by the pull-back of the
 relative diagonal section $\Delta_{M_i}M_{i+1}\subset M_{i+1}\times_{M_i}M_{i+1}$.

 Therefore, one may identify the complex rank two normal bundle 
${\bf N}_{{\bf s}_i(M_k)}{\cal X}_i$ to be 
$f^{\ast}_{k-1, i+1}{\bf T}_{M_{i+1}/f_i^{\ast}M_i}$.

The class $\eta$ should be a combination of cup products of various
 Chern classes of 
 the bundles ${\bf E}_{C-\sum_{j\leq i-1}m_jE_j}\otimes
 {\bf S}^{m_i-1}(f^{\ast}_{k-1, i}{\bf T}_{M_i/f_i^{\ast}M_{i-1}}\oplus 
 {\bf C}_{M_k})$.

  By using the product rule of the total Chern
 classes under bundle addition, $\eta$ can be identified with 
  $$c_{total}\bigl(\oplus_{i\leq k}{\bf E}_{C-\sum_{j\leq i-1}m_jE_j}\otimes 
{\bf S}^{m_i-1}(f^{\ast}_{k-1, i}{\bf T}_{M_i/f_i^{\ast}M_{i-1}}\oplus 
 {\bf C}_{M_k})\bigr).$$

  However, because ${\cal X}_0=M\times M_k\mapsto M_k=B$ is a 
product fiber bundle,
only the grade $4k$ component of $\eta$ contributes to the
 mixed invariant. The pure family invariant is equal to

$$\int_{M_k}c_{2k}\bigl(\oplus_{i\leq k}{\bf E}_{C-\sum_{j\leq i-1}m_jE_j}\otimes 
{\bf S}^{m_i-1}(f^{\ast}_{k-1, i}{\bf T}_{M_i/f_i^{\ast}M_{i-1}}\oplus 
 {\bf C}_{M_k})\bigr)\cdot FSW_{M_k}([M_k], {\bf E}_{2C}\otimes {\bf K}^{-1}_M)$$

$$= \int_{M_k}c_{2k}\bigl(\oplus_{i\leq k}{\bf E}_{C-
\sum_{j\leq i-1}m_jE_j}\otimes 
{\bf S}^{m_i-1}(f^{\ast}_{k-1, i}{\bf T}_{M_i/f_i^{\ast}M_{i-1}}\oplus 
 {\bf C}_{M_k})\bigr)\cdot SW({\bf E}_{2C}\otimes {\bf K}^{-1}_M).$$

 This ends the proof of the proposition. $\Box$

\bigskip

\begin{rem}\label{rem; poly}
  One may calculate the integral of the $2k-th$ Chern
 class $c_{2k}$ over $M_k$ by pushing forward along
 $f_i: M_{i+1}\mapsto M_i$ consecutively. By using the 
blowup formula of the Chern classes of tangent bundles,
 the final answer depends on $C=c_1({\bf E}_C), c_1({\bf T}_M),
 c_2({\bf T}_M)$ and all the multiplicities
 $m_1, m_2, \cdots, m_k$. It has to be a universal
 (manifold independent) polynomial of 
 $C^2[M], C\cdot c_1(M)[M], c_1^2(M)[M], c_2(M)[M]$, 
$P_{m_1, m_2, \cdots, m_k}(C^2[M], C\cdot c_1(M)[M], c_1^2(M)[M], c_2(M)[M])$. 

\end{rem}

\medskip

\subsection{\bf The Existence of Pseudo-Holomorphic Singular
 Curves with Prescribed Singular Multiplicities}
\label{subsection; mul}

\bigskip

Given a symplectic four manifold $M$ and a cohomology 
class $C$ whose Gromov moduli space has non-negative dimension. 
If $FSW=FGT$  is known to hold for those general families, then the
number of pseudo-holomorphic curves singular (with restriction on the orders
 of the singularities) at a finite number of 
points are symplectic invariants and the invariants are determined by the
 Gromov invariants of the class and the Chern classes information of the
 tangent bundles, etc. 
If  $M$ is a symplectic four manifold with $b^+_2>1$,
 Taubes [T1,T2, T3] proves that $M$  
has simple type. It means that once the
 moduli space dimension is positive, its invariant vanishes. On the other
 hand the dimension of the 
 moduli spaces of singular curves are strictly smaller than that of
the original moduli spaces. 
Hence it reasonable to speculate that
 the singular curves invariants  actually all vanish. 

 At first the speculation may look incompatible with the intuition we have.
 On every algebraic surface $X$ we are able to count e.g.\ nodal curves on 
$X$. The enumerative problem is well known in algebraic geometry. It does
 gives nontrivial answers in general. However our theorem tells us that
 the counting problem in algebraic geometric way does not give us the 
symplectic invariants when $b^+_2>1$. Only when $b^+_2=1$, the algebraic
calculation and the symplectic calculation coincide completely.  On the
 other hand, it does not mean that the symplectic counting problem
 does not
make sense for $b^+_2>1$ symplectic four manifolds. If fact, if we can
construct some family of symplectic four manifolds which become non-simple
type, then we can still make sense of the symplectic singular curves 
counting and get a nonzero answer. A good example is K3. As $K3$ is of 
$b^+_2=3$, it falls into the category discussed in the corollary. If we
count the singular invariants in the usual way, we always get zero. 
 As K3 has hyperkahler structures, a generic
 K3 has no integral $(1,1)$ class, not mentioning any effective holomorphic
curves. However we can still consider the $S^2$-family and the
 corresponding family invariants.
 In the final section, we make a systematic study upon the algebraic
 family Seiberg-Witten ``invariants''.  It will be indicated that the
 algebraic ``invariants'' and the smooth invariants coincide only when
 $b^+_2=1$.  Otherwise the algebraic ``invariants'' correspond to a
 ``local'' $p_g$ dimensional smooth family invariant.  The existence of
 the local family invariants was previously
  speculated independently by T.J. Li and the author. This also gives
 a philosophical explanation why the algebraic geometers can count the
 curves when the symplectic geometers claim the triviality of the invariants.

From the previous discussion we learn that in the Taubes' chamber, the
 invariants are all $\pm 1$ if $C^2\geq -2$. Using the previous formula
 we learn that if we count the singular curves in the $S^2$ family, we
get nontrivial results.  To demonstrate the power of the blowup formula, let
 us discuss some other explicit examples.

 Let us discuss another interesting applications of the family blowup
 formulas to the special families considered previously. Let $M$ be a symplectic 
four-manifold
 with $b^+_2=1$.
Take $B$ to be $M_k$, the $k-$th 
universal space (see [Liu1]). Take ${\cal X}_0=M\times M_k\mapsto M_k=B$ and
 ${\cal X}_k=M_{k+1}\mapsto M_k$ which can be constructed from ${\cal X}_0$ by
 $k$ consecutive blowing ups of cross sections. Given a class 
$C\in H^2(M, {\bf Z})$, it determines
 a smooth line bundle ${\bf E}_C$ on $M$. 
 Consider the $spin^c$ determinant line bundle $2{\bf E}^2-
\sum_{i\leq k}(2m_i+1)E_i-{\bf K}_M$ (in 
additive notation).
 
 Suppose one chooses the class $\eta=[M_k]\in H^{top}(M_k,{\bf Z})$
 and considers the corresponding blowup formula of mixed family
 invariants expressing
 $FSW_{M_k}(\eta, 2{\bf E}^2_C-\sum_{i\leq k}(2m_i+1)E_i-{\bf K}_M)$ in terms of
the mixed invariants from $M\times M_k\mapsto M_k$.
   
 One concludes the following
 interesting corollary,

\begin{cor} \label{cor; sym}

 Let $M$ be a symplectic four-manifold with $b^+_2=1$.
 Let $C$ be a cohomology class $\in H^2(M, {\bf Z})$ whose Gromov-Taubes 
invariant $Gr(C)$ is nonzero (see [T1], [T2]).
 Given any tuple of positive integers $m_i, 1\leq i\leq k$ with
 ${C^2-C\cdot c_1(K_M)\over 2}-{\sum_i (m_i^2-m_i)\over
 2}\geq 0$ and $k$ distinct points on $M$,
 then there exists a pseudo-holomorphic curve Poincare dual to $C$
 passing through these $k$ distinct points of $M$ with multiplicities at 
the $i-th$ point
 not less than $m_i$.
\end{cor}

\medskip

\noindent Proof: Continue the discussion before the statement of the corollary.
 As $\eta$ has exhausted the dimension of the base $M_k$,
 the extra Chern classes of the obstruction bundle will not be able
 to contribute to the invariants. 

 Thus, we conclude $FSW_B(\eta, 2{\bf E}_C^2-
\sum_{i\leq k}(2m_i+1)E_i-{\bf K}_M)$ is equal to
 $FSW_B(\eta, 2{\bf E}_C-{\bf K}_M)$ on $M\times M_k\mapsto M_k$.
On the other hand, the particular
 mixed invariant of $M\times M_k\mapsto M_k$ is nothing but the usual 
$SW(2{\bf E}_C-{\bf K}_M)$, which
 according to Taubes theorem 'SW=Gr', is equal to $Gr(C)$. 

 Thus one concludes that $FSW_B(\eta, 2{\bf E}^2-
\sum_{i\leq k}(2m_i+1)E_i-{\bf K}_M)$ is always
 nonzero if 

(1). The family Seiberg-Witten dimension of 
 $2{\bf E}^2_C-\sum_{i\leq k}(2m_i+1)E_i-{\bf K}_M$ is 
non-negative, which is reduced to the condition
 ${C^2-C\cdot c_1(K_M)\over 2}-{\sum_i (m_i^2-m_i)\over
 2}\geq 0$

(2). The ordinary Seiberg-Witten invariant of the class $2{\bf E}_C-{\bf K}_M$
  in Taubes' chamber $=Gr(C)$ is nonzero. 

 For any $k$ distinct points  $p_1, p_2, \cdots, p_k$ in $M$, it 
determines a point in the top (open) stratum of $M_k$.

  As the particular choice of $\eta=[M_k]$ is
 poincare dual to the zero cycle on $M_k$, the family invariant 
$FSW_{M_k}([M_k], 2{\bf E}^2_C-\sum_{i\leq k}(2m_i+1)E_i-{\bf K}_M)$ 
can be re-interpreted as the
 a counting of Seiberg-Witten solutions on the fiber above the point 
$\in M_k$. The fiber above the given point
 in $M_k$ is a symplectic blowing up $\tilde{M}$ of $M$ at these $k$ 
distinct points.
Then the family blowup formula
 asserts that if we deform the family Seiberg-Witten equations by a 
large deformation of a given family of
  symplectic forms, there is always some solutions of the Seiberg-Witten equations
 above the given point in $M_k$. Otherwise, the invariant count
 is zero for an empty moduli space.
 In particular, the same conclusion still holds if one takes
 the $r\mapsto \infty$ limit
 in the family Seiberg-Witten equations.

 On the other hand, by Taubes' analysis the solution of 
Seiberg-Witten equation of large symplectic form 
perturbation will converges (in the $r\mapsto \infty$ limit) to a
 $(1, 1)$ current which can be 
 regularized to be a pseudo-holomorphic curve dual to 
$C-\sum m_i E_i$. The readers can consult [T1, T2, T3]
 for details.

   Finally, a pseudo-holomorphic curve dual to 
$C-\sum m_i E_i$ in the blown up manifold $\tilde M$ gives rise
 to a pseudo-holomorphic curve dual to $C$ in $M$. The curve
 has to develop a multiplicity $m_i$ singularity
 at $p_i, i\leq k$ as the intersection number of its 
'proper transformation'' in $\tilde M$ has intersection 
 multiplicity $m_i$ with the exceptional curve $E_i$ above $p_i$. $\Box$

  The corollary tells us that one can construct a singular
 divisor(curve) on $M$ with
 prescribed singular multiplicities and singular points.  
Notice that the same type of
 conclusion was known for algebraic surfaces and the proof relies
 on the sheaf sequences for the ideal sheaves and the vanishing result
 for the higher sheaf cohomologies. 

 Let $M$ be an algebraic surface with $p_g=0, q(M)=0$ 
 and let ${\cal I}_Z$ be the ideal sheaf of $Z=\sum_{i\leq k} m_i p_i, p_i\in M$,
 $m_i\in {\bf N}$, then one has the following sheaf short exact sequence 

 $$0\mapsto {\cal I}_Z\mapsto {\cal O}_M\mapsto {\cal O}_Z\mapsto 0.$$

 For a locally free ${\cal E}$ with $c_1({\cal E}=C$, 

$$0\mapsto {\cal I}_Z\otimes {\cal E}\mapsto 
{\cal E}\mapsto {\cal O}_Z\otimes {\cal E}\mapsto 0$$

 induces the long exact sequence 

$$0\mapsto H^0(M, {\cal I}_Z\otimes {\cal E})\mapsto
 H^0(M, {\cal E})\mapsto H^0(Z, {\cal O}_Z\otimes {\cal E})
 \mapsto H^1(M, {\cal I}_Z \otimes {\cal E})\cdots.$$ 

 If ${\cal E}$ is sufficiently very ample to imply the
 vanishing of $h^1(M, {\cal I}_Z \otimes {\cal E})$, the 
 long exact sequence truncates to a short exact sequence.
 A curve in the linear system 
 ${\bf P}(H^0(M, {\cal E}))$ has the prescribed singular
 multiplicities on $p_i, i\leq k$ if and only if 
 the restriction of the corresponding defining section to 
the non-reduced $Z$ vanishes.  

 Given that $h^0(Z, {\cal O}_Z\otimes {\cal E})=\sum_i {m_i(m_i-1)\over 2}$ 
and
 $h^0(M, {\cal E})\geq \chi(M, {\cal E})=1+{C^2-c_1(K_M)\cdot C\over 2}$
 (assuming $h^2(M, {\cal E})=0$), 
 $h^0(M, {\cal I}_Z\otimes {\cal E})\geq \chi(M,
 {\cal E})-\sum_{i\leq k}{m_i^2-m_i\over 2}\geq 1$ if
 ${C^2-c_1(K_M)\cdot C\over 2}-\sum_i {m_i(m_i-1)\over 2}\geq 0$.

By applying the family blowup
 formula and Taubes' technique in symplectic geometry, corollary 
\ref{cor; sym} can be viewed as 
 the symplectic version of the theorem in $b^+_2=1$ category.   
As the corresponding
 technique in algebraic geometry has play an essential role in the
 study of linear system, one believes that the present symplectic
 version should also plays a similar role.  
The parallelism between the symplectic and the algebraic argument
 also suggests that the family blowup formula 
should have its algebraic geometric origin. This also motivates
 the definition of algebraic $SW$ invariant 
in section \ref{section; define} and the algebraic proof of family
 blowup formula in section \ref{section; ap}.

\begin{mem}\label{mem; cp}
 Let $M$ be ${\bf CP}^2$. Let the 
cohomology class $C$
 denote $dH, d\geq 1$. 
Then the expected dimension of $SW$
 (or Gromov-Taubes) moduli spaces of $(2d+3)H$ ($dH$) are of
 ${d^2+3d\over 2}$ dimension. 
Let us consider the following enumeration problem. We are
 interested in counting the curves
 which have one nodal(ordinary double)point in ${\bf CP}^2$. It is well
 known that a nodal condition decreases the moduli spaces by real two dimension.
 Therefore the nodal moduli spaces(which consist of nodal curves or their
compactifications) are real $d^2+3d-2$ dimensional. As $d^2+3d-2=(d+1)(d+2)-4$
 is always even, we can consider $(d+1)(d+2)/2-2$ generic points on
 ${\bf CP}^2$
and require the nodal curves to pass through these generic points. After imposing
these point passing conditions, generically (with respect to the almost 
complex structures
 and the points) there are finite number of nodal curves passing through
 these points.
 
 To calculate the number, we notice that by the family blowup formula we find
 the answer to be $\int_{{\bf CP}^2}c_2({\bf H}^d\otimes({\bf T}^{\ast}_{\bf
 CP^2}\oplus {\bf C}_{\bf CP^2}))\cdot
SW({\bf H}^{2d+3})$.  As we know that $SW({\bf H}^{2d+3})=1$ by the wall 
crossing formula, we
 conclude that the invariants are calculated by $c_2({\bf H}^d
\otimes({\bf T}^{\ast}_{\bf CP^2}\oplus 
{\bf C}_{\bf CP^2}))[{\bf CP}^2]$. 
After some simple calculation one gets $3(d-1)^2$.
 This answer is well known to algebraic geometers and is
 a very special case of Severi degrees.
 
\end{mem}

  Next let us discuss another simple but interesting example.
 
\begin{mem} \label{mem; k3}
 Let $M$ be a minimal symplectic four manifold with $b^+_2=3$. 
Suppose that there exists an ${\bf S}^2$ family of
 symplectic forms $\omega_{x}, x\in S^2$ on $M$ such that the projection to
 $H^2_+(M, {\bf R})-\{0\}\cong S^2\times {\bf R}^{\ast}$ is of mapping degree
 $1$. The family of symplectic 
 forms is called a hyper-twisting family of symplectic forms on $M$. 

 Such manifolds must have $c_1(M)=0$. $K3$ and $T^4$ are the only 
known examples in the Kahler category
 which have hyperwinding families of symplectic forms. These $S^2$
 families of forms can be constructed from
 hyperkahler families of Kahler forms.

 Let us consider the $S^2$ hyperwinding family of symplectic forms 
of $M$. Consider
 a primitive cohomology class $C\in H^2(M, {\bf Z})$ with square zero. 
There are an infinite number
 of these classes on such $M$. Firstly, consider the fiber bundle
 $M\times S^2\mapsto S^2$ with a
hyperwinding family of symplectic forms on the fibers.
 The family wall crossing
 formula implies that the family invariant $FSW_{S^2}(1, 2C)$ over 
$M\times S^2\mapsto S^2$ 
 is equal to $\pm 1$ in the first winding
 chamber [LL1]. 

  Consider an $S^2$ family of almost complex structures compatible with an $S^2$
 family of fiberwise Riemannian metrics.
  Consider the blown up fiber bundle ${\cal X}'=S^2
\times M_2\mapsto S^2\times M_1=B$ using the $S^2$ family of almost
 complex structures. It is easy to see that one can 
use the symplectic blowup construction to construct 
 an $S^2$ family of symplectic forms on $S^2\times M_2\mapsto 
S^2\times M_1=B$.
  By a similar calculation as in example \ref{mem; com-sur}, 
the family invariant $FSW_{S^2}(1, 2C-\sum 5E)$, evaluated 
in the winding chamber, 
is equal to $$c_2({\bf E}_C\otimes
  ({\bf T}^{\ast}M \oplus {\bf  C}_{M}))
[M]\cdot FSW_{S^2}(1, 2C)=\pm c_2({\bf E}_C\otimes
  ({\bf T}^{\ast}M \oplus {\bf  C}_{M}))$$
$$=c_2(M)[M]=\chi(M).$$ We have used $C^2=0$, 
$c_1({\bf T}^{\ast}M\oplus {\bf C}_M)=0$ and
 $FSW_{S^2}(1, 2C)=\pm 1$ in the winding chamber 
over $S^2\times M\mapsto S^2=B$.

  On the other hand, the Riemann-Roch formula for 
almost complex manifolds implies that on $M$ we have
 $$2-{b_1(M)\over 2}={1-b_1(M)+b^+_2(M)\over 2}=
ind(\bar\partial)={\int_M(c_1^2+c_2)\over 24}={\chi(M)\over 24}.$$

 Thus we find that the family invariant $FSW_{S^2}(1, 2C-5E)$ 
over $S^2\times M_2\mapsto S^2\times M$ is
 $\pm 24$ for $b_1(M)=0$ manifold $M$.

The answer $24$ does not come out accidentally.
 If we consider a generic 
 elliptic $K3$ surface elliptic fibered over
${\bf CP}^1$, then it is well known there are 24 singular nodal fibers.
 
The tool of the family blowup formula implies that on a 
simply connected hyperwinding family
 of symplectic four-manifold $M$, one can recover the number of 
single node nodal rational curves 
 within the $S^2$ family as the family invariant 
$FSW_{S^2\times M_1}(1, {\bf E}_{2C-3E})$, 
 and is identical to the number of singular nodal fibers of an
 elliptic $K3$ surface.
\end{mem}

 This above picture supports the following conjecture,

\begin{conj} \label{conj; k3}
 Let $M$ be a simply connected symplectic four-manifold with an 
$S^2$ family of hyper-winding
 symplectic forms, then $M$ is diffeomorphic to the underlying 
smooth manifold of the $K3$ surface
 and the hyperwinding family of symplectic forms is homotopic to 
the $S^2$ families of hyperkahler 
structures of the $K3$ surfaces through $S^2$ families of symplectic forms.
\end{conj} 

 One can formulate a similar conjecture for $T^4$ or other primary Kodaira
 surfaces (with $b^+_2=2$). The uniqueness of the $S^2$ or $S^1$ families up to
 homotopies  plays a crucial
 role in understanding the symplectic structures of these $c_1=0$ symplectic
 four-manifolds.

\bigskip

\section{\bf The Definition of Algebraic Seiberg-Witten Invariants}\label{section;
 define}

\bigskip

In this section we discuss the definition of algebraic Seiberg-Witten
 invariant on algebraic surfaces. The proof of main theorem \ref{main; ASW}
occupies the whole section. Because of the algebraic nature of the discussion,
 we will make use of algebraic (holomorphic) vector bundles and locally free
sheaves frequently. We adopt the convention that if we use the bold character
to denote an algebraic vector bundle, e.g. ${\bf E}$, the calligraphic
 character, ${\cal E}$, will denote the locally free sheaf of sections of
 ${\bf E}$ and vice versa. 

Given an algebraic surface $M$, there are two
important holomorphic invariants associated with $M$, $q(M)$, the irregularity of
 $M$ and $p_g$, the geometric genus of $M$. They are related to the 
homological invariant $b_1(M), b^+_2(M)$ by the relationship

$$b_1(M)=q(M), b^+_2(M)=1+2p_g.$$

   Both of the invariants are essential in defining the algebraic version of
 (family) Seiberg-Witten invariants. Let us discuss briefly before
 the formal mathematical treatment. If the geometric genus of the surface is
 $0$, then the algebraic family invariant coincides with the topological
family Seiberg-Witten invariant defined in [LL1]. On the other hand,
 the algebraic Seiberg-Witten invariant differ from the usual topological
 Seiberg-Witten invariant in that it 'formally' corresponds to topological
 family Seiberg-Witten invariant of a germ of a high dimensional 
 family.  

 Given an algebraic surface $M$, the set of $spin^c$ structures on the 
underlying smooth four-manifold of $M$ is isomorphic (up to torsions) to 
 $H^2(M, {\bf Z})$, which is isomorphic to the set of isomorphism classes
 of ${\cal C}^{\infty}$ line bundles on $M$. Thus, the algebraic Seiberg-Witten
 invariant of $M$ can be viewed as a map 

$${\cal ASW}:H^{\ast}(M, {\bf Z})\mapsto {\bf Z}.$$

 For generic classes, the formal base dimensions are $p_g$. But for
some non-generic classes (defined slightly later), the formal base dimensions
 of the infinitesimal germs are in-between $0$ and $p_g$.

 Let $M$ be an algebraic surface with $p_g=0$, then algebraic and topological
 Seiberg-Witten invariants coincide.  If $q(M)=0$, then a 
${\cal C}^{\infty}$ 
topological line bundle can be given a unique holomorphic structure. On the
 other hand, the holomorphic structures of a fixed ${\cal C}^{\infty}$ 
line bundle  on a $q(M)>0$ surface have non-trivial moduli and the
 algebraic Seiberg-Witten invariant counts holomorphic curves 
 from all the different holomorphic structures of a fixed topological 
line bundle. If one would like to enumerate the holomorphic curves from 
 a particular holomorphic structure of the ${\cal C}^{\infty}$ line bundle, 
 additional cohomological insertion has to be made on the family invariant.

\medskip
\begin{rem}\label{rem; nonlinear}
 The readers with an algebraic geometric background may notice that 
the enumeration of holomorphic curves  from the zero sections of 
different holomorphic structures
 of a fixed ${\cal C}^{\infty}$ line bundle corresponds to curve counting
 in a non-linear system, while enumerations of curves from a fixed 
 holomorphic structure corresponds to curve counting in a linear system. 
\end{rem}

\medskip

  There have been different versions of Symplectic or Algebraic 
 Gromov type invariants aiming at 
 curve enumerations ([Be], [LiT1], [LiT2], [R], [RT1], [RT2], [S], [T3]).
 It may be desirable to clarify the difference
 of ${\cal ASW}$ from the usual Gromov-Witten invariants.

\medskip

 To summarize, 

I. the algebraic (family) Seiberg-Witten invariant is
 an algebraic device used to enumerate curves as the 
divisors on a given algebraic surface 
 than holomorphic maps from domain curves to the target $M$.

II. Because usual Seiberg-Witten theory has
compact moduli spaces, the algebraic (family) 
Seiberg-Witten invariants are defined using compact moduli spaces as well.

\medskip

III. As 
its definition does not involve the domain curves, the (compactification of)
Deligne-Mumford moduli spaces ${\cal M}_{g, n}$ do not come into our picture.
 Thus, the algebraic Seiberg-Witten invariants do not have nice combinatorial
 structures from the domain curves as the usual Gromov-Witten invariants do.

\medskip

IV.  On the other hand, the algebraic Seiberg-Witten invariants are related to
 Surface Riemann-Roch formula closely. It can be seen directly 
from the dimension formula
 of ${\cal ASW}$ (see the discussion in the following section).

\medskip

V. Unlike the usual Seiberg-Witten invariants (usually defined by
 perturbation argument using ${\cal C}^{\infty}$ topology),
 the ${\cal ASW}$ are defined
 as the intersection numbers of various Chern classes on a neighborhood
 of the algebraic Seiberg-Witten moduli spaces based on the construction of
 algebraic Kuranishi models.  Unless the reduced 
algebraic Seiberg-Witten moduli space (cut off by a finite number of
codimension one cycles, determined by its formal dimension formula) 
 happens to be of zero dimensional, the invariant usually does not 
correspond to actual counting of the number of curves dual to a class 
 $C$. Thus, one should view ${\cal ASW}$ as a formal enumeration 
(in the sense of intersection theory).

\medskip

VI. In most of the cases, the data of algebraic Kuranishi models involves
 algebraic vector bundles and algebraic bundle maps between these vector bundles.
In some minor cases (discussed later), it involves algebraic vector bundles
 and non-algebraic bundle maps between them. We still call the corresponding
 defined invariant ``algebraic'' as the algebraic Seiberg-Witten moduli space 
itself 
 and the Chern classes of the algebraic vector bundles involved in 
defining ${\cal ASW}$ are algebraic objects.
\medskip
VII. Finally, by its definition the ${\cal ASW}$ is ``NOT'' an invariant
 in the traditional sense. Namely, it is not transparent that it is
independent of the complex structure of $M$. But we still call it 
an ``invariant'' because 

(i). it is independent to the choices of the 
algebraic Kuranishi models chosen to define the invariant.

\medskip

(ii). In many situations, it is related to the topological Seiberg-Witten
 invariant.

\medskip

(iii). As will be proved in section \ref{section; ap}, it shares the same 
functorial
 properties under blowing ups as the usual (family) Seiberg-Witten invariants.

\medskip

(iv). After explicit calculation (e.g. by wall crossing formula, etc), usually
 one can compute ${\cal ASW}$ explicitly and find it to be independent of
 the analytic information (like the complex structure) of $M$. 
\bigskip
\section{\bf The Definition of Algebraic Seiberg-Witten Invariant for
  Algebraic Surfaces with zero Geometric Genera}\label{section; da}
\bigskip

 Recall that the original Seiberg-Witten invariant $SW$ is defined
 for all $spin^c$ structures on a given smooth four-manifold with $b^+_2\geq 1$
 (with dependence on
 chamber structures when $b^+_2=1$.
 On the other hand, the usage of $spin^c$ structure is not particularly
 convenient for our discussion of algebraic Seiberg-Witten invariants. Therefore,
 we will
 adopt a slightly different notation in the holomorphic or algebraic category.

 To begin our discussion, we start by recalling (see e.g [FM]) how is the 
Seiberg-Witten
 invariant defined in terms of the holomorphic data.

  Let $(M, \omega_M)$ be a Kahler surface and let $\omega_M$ be the Kahler form.
Then $\omega_M$ splits the $spin^c$ spinor vector bundle into 

 $${\cal S}_{\cal L}^+\cong {\bf E}\oplus {\bf E}\otimes {\bf K}_M^{-1}.$$

 Following the usual convention, one can denote $\alpha\in \Gamma({\bf E})$ and
 $\beta\in \Gamma({\bf E}\otimes {\bf K}_M^{-1})$ the smooth sections of the
 ${\cal C}^{\infty}$ line bundles and the corresponding Dirac equation and
 Seiberg-Witten equation are reduced to scalar valued equation and  read as

  $$\overline{\partial_a}\alpha+\overline{\partial_a}^{\ast}\beta=0,$$
 $$ F^{0,2}_a=\overline{\partial_a}^2=\beta\cdot \overline{\alpha},$$
 $$F^{1,1}_a\wedge \omega_M={i}{|\alpha|^2-|\beta|^2-1\over 2}\omega_M^2.$$

  A standard argument implies that $\alpha\equiv 0$ or $\beta\equiv 0$ 
identically (which depends on whether $c_1({\bf E})\cdot \omega_M>0$ or
 $c_1({\bf E})\cdot \omega_M<0$.  If $c_1({\bf E})\cdot \omega_M>0$,
 the smooth section $\beta$ is identically zero and 

$$\overline{\partial_a}:\Gamma({\bf E})\longrightarrow \Gamma({\bf E}\otimes
 {\omega_M}^{0,1})$$

 satisfies $\overline{\partial_a}^2=0$ and thus  defines a holomorphic structure
 on the ${\cal C}^{\infty}$ line bundle ${\bf E}$. The Dirac equation
 $\overline{\partial_a}\alpha+\overline{\partial_a}^{\ast}\beta=0$ is then reduced
 to the d-bar equation $$\overline{\partial_a} \alpha=0$$
 and $\alpha$ is a non-zero holomorphic section of the holomorphic
 structure induced on the line bundle ${\bf E}$.  The $(1, 1)$-projection of the
Kahler-Seiberg-Witten equation can be reduced to the Kazdan-Warner equation on $M$.

\begin{defin}\label{defin; GT}
 Given a cohomology class $C\in H^2(M, {\bf Z})$, the Gromov-Taubes dimension of 
the class is
 defined to be 

 $$d_{GT}(C)={C\cdot C-c_1(K_M)\cdot C\over 2}.$$
\end{defin}

  The formula is nothing but the dimension formula of Seiberg-Witten invariant, 
formulated in terms of
 $C$ instead of ${\cal L}=2C-c_1(K_M)$ (in additive notation). It was first 
discovered by C. Taubes [T].
\medskip
\subsection{\bf The $p_g=0$ and $q(M)=0$ Case}\label{subsection; 00}
\medskip
 In the following, we discuss the algebraic Seiberg-Witten invariant in the 
easiest case. Later we will extend our
discussion to the $q(M)>0$ case and then to the $p_g>0$ case as well.

Let us fix a $c_1({\bf E})=C\in H^2(M, {\bf Z})$ and define the 
algebraic Seiberg-Witten invariant of $C$ on the algebraic surface $M$. Because
 $h^{2,0}=p_g=0$, the class 
$C$ is automatically a $(1, 1)$ class.

 Because our major interest is to study the holomorphic curves poincare dual to 
$C$, we
 require $C\cdot \omega_M>0$. Otherwise, $C$ cannot be represented by holomorphic
 curves and we simply
 define ${\cal ASW}(C)=0$.
The major tool for the definition of the algebraic Seiberg-Witten invariant is
the existence of the algebraic Kuranishi model, defined in definition 
\ref{defin; tuple} and definition
 \ref{defin; family algebraic}, etc, 

If $q(M)=0$, there exists a unique holomorphic structure on any given  
${\cal C}^{\infty}$
 line bundle ${\bf E}$. By abusing the notation, we use the same symbol
 ${\bf E}$ to denote the holomorphic line bundle and the underlying
 ${\cal C}^{\infty}$ line bundle.

 We start by defining ${\cal ASW}(C)={\cal ASW}_{pt}(1, C)$ in this special 
case. Heuristically speaking, 
the number ${\cal ASW}(C)$ should enumerate the number of
 holomorphic curves dual to $C=c_1({\bf E})$, passing through the 
 $d_{GT}(C)={C^2-C\cdot c_1(K_M)\over 2}$ number of generic points on $M$. 

\medskip
 
 Let $D_C$ be an effective Weil divisor defined by the zero locus of a non-zero 
holomorphic section of 
 ${\bf E}$.  Then
 one has the following short exact sheaves sequence,

$$0\mapsto {\cal O}_M\mapsto {\cal O}_M(D_C)\mapsto 
 {\cal O}_{D_C}(D_C)\mapsto 0.$$

 We have the following vanishing lemma for $h^2(M, {\bf E})$.

\medskip
\begin{lemm}\label{lemm; vanishing}
 Let $M$ be an algebraic surface with $p_g(M)=q(M)=0$. Suppose that
 $h^0_{\overline{\partial}}(M, {\bf E})>0$, 
then $h^2_{\overline{\partial}}(M, {\bf E})=0$.
\end{lemm}

\noindent Proof: 
 Suppose that $h^0(M, {\bf E})>0$, there is a non-zero
 holomorphic section of the holomorphic line bundle ${\bf E}$. Pick one
 non-trivial section and denote the corresponding zero locus by $D_C$.

Then it is well known by spectral sequence argument that
the sheaf cohomology group $H^i(M, {\cal O}_M(D_C))$ is isomorphic to
 d-bar cohomology group $H^i_{\overline{\partial}}(M, {\bf E})$.

To argue the vanishing of $H^2_{\overline{\partial}}(M, {\bf E})$, it 
suffices to show that $H^2(M, {\cal O}_M(D_C))={\bf 0}$.

 Let us take the right derived long exact sequence of the above 
short exact sequence and consider the last three terms,

 $$H^2(M, {\cal O}_M)\mapsto H^2(M, {\cal O}_M(D_C))\mapsto 
 H^2(D_C, {\cal O}_{D_C}(D_C))\mapsto 0.$$

  One observes that the first term and the third term are both trivial.
 The $H^2(D_C, {\cal O}_{D_C}(D_C))$ is trivial because $D_C$ is complex one 
dimensional. On the other hand, $dim_{\bf C}H^2(M, {\cal O}_M)=p_g(M)$ is
 assumed to be zero, then the middle term $H^2(M, {\cal O}_M(D_C))$ must
 be trivial. $\Box$

\medskip

 Knowing the triviality of the second $d-bar$ cohomology, we have the following
 dimension count

$$h^0_{\overline{\partial}}(M, {\bf E})-h^1_{\overline{\partial}}(M, {\bf E})=
h^0(M, {\cal O}_M(D_C))-h^1(M, {\cal O}_M(D_C))$$
$$=\chi({\cal O}_M(D_C))={D_C\cdot D_C-K_M\cdot D_C\over 2}+1=
{C\cdot C-c_1(K_M)\cdot C\over 2}+1.$$

 We have abused the notation $\cdot$ a bit.
 The pairing $D_C\cdot D_C$ denotes the intersection pairing on divisors, while
 the pairing $C\cdot C$ means the cohomology pairing on $H^2(M, {\bf Z})$.

  If $h^1(M, {\cal O}_M(D_C))=0$, then $h^0(M, {\cal O}(D_C))={C^2-c_1(K_M)\cdot
 C\over 2}+1$ and the complete linear system $|D_C|$ is a projective space
 of the expected Gromov-Taubes dimension ${C^2-c_1(K_M)\cdot C\over 2}$.
 The counting indicates that when $h^1(M, {\cal O}_M(D_C))\not=0$, 
the vector space dimension $h^0(M, {\cal O}_M(D_C))$ may be 
different from the ``expected dimension'' ${C^2-c_1(K_M)\cdot C\over 2}+1$ and
 $H^1(M, {\cal O}_M(D_C))\cong H^1(D_C, {\cal O}_{D_C}(D_C))$ represents the 
obstruction space.

  We define the algebraic moduli space of curves dual to $C$ to be
 the projective space formed by the vector space $H^0(M, {\cal O}_M(D_C))\cong
 H^0_{\overline{\partial}}(M, {\bf E})$, which is of 
 $h^0(M, {\cal O}_M(D_C))-1$ dimension.

 We adopt the following formulation which will be extended to $p_g>0$ or 
$q(M)>0$ cases.
 As the algebraic family moduli space is not of the expected dimension, 
 an insertion of the top Chern class of the obstruction bundle is necessary.

\medskip
\begin{defin}\label{defin; tuple}
 Consider the triple $(V, W, \Phi_{VW})$ with $V, W$
 be finite dimensional vector spaces
 and let $\Phi_{VW}$ be a linear map from $V$ to $W$, 

$$\Phi_{VW}:V\longrightarrow W.$$

It is said to be an algebraic Kuranishi model of the class $C$ over $M$
 if $Ker(\Phi_{VW})\cong H^0(M, {\cal O}_M(D_C))$ and
 $Cokernel(\Phi_{VW})\cong H^1(M, {\cal O}_M(D_C))$.
\end{defin}
\medskip

Given an algebraic Kuranishi model of $C$, one defines the algebraic family
Seiberg-Witten invariant by the following recipe.

 Consider the projective space $\pi_{{\bf P}(V)}:{\bf P}(V)\mapsto pt$ and the 
corresponding
 tautological line bundle ${\bf H}^{\ast}$, the dual of hyperplane line bundle. 
The vector space morphism
 $\Phi_{VW}$ induces a bundle map ${\bf H}^{\ast}\mapsto \pi^{\ast}_{{\bf P}(V)}W$ 
and
therefore a canonical section of the obstruction bundle 
${\bf H}\otimes_{\bf C} \pi^{\ast}_{{\bf P}(V)}W$.

 The algebraic cycle of the moduli space of curves poincare dual to $C$ is
represented by the top Chern class 
$c_{top}({\bf H}\otimes_{\bf C}\pi^{\ast}_{{\bf P}(V)}W)$.
 Requiring the curve to pass through a generic point imposes an additional
 $c_1({\bf H})$ insertion.

 One defines ${\cal ASW}(C)$ to be

\begin{defin}\label{defin; ASW}
$${\cal ASW}(C)=
\int_{{\bf P}(V)}c_1^{dim_{\bf C}V-dim_{\bf C}W-1}
({\bf H}) c_{top}({\bf H}\otimes_{\bf C} W).$$
\end{defin}
 
\medskip
\begin{rem}\label{rem; indep}
 It is easy to see that the definition of ${\cal ASW}_{pt}$ is independent of
 the choices of $(V, W, \Phi_{VW})$ and is equal to $1$, which is up to a sign 
 equal to the wall crossing number [KM] of a $b^+_2=1, b_1(M)=0$ four-manifold. 
\end{rem}
\bigskip
 In the following, we extend the definition of algebraic Seiberg-Witten 
invariants to
 $q(M)\not=0$ case.

\medskip
\subsection{\bf The Definition of {\cal ASW} for $q(M)\not=0, p_g=0$ Algebraic
 Surfaces}\label{subsection; q>0}
\bigskip

    Let $M$ be an $p_g=0$ algebraic surface with $q(M)\not=0$ and let $C$
 be an element
 in $H^2(M, {\bf Z})$. As $H^{2,0}(M, {\bf C})=H^{0,2}(M, {\bf C})={\bf 0}$,
 $C$ is automatically a $(1, 1)$ class, which becomes the first Chern class of
 a holomorphic line bundle ${\bf E}$.  However, the holomorphic structures on
 ${\bf E}$  are not unique and form a continuous moduli known to be $Pic^0(M)$, the
 connected component of the Picard group $Pic(M)$ containing identity.

  Consider the Hodge decomposition 

  $$H^1(M, {\bf R})\otimes_{\bf R} {\bf C}=H^{1,0}(M, {\bf C})\oplus
 H^{0, 1}(M, {\bf C}),$$  which induces an real vector space isomorphism

  $$i: H^1(M, {\bf R})\cong H^{0, 1}(M, {\bf C}) \cong H^1(M, {\cal O}_M).$$

The following is well known to algebraic geometers,
\medskip

\begin{prop}\label{prop; torus}
  Let $M$ be a Kahler surface with $q(M)\not=0$, then the connected 
component of the Picard variety $Pic(M)$ containing identity, $Pic^0(M)$,
  can be identified as a complex variety 
to be the following complex torus

  $$T(M)=H^1(M, {\cal O}_M)/i(H^1(M, 
{\bf Z}))\cong H^{0, 1}(M, {\bf C})/i(H^1(M, {\bf Z})).$$

\end{prop}

\medskip
 For a discussion on $Pic^0(M)=T(M)$, please consult [BPV], page 36, section 13.
\medskip 
 Given a ${\cal C}^{\infty}$ topological line bundle, 
there is a universal holomorphic line bundle ${\bf E}$ over $M\times T(M)$ 
such that the restriction, ${\bf E}|_{M\times\{t\}}$, $t\in T(M)$, is the
 holomorphic line bundle over $M$, $c_1({\bf E}|_{M\times\{t\}})=C$, with
 the specific holomorphic structure 
 parametrized by $t\in T(M)$. One way to construct this is to consider the
 Poincare line bundle over $M\times T(M)$
 and tensor it with any holomorphic line bundle over $M$ with first Chern
 class $=C$. 

 Let us consider a class $C\in H^2(M, {\bf Z})$ which satisfies
 $(i). C\cdot C-c_1(K_M)\cdot C\geq 0$, $(ii).$ $C\cdot \omega_M>0$.

 Then the surface Riemann-Roch formula for ${\bf E}$, $c_1({\bf E})=C$, gives 
$$\chi({\bf E})={c_1({\bf E})\cdot c_1({\bf E})-c_1(K_M)\cdot c_1({\bf E})\over 2}+
1-q(M)+p_g$$
$$={C\cdot C-c_1(K_M)\cdot C\over 2}-q(M)+1=d_{GT}(C)-q(M)+1,$$
 which differs from the Gromov-Taubes dimension formula by $1-q(M)$.
\medskip
 Consider the projection $\pi: M\times T(M)\mapsto T(M)$ and 
push forward the sheaf of holomorphic sections of ${\bf E}$, ${\cal E}$, 
along $M\times T(M)\mapsto T(M)$.

 There are three right derived image 
sheaves ${\cal R}^i\pi_{\ast}\bigl({\cal E}\bigr)$,
 $0\leq i\leq 2$, to consider. 

\medskip

\noindent Case 1. $(c_1(K_M)-C)\cdot \omega_M<0$. For most of the classes $C$
 with high enough degree (energy) $C\cdot \omega_M\gg 0$, 
the pairing $c_1(K_M)\cdot \omega_M-C\cdot \omega_M$ is negative. \label{case}

 Under this additional assumption on $C$, it follows from relative Serre
 duality that the 
sheaf ${\cal R}^2\pi_{\ast}\bigl({\cal E}\bigr)$ over $T(M)$ vanishes. 
Otherwise $K_M\otimes {\bf E}^{\ast}$ has non-zero sections, which implies
 that the degree (energy)
$c_1(K_M\otimes {\bf E}^{\ast})\cdot \omega_M=c_1(K_M)\cdot \omega_M
 -C\cdot \omega_M>0$.

The sheaves ${\cal R}^0\pi_{\ast}\bigl( {\cal E}\bigr)$ and
 ${\cal R}^1\pi_{\ast}\bigl({\cal E}\bigr)$ may not be locally free.
When $t\in T(M)$ moves, the dimension of $H^0(M, {\bf E}|_{M\times t})$
 may vary and can even be zero at the generic points of $T(M)$.

 Given an effective divisor $D_C$, which is the zero locus of some
 holomorphic section of ${\bf E}_{M\times t}, t\in T(M)$, the 
 short exact sequence 

$$0\mapsto {\cal O}_{M}\mapsto {\cal O}_{M}(D_C)\mapsto
 {\cal O}_{D_C}(D_C)\mapsto 0$$
 
 has a corresponding long exact sequence for $p_g=0$ surfaces,

$$0\mapsto {\bf C}\mapsto H^0(M, {\cal O}_M(D_C))\mapsto
 H^0(D_C, {\cal O}_{D_C}(D_C))\mapsto H^1(M, {\cal O}_M)$$
$$\mapsto
 H^1(M, {\cal O}_M(D_C))\mapsto H^1(D_C, {\cal O}_{D_C}(D_C))\mapsto 0.$$

 The space $H^0(D_C, {\cal O}_{D_C}(D_C))$ represents the infinitesimal
 deformations of the curve $D_C$ in $M$. 

 If $H^1(M, {\cal O}_M)$ is identified with the tangent space of
 $T(M)$ at $t$, the exactness of

$$H^0(M, {\cal O}_M(D_C))\mapsto
 H^0(D_C, {\cal O}_{D_C}(D_C))\mapsto H^1(M, {\cal O}_M)$$

 indicates that the infinitesimal deformation of the curve $D_C\subset M$
 comes from two sources. Either it comes from the deformation of the
 holomorphic sections in the same linear system ${\bf P}(H^0(M, {\cal O}_M(D_C)))$,
 or it comes from the variation of holomorphic structures parametrized by $t$.

 On the other hand, the exactness of 
$$ H^0(D_C, {\cal O}_{D_C}(D_C))\mapsto H^1(M, {\cal O}_M)\mapsto
 H^1(M, {\cal O}_M(D_C))$$

 indicates that the tangent direction in ${\bf T}_t(T(M))$ which gives rise
 to non-trivial infinitesimal curve deformations is in the kernel of the
 tangent obstruction map $${\bf T}_t(T(M))\cong H^1(M, {\cal O}_M)\mapsto
 H^1(M, {\cal O}_M(D_C)).$$

\medskip
 Following the same idea as before, we define the algebraic 
family Kuranishi model of the class $C$,

\medskip
\begin{defin}\label{defin; family algebraic}
Let ${\bf V}, {\bf W}$ be two algebraic vector bundles over the $q(M)$ 
dimensional torus $T(M)$ and let $\Phi_{\bf VW}:{\bf V}\mapsto {\bf W}$
 be an algebraic bundle map from ${\bf V}$ to ${\bf W}$. 
We denote the corresponding sheaves of sections and the sheaf morphism
 by ${\cal V}, {\cal W}, \Phi_{\cal VW}$, respectively.
  
 The data $({\bf V}, {\bf W}, \Phi_{\bf VW})$ is said to be an algebraic 
 Kuranishi model of $C$ over $M\times T(M)$ if there exist sheaf 
isomorphisms
$Ker(\Phi_{\cal VW})\cong {\cal R}^0\pi_{\ast}\bigl({\cal E}\bigr)$ and 
$Coker(\Phi_{\cal VW})\cong {\cal R}^1\pi_{\ast}\bigl({\cal E}\bigr)$.
\end{defin}

\medskip
In the following, we prove the existence of the algebraic Kuranishi model of
 $C$.
\medskip

\begin{prop}\label{prop; existence}
  Suppose that $C\in H^2(M, {\bf Z})$ on an algebraic surface $M$ with
 $p_g(M)=0, q(M)>0$ satisfies

 $(i). C\cdot C-c_1(K_M)\cdot C\geq 0$,

 $(ii). \bigl(c_1(K_M)-C\bigr)\cdot \omega_M<0$.

 Then there exists an algebraic Kuranishi model of $C$ over $M\times T(M)$.
\end{prop}

\medskip

\noindent Proof of the proposition: Because $M$ is algebraic, we consider an 
ample effective 
divisor $D$ on $M$. Then ${\bf D}=D\times T(M)$ is relative ample along 
$\pi:M\times T(M)\mapsto T(M)$.

 We consider the short exact sequence 

$$0\mapsto {\cal O}_{M\times T(M)}\mapsto {\cal O}_{M\times T(M)}(n{\bf D})\mapsto
{\cal O}_{n{\bf D}}(n{\bf D})\mapsto 0,$$

 for a large $n$, whose value is yet to be determined.

 By tensoring with ${\cal E}$ and by taking the right derived long exact 
sequence,

$$0\mapsto {\cal R}^0\pi_{\ast}\bigl({\cal E}\bigr)\mapsto 
 {\cal R}^0\pi_{\ast}\bigl( {\cal E}(n{\bf D})\bigr)\mapsto $$
$${\cal R}^0\pi_{\ast}\bigl({\cal O}_{n{\bf D}}\otimes {\cal E}(n{\bf D})\bigr)
\mapsto {\cal R}^1\pi_{\ast}\bigl({\cal E}\bigr)$$
$$\mapsto {\cal R}^1\pi_{\ast}\bigl( {\cal E}(n{\bf D})\bigr)\mapsto
{\cal R}^1\pi_{\ast}\bigl({\cal O}_{n{\bf D}}\otimes {\cal E}(n{\bf D})\bigr)
\mapsto 0,$$

 where we have used the vanishing of ${\cal R}^2\pi_{\ast}\bigl( {\cal E} \bigr)$
 over $T(M)$.

 To show that ${\cal R}^0\pi_{\ast}\bigl( {\cal E}(n{\bf D})\bigr)$
 and ${\cal R}^0\pi_{\ast} \pi_{\ast}\bigl({\cal O}_{n{\bf D}}\otimes 
{\cal E}(n{\bf D})\bigr)$
 are locally free and the sequence is truncated to a four-term long exact 
sequence, it suffices to prove that the sheaf 
 ${\cal R}^1\pi_{\ast}\bigl( {\cal E}(n{\bf D})\bigr)$ vanishes on $T(M)$.
 Once this has been done, the sheaf exact sequence truncates to a four-term
 exact sequence and 
${\cal R}^1\pi_{\ast}\bigl({\cal O}_{n{\bf D}}\otimes {\cal E}(n{\bf D})\bigr)$
 vanishes as well. 

 Then ${\cal R}^0\pi_{\ast}\bigl({\cal O}_{n{\bf D}}\otimes {\cal E}(n{\bf D})
\bigr)$ is locally free as 
 $H^0(n{\bf D}, {\cal O}_{n{\bf D}}(n{\bf D})\otimes {\cal E}|_{M\times t})$ 
is of constant rank
 independent of $t\in T(M)$ (see [Ha] theorem 21.11).

 To show that ${\cal R}^1\pi_{\ast}\bigl( {\cal E}(n{\bf D})\bigr)$ vanishes on 
$T(M)$,
 it suffices to check that 
$h^1(M, {\cal O}_M(n{\bf D})\otimes_{{\cal O}_M} {\cal E}|_{M\times t})$ is zero 
for all $t\in 
 T(M)$.

  By Nakai criterion, one can make ${\cal O}(nD)\otimes K_M^{-1}\otimes 
{\cal E}|_{M\times \{t\}}$ ample
 if $n$ is chosen to be large enough. Then the vanishing of the
 $h^1(M, {\cal O}_M(nD)\otimes {\cal E}|_{M\times t})$ is a simple 
consequence of Kodaira vanishing theorem. $\Box$
\medskip

  We define the algebraic 
Seiberg-Witten invariant of $C$ to be 

 $${\cal ASW}(C)=\int_{{\bf P}({\bf V})}c_1^{d(C)}({\bf H})c_{top}(
 {\bf H}\otimes \pi^{\ast}_{{\bf P}({\bf V})}{\bf W})$$

 $$={\bf c}_1^{d(C)}({\bf H})\cap {\bf c}_{top}(
 {\bf H}\otimes \pi^{\ast}_{{\bf P}({\bf V})}{\bf W})\in 
{\cal A}_0({\bf P}({\bf V})),$$

 where the bold face ${\bf c}_i$ denotes the Chern classes as algebraic 
cycle classes 
$\in {\cal A}_{\ast}({\bf P}({\bf V}))$. The $\cap$ is the intersection 
product on the cycles.

 The definition is independent to the choices of the algebraic Kuranishi models 
 of $C$ as we notice that ${\cal ASW}$ can be computed by using 
Grothendieck-Riemann-Roch formula on $[{\bf V}-{\bf W}]\in K(T(M))$ 
and is identical to the wall crossing calculation of $2C-K_M$ performed in [LL2].
 \medskip
\begin{rem}\label{rem; restrict}
  The algebraic Seiberg-Witten invariant defined in this way corresponds to
 counting curves from all the possible holomorphic structures on the line bundle
  with first Chern class $C$. To count curves from a fixed holomorphic 
structure (linear system) we have to modify the definition of ${\cal ASW}$
 by the following recipe

 $$\int_{{\bf P}({\bf V})}c_1^{d(C)-q(M)}({\bf H})\pi^{\ast}_{{\bf P}({\bf V})}
[T(M)]c_{top}({\bf H}\otimes \pi^{\ast}_{{\bf P}({\bf V})}
 {\bf W}),$$ by inserting the pull back of the fundamental class $[T(M)]$ and
 dropping
 the power of $c_1({\bf H})$.

 In general, one can consider a finite number of classes in $H_1(M, {\bf Z})$
 which generate a finite number of elements $\eta_i$ in $H^1(T(M), {\bf Z})$
 through the isomorphism $H_1(M, {\bf Z})\cong H^1(T(M), {\bf Z})$. One
 may insert the cup product of these elements $\wedge \eta_i$ 
into the integral (intersection pairing) and the
 corresponding algebraic Seiberg-Witten invariant counts the algebraic 
curves from the holomorphic structures in a locus of $T(M)$ poincare dual to  
$\wedge \eta_i\in H^{\ast}(T(M), {\bf Z})$.

 The corresponding wall crossing number has been calculated in [L][LL1], and
 we omit the details here.

\end{rem}

\bigskip

\noindent Case 2: $(c_1(K_M)-C)\cdot \omega_M>0$.
\medskip
 If both pairings $C\cdot \omega_M$ and $(c_1(K_M)-C)\cdot \omega_M$ are positive,
 ${\cal R}^2\pi_{\ast}\bigl({\cal E}\bigr)$ may not be trivial and the
 previous argument is not applicable.

 Nevertheless, we still have the following proposition,

\begin{prop}\label{prop; exclude}
 The supports of the coherent sheaves ${\cal R}^0\pi_{\ast}\bigl({\cal E}\bigr)$
 and ${\cal R}^2\pi_{\ast}\bigl({\cal E}\bigr)$ do not intersect on $T(M)$.
\end{prop}
\medskip
\noindent Proof: 
The proposition is in fact a generalization of lemma 
\ref{lemm; vanishing}. 
Firstly, we mark that by [Ha] chapter 2, section 5, exercise 5.6(c), 
the support of a coherent sheaf is closed.
 If the intersection of  these two supports 
 is non-empty, then by base change there exists at least a $t\in T(M)$ such that
$H^0(M, {\cal E}|_{M\times t})$ and $H^2(M, {\cal E}|_{M\times t})$ are both
 non-trivial. Then the chosen $t\in T(M)$ defines a holomorphic line bundle
 structure on the
 ${\cal C}^{\infty}$ line bundle with first Chern class $C$. Then
 one applies lemma \ref{lemm; vanishing} to this situation and finds the
 contradiction! Thus, the supports of the two coherent sheaves can never
 intersect. $\Box$

\bigskip

 The positivity of both $(c_1(K_M)-C)\cdot \omega_M $ and $C\cdot \omega_M$ 
implies
 that $c_1(K_M)\cdot \omega_M>0$.

 Given the class $c_1(K_M)$ with positive pairing with $\omega_M$, we separate
 into two cases depending on whether $c_1(K_M)$ is poincare dual to a holomorphic
 curve in $M$ or not.

\medskip
\noindent (i). $c_1(K_M)$ is poincare dual to a holomorphic curve in $M$.

\medskip

 Because $p_g=dim_{\bf C}H^0(M, K_M)=0$, the holomorphic curve is an effective
 divisor of some holomorphic line bundle which has an
 underlying ${\cal C}^{\infty}$
line bundle as $K_M$.

 According to surface classification result, such algebraic surface is 
either a general type surface or an elliptic surface. For both types of surfaces,
 $K_{M_{min}}$ is numerically effective. Namely, $c_1(K_{M_{min}})$ has 
non-negative pairings with effective classes.

  By repeating blowing down $M$, one may get a unique minimal model of $M$,
 denote by $M_{min}$, and the adjunction equality states that

$$c_1(K_M)=c_1(K_{M_{min}})+\sum E_i,$$

where $E_i\in H^2(M, {\bf Z})$
 are the exceptional $-1$ class of the blowing down map
 $M\mapsto M_{min}$.

 \medskip
 
We may write $C$ as $C_{min}+\sum m_i E_i$, with
$C_{min}\in H^2(M_{min}, {\bf Z})\subset H^2(M, {\bf Z})$.
A direct calculation shows that 
$$C_{min}^2-c_1(K_{min})\cdot C_{min}=C^2-c_1(K_M)\cdot C+\sum (m_i^2-m_i)
\geq C^2-c_1(K_M)\cdot C\geq 0,$$
 independent of the signs of these $m_i$.

  If $C_{min}=0\in H^2(M_{min}, {\bf Z})$, then $C=\sum m_i E_i$ is a multiple
 of exceptional classes. The condition $C^2-C\cdot K_M=-\sum (m_i^2-m_i)\geq 0$
 forces $m_i(m_i-1)=0$ for all $i$ and therefore either $m_i=1$ or $m_i=0$.

  Define $J=\{i|m_i=1, 1\leq i\leq n\}$.
Thus, $C=\sum_{i\in J} E_i$ can and only can be
 represented by the sum of holomorphic $-1$ curves.  
In this case, we define ${\cal ASW}(\sum E_i)$ to be 
$1$.
\medskip
  If $C_{min}\not= 0$, then the image of the holomorphic curve
 dual to $C$ under the blowing down map $M\mapsto M_{min}$ is a holomorphic
 curve dual to $C_{min}$.

 Denote $L=C_{min}-c_1(K_{M_{min}})\in H^2(M_{min}, {\bf Z})$.
  If $L=0$, $C$ has to be equal to $c_1(K_{M_{min}})+\sum_{i\in J} E_i$ for some 
 sub-collection $J$ of $-1$ classes.
 If $-L$ is not
 poincare dual to any holomorphic curve in $M_{min}$, then 
$c_1(K_M)-C=-L+\sum_{1\leq i\leq n}
 (1-m_i)E_i$ is not representable by holomorphic curves in $M$, either.
 In this case, ${\cal R}^2\pi_{\ast}\bigl({\cal E}\bigr)$ vanishes. The
 definition of ${\cal ASW}(C)$ is identical to the 
 $(c_1(K_M)-C)\cdot \omega_M<0$ case. 

 So we may suppose that $-L$ is representable by holomorphic curves in 
$M_{min}$ and therefore $-L\cdot \omega_{M_{min}}>0$.

  Because $K_{M_{min}}$ is numerically effective and both $C_{min}$ and
 $c_1(K_{M_{min}})-C_{min}$ are represented by holomorphic curves in $M_{min}$,
  it follows that $c_1(K_{M_{min}})\cdot C_{min}$ and $c_1(K_{M_{min}})\cdot
 (c_1(K_{M_{min}})-C_{min})$ are non-negative. 

 Therefore 
$$C^2_{min}=(C^2_{min}-c_1(K_{M_{min}})\cdot C_{min})+c_1(K_{M_{min}})\cdot C_{min}
 \geq 0+0=0,$$

$$(c_1(K_{M_{min}})-C_{min})^2=\bigl((c_1(K_{M_{min}})-C_{min})^2-c_1(K_{M_{min}})
\cdot  (c_1(K_{M_{min}})-C_{min})  \bigr)$$
$$+c_1(K_{M_{min}})\cdot 
(c_1(K_{M_{min}})-C_{min})$$
$$=(C_{min}^2-c_1(K_{M_{min}})\cdot C_{min})+c_1(K_{M_{min}})\cdot 
(c_1(K_{M_{min}})-C_{min})\geq 0+0=0.$$

Because both the classes $c_1(K_{M_{min}})-C_{min}$ and $C_{min}$ are
 in the forward light cone, 
then by the line cone lemma [LL3], their intersection pairing

$$\bigl(c_1(K_{M_{min}})-C_{min}\bigr)\cdot
 C_{min}=-(C_{min}^2-C_{min}\cdot c_1(K_{M_{min}}))\geq 0.$$  

 Thus $\bigl(c_1(K_{M_{min}})-C_{min}\bigr)\cdot C_{min}=0$, which can
 only occur either
 if $c_1(K_{M_{min}})-C_{min}$ and $C_{min}$ are parallel to each other 
(up to torsions) in 
 $H^2(M, {\bf Z})$ and both $c_1(K_{M_{min}})-C_{min}$ and $C_{min}$ lie
 on the boundary
 of the light cone, i.e. $M_{min}$ is a minimal elliptic surface with 
$c_1^2(K_{M_{min}})=0$
 and $C_{min}=rc_1(K_{M_{min}})$ (up to torsions) for some $r\in {\bf Q}, 
|r|\leq 1$.

 The definition of the algebraic Seiberg-Witten invariant ${\cal ASW}(C)$ in 
this situation deserves some additional discussion.  Because of the presence 
of the sheaf 
 ${\cal R}^2\pi_{\ast}\bigl({\cal E}\bigr)$, the algebraic 
Kuranishi model for $C$ is not defined 
over the whole $T(M)$ as in the $(c_1(K_M)-C)\cdot \omega_M<0$ case.

 Instead one considers the support of 
${\cal R}^0\pi_{\ast}\bigl({\cal E}\bigr)$ and
 ${\cal R}^2\pi_{\ast}\bigl({\cal E}\bigr)$ and denote
 them by $Z_0\subset T(M)$ and $Z_2\subset T(M)$, respectively. 
 According to [Ha] section 2.5, exercise 5.6., $Z_i, i=0, 2$ are  
compact sub-varieties of $T(M)$.

 As before consider the derived long exact sequence of the short exact sequence 

$$0\mapsto {\cal E}\mapsto {\cal E}(n{\bf D})\mapsto 
{\cal O}_{n{\bf D}}\otimes {\cal E}(n{\bf D})\mapsto 0,$$

 for an ample divisor $D\subset M$, ${\bf D}=D\times T(M)\subset M\times T(M)$, 
and a sufficiently large $n$.

 As in the $ (c_1(K_M)-C)\cdot \omega_M<0$ case,
 the sheaf ${\cal R}^2\pi_{\ast}\bigl({\cal E}\bigr)$
 is trivial over Zariski open $Z^c_2=T(M)-Z_2$. Thus, one still
 can prove the existence of

  $$\Phi_{\cal VW}:{\cal V}\longrightarrow {\cal W}$$
  over $Z^c_2$, where ${\cal V}$ and ${\cal W}$ are locally free
 sheaves over $Z^c_2$ such that
 $Ker(\Phi_{\cal VW})\cong {\cal R}^0\pi_{\ast}\bigl({\cal E}|_{M\times 
Z^c_2}\bigr)$ and
 $Coker(\Phi_{\cal VW})\cong {\cal R}^1\pi_{\ast}\bigl({\cal E}_{M\times 
Z^c_2}\bigr)$.

  Let ${\bf V}, {\bf W}$ be the corresponding algebraic 
vector bundles over $Z^c_2$, then we still 
 have $rank_{\bf C}({\bf V}-{\bf W})={C^2-c_1(K_M)\cdot C\over 2}-q(M)+1$. 
 
  Consider the projective space bundle 
$\pi_{{\bf P}}:{\bf P}_{Z^c_2}({\bf V})\mapsto Z_2^c$ as the 
ambient space. Then $c_{top}({\bf H}\otimes \pi_{{\bf P}}^{\ast}{\bf W})$
 determines an algebraic cycle class $[{\cal M}_C]\subset {\bf P}_{Z^c_2}({\bf V})$
 of dimension $q(M)+(rank_{\bf C}{\bf V}-1)-rank_{\bf C}{\bf W}=0$ 
representing $c_{top}({\bf H}\otimes \pi_{{\bf P}}^{\ast}{\bf W})$. 
The zero dimensional cycle class
 is an integral multiple $m_C$ of the generator $[pt]\in
 {\cal A}_0({\bf P}_{Z^c_2}({\bf V}))$.

   We define ${\cal ASW}(C)$ to be $m_C$.
\medskip
\begin{rem}\label{rem; fiber}
  One may construct an algebraic Kuranishi model for both ${\cal 
R}^i\pi\bigl({\cal E}\bigr), 
i=0, 2$ over $T(M)$, but the $c_{top}({\bf H}\otimes {\bf W})$ over
 the projective space bundle
 ${\bf P}_{T(M)}({\bf V})$ represents the wall
 crossing number of the associated 
$spin^c$ structure and is calculable by the universal wall crossing
 formula [LL2]. It is equal to zero for $C=\alpha c_1(K_M)+\sum E_i$ on $p_g=0$ 
elliptic surfaces, which indicates that besides the algebraic cycle class
 associated with
 ${\cal R}^0\pi_{\ast}\bigl({\cal E}\bigr)$, the algebraic cycle class 
associated with
 ${\cal R}^2\pi_{\ast}\bigl({\cal E}\bigr)$ gives an opposite contribution.

 This indicates that the algebraic Seiberg-Witten invariant, which is 
supposed to be 
equal to the Seiberg-Witten invariant under large deformation on the Seiberg-Witten
 equations by the symplectic form, is identical to the Seiberg-Witten
 invariant in the
 metric chamber (because the wall crossing number is $0$). 

 In fact, for an elliptic surface $M \mapsto \Sigma$ over a higher 
genus curve $\Sigma$,
 $T(M)$ can be identified with $J(\Sigma)$, the Jacobian variety of 
$\Sigma$ and the
 Seiberg-Witten invariant of the $spin^c$ class $2C-K_{M}$ (in 
additive notation) was
 determined in [FM2] and was closed related to the intersection theory on 
the symmetric product ${\bf S}^d\Sigma$ for some $d\in {\bf N}$. The
 space $S^d\Sigma$ can be 
 identified with $\coprod_{t\in Z_0}{\bf P}\bigl(H^0(M, {\bf E}_{M\times
 \{t\}})\bigr)$ in 
 our picture and the support $Z_0$ can be identified with the image of 
$S^d\Sigma\mapsto
 J(\Sigma)=T(M)$.

 For the details of the enumeration of ${\cal ASW}(C)$ based on curve theory, 
please consult [FM2].
\end{rem}

\medskip
\noindent (ii). $c_1(K_M)$ is not poincare dual to any holomorphic curve in $M$.

\medskip

If ${\cal R}^0\pi_{\ast}\bigl({\cal E}\bigr)$ is trivial, namely, for
 all $t\in T(M)$ the global sections 
$\Gamma(M\times \{t\}, {\cal E}_{M\times \{t\}})={\bf 0}$, then 
$C$ is not 
 poincare dual to any holomorphic curve in $M$. 
 In this case, we simply define ${\cal ASW}(C)$ to be $0$.

 On the other hand, suppose ${\cal R}^0\pi_{\ast}\bigl({\cal E}\bigr)$
 is not trivial, then $C$ is poincare dual to some holomorphic curve in $M$.
 This implies that $c_1(K_M)-C$ can not be poincare dual to any holomorphic
 curve in $M$. Otherwise, the sum of the homology classes of the holomorphic
 curves dual to $C$ and $c_1(K_M)-C$ is poincare dual to $c_1(K)$, violating
 the assumption that $c_1(K_M)$ is not represented by holomorphic curves in $M$.

 Thus, the sheaf ${\cal R}^2\pi_{\ast}\bigl({\cal E}\bigr)$ must vanish. Then
we define ${\cal ASW}(C)$ by the same recipe as in $(c_1(K_M)-C)\cdot
 \omega_M<0$ case on page \pageref{case}.

\medskip

\bigskip
\subsection{The Algebraic Family Seiberg-Witten Invariants for $p_g>0$ Algebraic 
Surfaces}\label{subsection; pg>0}
\bigskip
  In the previous subsection, we defined the algebraic Seiberg-Witten 
invariants for $p_g=0$ surfaces, based on discussions on the class $C$. 
The invariants defined as intersection numbers of (subspaces) of projective
 space bundles are equal to the 
 topological Seiberg-Witten invariants of these surfaces in the specific
 chambers corresponding to large deformations of symplectic ( Kahler) forms.
 
 In general, their values are directly related to the wall crossing formula
 calculated in [KM], [LL2].

  However, the ordinary Seiberg-Witten invariant (which are equal to 
Gromov-Taubes invariant through Taubes' theorem 'SW=Gr') for $p_g>0$
 algebraic surface has the simple type property. Namely, the invariant
 vanishes if the expected dimension $d_{GT}(C)$ of the $spin^c$ structure is 
positive.

 On the other hand, we will define a version of algebraic Seiberg-Witten
 invariant of $C\in H^2(M, {\bf Z})$, which is non-zero for most classes
 $C$ with $C\cdot\omega_M>0, C^2-c_1(K_M)\cdot C=2d_{GT}(C)\geq 0$.  The readers 
should
 be aware that ${\cal ASW}$ defined in this section is not directly related
 to the usual Seiberg-Witten invariant $SW(2C-K_M)$ of the $spin^c$ class
 $2C-K_M$.

 It turns out that for most of the $C$ specified above, we can construct algebraic
 Kuranishi models by algebraic vector bundles and algebraic bundle morphisms.
 But for a
 few minor cases, we
 have to move out of the algebraic category and consider non-algebraic bundle 
maps, despite
 that the moduli space of curves dual to $C$ is still an algebraic object. 

  Firstly, We begin by addressing the definition of 
 the invariant ${\cal ASW}$ and then at the end
 we give a brief discussion about its relationship with the family Seiberg-Witten
 invariant.

  The key difference from the usual Seiberg-Witten invariant is the dimension 
formula
 of ${\cal ASW}$.

 For a class $C\in H^2(M, {\bf Z})\cap H^{1,1}(M, {\bf C})$ with $C\cdot 
\omega_M>0$, 
its expected family
 moduli space dimension
 is ${C^2-C\cdot c_1(K_M)\over 2}+fbd(C, M)$, where $fbd(C, M)\in {\bf N}\cup
 \{0\}$,
 the abbreviation of the formal base dimension, is a correction
 term depending on $C$.

 As usual, let ${\cal E}\mapsto M\times T(M)$ denote the locally free sheaf 
associated to ${\bf E}$,
 with the holomorphic structures parametrized by $T(M)$.
\medskip
\begin{defin}\label{defin; r}
  If ${\cal R}^0\pi_{\ast}\bigl({\cal E}\bigr)$ is the zero sheaf over $T(M)$, 
then define $fbd(C, M)=0$.
  Otherwise, $C$ is poincare dual to the fundamental class of a holomorphic 
curve, which is
 the zero locus of some $s_t\in H^0(M, {\cal E}_{M\times \{t\}})$.

  Given the pair $(s_t, t)$, $t\in T(M)$ and $s_t\in H^0(M, {\cal E}_{M\times
 \{t\}})-\{0\}$, the
 tensor product with $s_t$ induces the sheaf morphism

 $$\otimes s_t:{\cal O}_M\longrightarrow {\cal E}_{M\times \{t\}}$$
 and then the induced morphism of cohomologies

 $$(\otimes s_t)_{\ast}: H^2(M, {\cal O}_M)\mapsto H^2(M, {\cal E}_{M\times 
\{t\}}).$$

  Then the number $fbd(C, M)$ is defined to be the dimension of $\cap_{(s_t, t)}
 Ker\bigl((
\otimes s_t)_{\ast}\bigr)$, with $t\in T(M)$ and $s_t\in H^0(M, {\cal E}_{M\times 
\{t\}})-\{0\}$.
\end{defin}
\medskip
 It follows from the definition that $0\leq fbd(C, M)\leq dim_{\bf C}H^2(M, 
{\cal O}_M)=p_g$.
\medskip
 In the following, we assume that $C\in H^2(M, {\bf Z})\cap H^{1,1}(M, 
{\bf C})$ satisfies 
 $C\cdot \omega_M>0$ and ${C^2-c_1(K_M)\cdot C\over 2}+fbd(C, M)\geq 0$.

\medskip
 As in the $p_g(M)=0$ case, we separate into different cases.

\medskip

\noindent Case 1: $(c_1(K_M)-C)\cdot \omega_M<0$:  This condition implies the
 vanishing of the second derived image sheaf ${\cal R}^2\pi_{\ast}
\bigl({\cal E}\bigr)$.
 It follows from the definition of $fbd$ that $fbd(C, M)=p_g$.

 As in the $p_g=0$ case, we can construct algebraic Kuranishi model
 $({\bf V}, {\bf W}, \Phi_{\bf VW})$ over $T(M)$, with the
 understanding that $T(M)=Pic^0(M)$ reduces to a point when $q(M)=0$. 
 Then one defines ${\cal ASW}(C)$ to be

$$\int_{{\bf P}_{T(M)}({\bf V})}c_1^{{C^2-C\cdot c_1(K_M)\over 2}+p_g}({\bf H})
\cup c_{top}({\bf H}\otimes \pi_{{\bf P}({\bf V})}^{\ast}{\bf W})$$

$${\bf c}_1^{{C^2-C\cdot c_1(K_M)\over 2}+p_g}({\bf H})
\cap {\bf c}_{top}({\bf H}\otimes \pi_{{\bf P}({\bf V})}^{\ast}{\bf W})\in 
{\cal A}_0({\bf P}_{T(M)}).$$

\bigskip

\noindent Case 2: $(c_1(K_M)-C)\cdot\omega_M>0$: If the cohomology class 
$c_1(K_M)-C$ is not represented as holomorphic curves, an easy calculation 
shows that 
$fbd(C, M)=p_g$ and  the discussion is identical to
the $p_g(M)=0$ case with the shifting of dimension formula 
 ${C^2-C\cdot c_1(K_M)\over 2}\mapsto {C^2-C\cdot c_1(K_M)\over 2}+p_g$.
 We focus upon the more interesting case that both $C$ and
 $c_1(K_M)-C$ are represented by holomorphic curves. 

\medskip

\noindent (i). 
$2d_{GT}(C)=C^2-c_1(K_M)\cdot C\geq 0$. Then the same argument as in the 
$p_g=0$ case shows
 that either $C$ is a sum of $-1$ classes or
 as before $M$ must be
 an elliptic surface and $C=\alpha c_1(K_{M_{min}})+\sum E_i$. 
 
 If $C=\sum E_i$, then $C^2-c_1(K_M)\cdot C=0$. On the other hand, the 
following lemma 
characterizes $fbd(C, M)$ uniquely.

\medskip
\begin{lemm}\label{lemm; =0}
 Let $M$ be an algebraic surface with $p_g(M)>0$ and $M_{min}$ 
denotes the unique minimal
 model of $M$. Let $E_1, E_2, \cdots E_n$ denote the $-1$ classes 
$\in H^2(M, {\bf Z})-H^2(M_{min}, {\bf Z})$ of the blowing down map 
$p:M\mapsto M_{min}$.
 Then for $C=\sum_{i\in I}E_i, I\subset \{1, 2,\cdots, n\}$, the 
number $fbd(C, M)$ is zero.
\end{lemm}
\medskip
\noindent Proof of the lemma:  By adjunction formula $c_1(K_M)=
c_1(K_{M_{min}})+
\sum_{i=1}^{i=n}E_i$.
 Let ${\cal K}_M ({\cal K}_{M_{min}})$ denote the canonical sheaf 
associated to $K_M 
$ $(K_{M_{min}})$, 
respectively.  
By abusing the notation, we use the same symbol $E_i$ to denote the 
$-1$ cohomology class 
 and the unique effective $-1$ divisor it associates with. As usual 
${\cal E}$ is the locally
 free sheaf over $M\times T(M)$ with $c_1({\cal E})=\sum_{i\in I}E_i$.
 
 Then $$H^0(M, {\cal K}_M\otimes {\cal O}(-\sum_{i\in I}E_i))
\longrightarrow H^0(M, {\cal K}_M),
 I\subset \{1, 2\cdots, n\}$$ is an isomorphism.  Denote $p:M\mapsto 
M_{min}$ to be
 the blowing down map.
The isomorphism follows from rewriting 
 ${\cal K}_M\otimes {\cal O}(-\sum_{i\in I}E_i)$ and ${\cal K}_M$ as 
$p^{\ast}{\cal K}_{M_{min}}\otimes {\cal O}(\sum_{i\notin I}E_i)$ and 
$p^{\ast}{\cal K}_{M_{min}}\otimes
 {\cal O}(\sum_{i=1}^nE_i)$ and there is a commutative diagram of isomorphisms

\[
\begin{array}{ccc}
H^0(M_{min}, {\cal K}_{M_{min}})&\longrightarrow & 
H^0(M, p^{\ast}{\cal K}_{M_{min}})) \\
\Big\downarrow & &\Big\downarrow \\
H^0(M, p^{\ast}{\cal K}_{M_{min}}\otimes {\cal O}(\sum_{i\not\in I}E_i)& 
\longrightarrow & H^0(M, {\cal K}_M)  
\end{array}
\]

   By Serre duality, this implies that for some $t\in T(M)$, $H^2(M, 
{\cal O}_M)\mapsto
 H^2(M, {\cal E}_{M\times \{t\}})$ is an isomorphism. Thus, the formal 
base dimension 
$fbd(C, M)=0$ for such classes.
$\Box$

 This implies that the expected dimension 
of the classes $\sum_{i\in I} E_i$ is $0$, the same as $p_g(M)=0$ case.
The moduli space of curves is a single regular point and 
${\cal ASW}(\sum_{i\in I}E_i)$
 is defined to be $1$.
\medskip

\noindent (i)'. If $C=\alpha c_1(K_{M_{min}})+\sum_{i\in I}E_i$, 
$0\not=|\alpha|\leq 1$, $I\subset \{1, 2, \cdots
, n\}$,
 on an elliptic surface $M$ with the minimal model $M_{min}$, then 
$f: M_{min}\mapsto \Sigma$ 
is an elliptic fibration without multiple fibers. $K_{M_{min}}=f^{\ast}K$, 
where $K$ is a 
line bundle over $\Sigma$, $deg_{\Sigma}K=g(\Sigma)-1+p_g(M)$. \label{1'}

 In this case we know that, 

\medskip

\begin{prop}\label{prop; fbd=0}
 Let $C=\alpha c_1(K_{M_{min}})+\sum_{i\in I}E_i$, $0\not=|\alpha|\leq 1$
 on an elliptic surface $M$, then
 $fbd(C, M)=0$.

\end{prop}
\medskip
\noindent Proof of the proposition: 
We first show that $fbd(C, M)=fbd(\alpha c_1(K_{M_{min}}), M_{min})$ 
for
 such classes. 

  Let ${\cal E}_{min}$ be the locally free sheaf over $M_{min}\times T(M_{min})$
 with $c_1({\cal E}_{min})=C_{min}=\alpha c_1(K_{M_{min}})\in H^2(M_{min},
{\bf Z})$. Then for the blowing down map $p:M\mapsto M_{min}$, ${\cal E}=
p^{\ast}{\cal E}_{min}\otimes
 {\cal O}_M(\sum_{i\in I}E_i)$.
 Pick an arbitrary $t\in T(M_{min})\cong T(M)$ and 
  consider the following commutative diagram 

\[
\begin{array}{ccc}
 H^2(M_{min}, {\cal O}_{M_{min}}) & 
\longrightarrow & H^2(M_{min}, {\cal E}_{min}|_{M_{min}\times \{t\}}) \\
 \Big\downarrow & & \Big\downarrow \\
 H^2(M, {\cal O}_M) & \longrightarrow & H^2(M, {\cal E}|_{M\times \{t\}})
\end{array}
\]

 It is not hard to show that both the vertical arrows are isomorphisms. 
The equality of the 
 formal base dimensions follows.

   From now on we may assume that $M=M_{min}$ is a minimal elliptic surface
 and ${\cal E}={\cal E}_{min}$.
   We would like to show that for any given array in $H^2(M, {\cal O}_{M})$,
 there exists a pair
 $(t, s), t\in T(M), s\in H^0(M, {\cal E}|_{M\times \{t\}})-\{0\}$ such that
  this array is mapped injectively under $\otimes s: 
H^2(M, {\cal O}_M)\mapsto H^2(M, {\cal E}_{M\times \{t\}})$.

 By Serre duality, it suffices to show that for any given ray in 
$H^0(M, {\cal K}_M)$,
 there exists a pair $(t, s)$ as above such that

 $$(\otimes s)^{\ast}:H^0(M, {\cal K}_M\otimes {\cal E}^{\ast}|_{M
\times \{t\}})\mapsto H^0(M, {\cal K}_M)$$

 maps onto this ray.  Because both $C$ and $c_1(K_M)$ are pulled back 
from the base $\Sigma$ of the
 elliptic fibration $f:M\mapsto 
\Sigma$, the above statement can be translated into showing that: For 
all non-zero sections in
 $H^0(\Sigma, {\cal K})$, ${\cal K}_M=f^{\ast}{\cal K}$,
 $deg_{\Sigma}{\cal K}=g(\Sigma)-1+p_g(M)$, there is an invertible sheaf 
${\cal D}$ over
$\Sigma$, $deg_{\Sigma}{\cal D}=(1-\alpha)deg_{\Sigma}{\cal K}$ and a 
section $\underline{s}$ in
 $H^0(\Sigma, {\cal K}\otimes {\cal D}^{\ast})$ such that the image of 
  $\otimes \underline{s}:H^0(\Sigma, {\cal D})\mapsto H^0(\Sigma, {\cal K})$  
contains these sections in $H^0(\Sigma, {\cal K})=H^0(\Sigma, 
{\cal O}_{\Sigma}(\sum_l m_l x_l))$.

 Let $\sum_l m_l x_l, x_l\in \Sigma$ be an effective divisor in the 
linear system of $|{\cal K}|$, we know that
 $\sum_l m_l=deg_{\Sigma}{\cal K}$. Choose a tuple $n_l\leq m_l, 
\forall l$ such that 
$\sum_l n_l=(1-\alpha) deg_{\Sigma}{\cal K}$. Then we can take 
 ${\cal D}={\cal O}_{\Sigma}(\sum_l n_l x_l)$ and $\sum_l (m_l-n_l)x_l$ 
defines a ray (up to ${\bf C}^{\ast}$ action) 
of sections $\underline{s}\in H^0(\Sigma, 
{\cal O}_{\Sigma}(\sum_l(m_l-n_l)x_l))\cong 
 H^0(\Sigma, {\cal K}\otimes {\cal D}^{\ast})$. It is 
apparent that the map induced by tensoring with
 $\underline{s}$ maps onto the ray of sections defining $\sum_l m_l x_l$. 
$\Box$

\medskip

\noindent (ii).  $-2 fbd(C, M)\leq C^2-c_1(K_M)\cdot C<0$. We do not
 classify $C$ in this situation. Any such
 $C$ does not correspond to basic classes for ordinary Seiberg-Witten 
invariants because of the negativity of its
 Gromov-Taubes dimension. These classes are candidates of the exceptional 
classes.
\label{2}

\medskip

 In the following, we discuss the construction of algebraic Kuranishi model 
of a $C\in H^2(M, {\bf Z})\cap
 H^{1, 1}(M, {\bf C})$ satisfying (i)' or (ii) above. 
In these cases the bundle map of the Kuranishi model will not be algebraic. 

\medskip

\begin{prop}\label{prop; exist}
 Let $M$ be an algebraic surface with $p_g(M)>0$ and let $C\in H^2(M,
 {\bf Z})\cap H^{1,1}(M, {\bf C})$
 be an integral $(1, 1)$ class satisfying $(i)'$ and $(ii)$ on page 
\pageref{1'}, \pageref{2}.
 Then there exist algebraic vector bundles ${\bf V}$,
${\bf W}$ and
 an algebraic bundle map $\Phi_{\bf VW}:{\bf V}\mapsto {\bf W}$ such that

\medskip

\noindent (i). There exist another pair of algebraic bundles $\tilde{\bf W}, 
\tilde{\bf V}$ and the
 bundle map
 $\Psi_{\tilde{\bf W}\tilde{\bf V}}: \tilde{\bf W}\mapsto \tilde{\bf V}$ 
such that
 
$$Coker(\Psi_{\tilde{\cal W}\tilde{\cal V}})
 \cong \bigl({\cal R}^2\pi_{\ast}\bigl({\cal 
 E}\bigr)\bigr)\cong {\cal R}^0\pi_{\ast}\bigl({\cal 
 E}^{\ast}\otimes {\cal K}_M\bigr)^{\ast}, $$ by relative Serre duality. 

\medskip

\noindent (ii). The locally free sheaf ${\cal W}$ associated with the 
vector bundle ${\bf W}$
 is a sub-sheaf of the locally free sheaf 
$\tilde{\cal W}\oplus {\cal O}^{p_g-fbd(C, M)}_{T(M)}$. And we have

$$Ker(\Phi_{{\cal V}{\cal W}})\cong {\cal R}^0\pi_{\ast}\bigl({\cal E}\bigr),$$ and
  $Coker(\Phi_{{\cal V}{\cal W}}$ contains a sub-sheaf isomorphic to 
${\cal R}^1\pi_{\ast}\bigl({\cal E}\bigr)$.

\medskip

\noindent (iii). $rank_{\bf C}{\bf V}-rank_{\bf C}{\bf W}={C^2-c_1(K_M)\cdot
 C\over 2}-q(M)+p_g(M)+fbd(C, M)+1$.
\end{prop}

\medskip 
 
\noindent Proof of the proposition: We prove $(i)$, $(ii)$, and $(iii)$
 step by step.

\medskip
\noindent Step one:
  As usual, we choose an effective ample divisor $D$ on $M$ and consider
 the following short exact sequence for
 a large enough $n$. By choosing $n$ large enough, we may assume 
$\Sigma\in |nD|$ to be a smooth curve in the complete 
linear system $|nD|$.
 In the following, we assume that $\Sigma=nD$ is a smooth curve.

$$0\mapsto {\cal O}_M({\cal E})\mapsto {\cal O}(\Sigma)\otimes {\cal E}
\mapsto {\cal O}_{\Sigma}(\Sigma)\otimes
 {\cal E}\mapsto 0.$$

   We have
 the following long exact sequences and sheaf isomorphisms,

 $$ 0\mapsto {\cal R}^0\pi_{\ast}\bigl( {\cal E} \bigr)\mapsto 
{\cal R}^0\pi_{\ast}\bigl(  
 {\cal O}(\Sigma)\otimes {\cal E}\bigr)\mapsto
{\cal R}^0\pi_{\ast}\bigl( {\cal O}_{\Sigma}(\Sigma)\otimes {\cal E} \bigr)$$
$$\mapsto {\cal R}^1\pi_{\ast}\bigl( {\cal E} \bigr)\mapsto 0,$$

$${\cal R}^1\pi_{\ast}\bigl({\cal O}_{\Sigma}(\Sigma)\otimes {\cal E}\bigr)
\cong {\cal R}^2\pi_{\ast}\bigl({\cal E} \bigr).$$

\medskip

 The difference from the previous situations is that ${\cal R}^0
\pi_{\ast}\bigl( {\cal O}_{\Sigma}
(\Sigma)\otimes {\cal E} \bigr)$
 and ${\cal R}^1\pi_{\ast}\bigl( {\cal O}_{\Sigma}(\Sigma)\otimes 
{\cal E} \bigr)$ may not be locally free over
 $T(M)$. 

 To remedy this, consider a sufficiently very ample divisor $\Delta$
 on the smooth curve 
$\Sigma\in |nD|$ and take the derived long 
exact sequence of the following short exact sequence,

 $$\hskip -.5in
0\mapsto {\cal O}_{\Sigma}(\Sigma)\otimes {\cal E}|_{\Sigma\times T(M)}
\mapsto {\cal O}_{\Sigma}(\Sigma+\Delta)\otimes
 {\cal E}|_{\Sigma\times T(M)}\mapsto 
{\cal O}_{\Delta}(\Sigma+\Delta)\otimes {\cal E}|_{\Sigma\times T(M)}
\mapsto 0.$$ Then we have
 the following four-terms long exact sequence on the derived sheaves 
when $\Delta$ is sufficiently very ample
 on $\Sigma$,

$$0\mapsto {\cal R}^0\pi_{\ast}\bigl({\cal O}_{\Sigma}(\Sigma)\otimes 
{\cal E}|_{\Sigma\times T(M)}\bigr)\mapsto
 {\cal R}^0\pi_{\ast}\bigl({\cal O}_{\Sigma}(\Sigma+\Delta)\otimes 
{\cal E}|_{\Sigma\times T(M)}\bigr)\mapsto$$

$${\cal R}^0\pi_{\ast}\bigl({\cal O}_{\Delta}(\Sigma+\Delta)\otimes 
{\cal E}|_{\Delta\times T(M)}\bigr)\mapsto
 {\cal R}^1\pi_{\ast}\bigl({\cal O}_{\Sigma}(\Sigma)\otimes 
{\cal E}|_{\Sigma\times T(M)}\bigr)\mapsto 0.$$ 

\medskip

Define ${\cal V}$, $\tilde{\cal W}$ and $\tilde{\cal V}$ to be the 
locally free sheaves

${\cal R}^0\pi_{\ast}\bigl({\cal O}(\Sigma)\otimes {\cal E}\bigr)$,
${\cal R}^0\pi_{\ast}\bigl({\cal O}_{\Sigma}(\Sigma+\Delta)\otimes 
{\cal E}|_{\Sigma\times T(M)}\bigr)$
 and ${\cal R}^0\pi_{\ast}\bigl({\cal O}_{\Delta}(\Sigma+\Delta)
\otimes {\cal E}|_{\Delta\times T(M)}\bigr)$
 over $T(M)$, respectively.
 As usual let ${\bf V}$, $\tilde{\bf W}$ and $\tilde{\bf V}$ be the 
corresponding algebraic vector bundles. 
  
 Then there are bundle maps $\Phi'_{{\bf V}\tilde{\bf W}}:{\bf V}\mapsto
 \tilde{\bf W}$, and 
$\Psi_{\tilde{\bf W}\tilde{\bf V}}:\tilde{\bf W}\mapsto \tilde{\bf V}$.

 The map $\Phi'_{{\bf V}\tilde{\bf W}}$ is induced by the sheaf morphism
 composition

$$\hskip -.9in {\cal R}^0\pi_{\ast}\bigl({\cal O}(\Sigma)\otimes 
{\cal E}\bigr)\mapsto 
{\cal R}^0\pi_{\ast}\bigl( {\cal O}_{\Sigma}(\Sigma)\otimes {\cal E}
 \bigr) \mapsto 
{\cal R}^0\pi_{\ast}\bigl({\cal O}_{\Sigma}(\Sigma+\Delta)\otimes 
{\cal E}|_{\Sigma\times T(M)}\bigr),$$

while $\Psi_{\tilde{\bf W}\tilde{\bf V}}$ is induced by the above
 four term long exact sequence.

We have $Ker(\Phi'_{{\bf V}\tilde{\bf W}})\cong {\cal R}^0
\pi_{\ast}\bigl({\cal E}\bigr)$.

 By construction $$Coker(\Psi_{\tilde{\bf W}\tilde{\bf V}})
\cong 
{\cal R}^1\pi_{\ast}\bigl({\cal O}_{\Sigma}(\Sigma)\otimes 
{\cal E}|_{\Sigma\times T(M)}\bigr)$$
$$\cong {\cal R}^2\pi_{\ast}\bigl({\cal E} \bigr)\cong 
{\cal R}^0\pi_{\ast}\bigl({\cal E}^{\ast}\otimes
{\cal K}_M \bigr).$$

 Thus, statement one has been proved.

 In the following, we will implicitly identify 
$Coker(\Psi_{\tilde{\bf W}\tilde{\bf V}})$ with 
 ${\cal R}^2\pi_{\ast}\bigl({\cal E} \bigr)$.

\medskip

\noindent Step two: The definition of the number $fbd(C, M)$ 
implies that there exists a surjective
 sheaf morphism 

$$\hskip -.7in {\cal O}_{T(M)}^{p_g-fbd(C, M)}
\subset {\cal O}_{T(M)}^{p_g}\cong
 {\cal R}^2\pi_{\ast}\bigl( {\cal O}_{M\times 
T(M)}\bigr)\mapsto {\cal R}^2\pi_{\ast}\bigl( {\cal E}\bigr)\mapsto 0.$$

We define the sheaf ${\cal W}$ to be
 the direct sum
$${\cal W}=Ker(\tilde{\cal W}\mapsto \tilde{\cal V})\oplus 
 Ker({\cal O}_{T(M)}^{p_g-fbd(C, M)}\mapsto
{\cal R}^2\pi_{\ast}\bigl( {\cal E}\bigr)).$$

 Thus, ${\cal W}$ is naturally a subsheaf of $\tilde{\cal W}\oplus
 {\cal O}^{p_g-fbd(C, M)}_{T(M)}$.

\begin{lemm}\label{lemm; onto}
  The sheaf ${\cal W}$ is locally free.
\end{lemm}

\medskip

\noindent Proof of the lemma: To show that ${\cal W}$ is locally free,
 it suffices to show that ([Ha], page 174, chapter II. 8.9.) for all
 $x\in T(M)$, the 
 $k(x)$ vector space ${\cal W}\otimes k(x)$ is of constant rank independent of
 $x$.

 But this follows from the equality

$$ rank_{k(x)}({\cal W}\otimes k(x))=
rank_{k(x)}Ker(\tilde{\cal W}\mapsto \tilde{\cal V})\otimes k(x)$$
$$+rank_{k(x)}Ker({\cal O}_{T(M)}^{p_g-fbd(C, M)}\mapsto
{\cal R}^2\pi_{\ast}\bigl( {\cal E}\bigr))\otimes k(x)$$
$$ =rank_{k(x)}\tilde{\cal W}\otimes k(x)-rank_{k(x)}\tilde{\cal V}\otimes k(x)+
rank_{k(x)}{\cal R}^2\pi_{\ast}\bigl( {\cal E}\bigr))\otimes k(x)$$
$$+rank_{k(x)}{\cal O}^{p_g-fbd(C, M)}_{T(M)}\otimes k(x)-
rank_{k(x)}{\cal R}^2\pi_{\ast}\bigl( {\cal E}\bigr))\otimes k(x)$$
$$=rank_{\bf C}\tilde{\cal W}-rank_{\bf C}\tilde{\cal V}+p_g-fbd(C, M),$$

and is independent of $x\in T(M)$
because $\tilde{\cal V}$, $\tilde{\cal W}$ are locally free.
 The lemma has been proved. $\Box$

\medskip

 Once we realize ${\cal W}$ to be a locally free sheaf, we are ready to
 construct $\Phi_{{\cal V}{\cal W}}:{\cal V}\mapsto {\cal W}$.

 Because the image $\Phi'_{{\cal V}\tilde{\cal W}}({\cal V})$ 
lies in the sub-sheaf 
${\cal R}^0\pi_{\ast}\bigl( {\cal O}_{\Sigma}(\Sigma)\otimes 
{\cal E} \bigr)$, the kernel of
 $\Psi_{\tilde{\bf W}\tilde{\bf V}}$, $Im(\Phi'_{{\cal V}
\tilde{\cal W}})$ actually lies inside
 ${\cal W}$. Thus, the map $\Phi'_{{\cal V}\tilde{\cal W}}$ 
factors through a sheaf morphism $\Phi_{{\cal V}{\cal W}}:{\cal V}\mapsto
{\cal W}$ and 
$$Ker(\Phi_{{\cal V}{\cal W}})=Ker(Phi'_{{\cal V}\tilde{\cal W}})\cong
 {\cal R}^0\pi_{\ast}\bigl({\cal E}\bigr).$$

 On the other hand, 

$$Coker(\Phi_{{\cal V}{\cal W}})\supset Ker(\Psi_{\tilde{\bf W}
\tilde{\bf V}})/Im(Phi'_{{\cal V}\tilde{\cal W}})$$
$$={\cal R}^0\pi_{\ast}\bigl({\cal O}_{\Sigma}(\Sigma)\otimes 
{\cal E}\bigr)/Im({\cal R}^0\pi_{\ast}bigl( 
{\cal O}(\Sigma)\otimes {\cal E}\bigr)\cong {\cal R}^1
\pi_{\ast}\bigl( {\cal E}\bigr).$$ 

This ends the proof of the second statement.

\medskip

\noindent Step Three:
We have $$rank_{\bf C}{\bf V}-rank_{\bf C}{\bf W}=rank_{\bf C}{\bf V}-
(rank_{\bf C}\tilde{\bf W}+(p_g-fbd(C, M))-rank_{\bf C}\tilde{\bf V})$$
$$=\chi({\cal O}(\Sigma)\otimes {\cal E})-
(\chi({\cal O}_{\Sigma}(\Sigma+\Delta) \otimes {\cal E})-
deg(\Delta))-p_g+fbd(C, M)$$
$$=\chi({\cal E})-p_g+fbd(C, M)={C^2-c_1(K_M)\cdot 
C\over 2}-q(M)+fbd(C, M)+1, $$ by surface Riemann-Roch theorem.
 This finishes the proof of step three and therefore the proposition. $\Box$

\medskip

   By the algebraic Kuranishi model $\Phi_{{\bf V}{\bf W}}:{\bf V}\mapsto 
{\bf W}$ constructed in 
proposition \ref{prop; exist}, 
the algebraic moduli space of curves dual $C$ can be realized as a
 projective cone determined by the coherent sheaf 
${\cal R}^0\pi_{\ast}\bigl({\cal E}\bigr)$ and is
 embedded inside ${\bf P}_{T(M)}({\bf V})$ as the zero locus of a canonical section
 $s_{\Phi_{{\bf V}{\bf W}}}:\Gamma({\bf P}_{T(M)}({\bf V}), 
{\bf H}_{{\bf P}_{T(M)}}\otimes {\bf W})$.

   Thus,  one defines 
${\cal ASW}(C)$
 to be the intersection product of the algebraic cycle 
classes

 $${\bf c}_1({\bf H})^{\cap {C^2-c_1(K_M)\cdot C\over 2}+
fbd(C, M)}\cap {\bf c}_{top}({\bf H}\otimes {\bf W})
 \in {\cal A}_0({\bf P}_{T(M)}({\bf V})).$$

\medskip

\bigskip

 At the end of this subsection, we explain why the above 
definitions of ${\cal ASW}$ in the various
 cases are independent to the choices of the algebraic Kuranishi models.

\medskip

\begin{prop}\label{prop; independent}
 The algebraic Seiberg-Witten invariants defined in this section
 are independent to the choices of
 algebraic Kuranishi models $({\bf V}, {\bf W}, \Phi_{\bf VW})$ 
used to define them.
\end{prop}

\medskip

\noindent Proof of the proposition: As has been pointed out 
earlier, in a few cases ($p_g=0$) 
the calculation of 
${\cal ASW}$ can be identified with the wall crossing formula
 of Seiberg-Witten invariants 
[LL2] or some
known calculation of Seiberg-Witten invariants [FM2]. 
Thus the answers are known to be 
independent to extra data like the choices of Kuranishi 
models. Nevertheless, we offer an algebraic proof
 for all the different cases.

 The algebraic Seiberg-Witten invariants for $q(M)=0$ 
algebraic surface is $\pm 1$ and the proof of the proposition is
 trivial in these cases. In the following, we assume 
$q(M)>0$ and separate into two different situations.

 \noindent (i). 
The sheaf ${\cal R}^2\pi_{\ast}\bigl({\cal E}\bigr)=0$ over $T(M)$. In this case,
 ${\cal ASW}(C)$ is equal to

 $$\int_{{\bf P}({\bf V})}\pi^{\ast}_{{\bf P}({\bf V})} 
c_1^{{C^2-c_1(K_M)\cdot C\over 2}+p_g}({\bf H})\cdot
 c_{top}({\bf H}\otimes \pi^{\ast}_{{\bf P}({\bf V})}{\bf W}),$$

  for some algebraic Kuranishi model $({\bf V}, {\bf W}, \Phi_{{\bf VW}})$.

  One may evaluate ${\cal ASW}(C)$ directly and find the answer to be 
$\int_{T(M)}c_{q(M)}({\bf W}-{\bf V})=\int_{T(M)}c_{q(M)}({\cal W}-{\cal V})$
 (following the calculation in the wall crossing formula [LL1], [LL2]).

  Then the independence to the algebraic Kuranishi models is due to the equality
 in the $K$ group of coherent sheaves on
 $T(M)$, $${\cal W}-{\cal V}={\cal R}^1\pi_{\ast}\bigl({\cal E}\bigr)-
{\cal R}^0\pi_{\ast}\bigl({\cal E}\bigr).$$

  An alternative way without evaluating ${\cal ASW}$ is by stabilization and 
lemma \ref{lemm; stablization} on page
 \pageref{lemm; stablization}. Suppose that $({\bf V}_a, {\bf W}_a, 
\Phi_{{\bf V_aW_a}})$, 
$({\bf V}_b, {\bf W}_b, \Phi_{{\bf V_bW_b}})$ are two different 
algebraic Kuranishi models of the class $C$.

 Then by lemma \ref{lemm; stablization} it is not hard to 
see that both the original integral expressions can be
 stabilized to 

 $$\int_{{\bf P}({\bf V}_a\oplus {\bf V}_b)}c_1^{{C^2-c_1(K_M)
\cdot C\over 2}+p_g}({\bf H})\cdot
 c_{top}({\bf H}\otimes \pi^{\ast}_{{\bf P}({\bf V}_a\oplus 
{\bf V}_b)}({\bf W}_a\oplus {\bf V}_b))$$

$$=\int_{{\bf P}({\bf V}_a\oplus {\bf V}_b)}c_1^{{C^2-c_1(K_M)
\cdot C\over 2}+p_g}({\bf H})\cdot
 c_{top}({\bf H}\otimes \pi^{\ast}_{{\bf P}({\bf V}_a\oplus 
{\bf V}_b)}({\bf W}_b\oplus {\bf V}_a)).$$
 
by using $[{\bf W}_b-{\bf V}_b]=[{\bf W}_a-{\bf V}_a]$ in the
 $K$ group of vector bundles.

\medskip

\noindent (ii). Suppose that ${\cal R}^2\pi_{\ast}\bigl({\cal
 E}\bigr)\not= 0$ on $T(M)$, we show the independence to
 the Kuranishi models by a stabilization argument.

\medskip

 Case I: $p_g(M)=0$. Denote the support of the coherent 
${\cal R}^2\pi_{\ast}\bigl( {\cal E}\bigr)$ to be $Z_2$.
 Suppose there are two algebraic Kuranishi models 
$({\bf V}_a, {\bf W}_a, \Phi_{\bf V_aW_a})$,
 $({\bf V}_a, {\bf W}_a, \Phi_{\bf V_aW_a})$ over 
$Z_2^c$ 
for ${\cal M}_C$, the identical stabilization argument as 
in (i). works except the Kuranishi models are defined
 over $Z_2^c$ instead of the whole $T(M)$.
 
\medskip

 Case II: $p_g(M)>0$. In this case the ${\cal ASW}(C)$ 
is defined to be an intersection number on the total space of
 an algebraic vector bundle (over a projective space bundle). 

  Suppose that we are given two different algebraic 
Kuranishi models of ${\cal M}_C$, with the corresponding 
algebraic vector bundles given by
 ${\bf V}_1, {\bf W}_1, \tilde{\bf W}_1, \tilde{\bf V}_1$ and
 ${\bf V}_2, {\bf W}_2, \tilde{\bf W}_2, \tilde{\bf V}_2$, 
respectively. 

 Then we have to show that

$${\bf c}_1({\bf H})^{\cap \{{C^2-c_1(K_M)\cdot C\over 2}+fbd(C, M)\}}\cap 
{\bf c}_{top}({\bf H}\otimes {\bf W}_1) 
\in {\cal A}_0({\bf P}_{T(M)}({\bf V}_1))\cong {\bf Z}$$

 and 

$${\bf c}_1({\bf H})^{\cap \{{C^2-c_1(K_M)\cdot C\over 2}+fbd(C, M)\}}\cap 
{\bf c}_{top}({\bf H}\otimes {\bf W}_2) \in 
{\cal A}_0({\bf P}_{T(M)}({\bf V}_2))\cong {\bf Z}$$

 coincide.

By lemma \ref{lemm; stablization}, we may replace the intersection numbers 
of algebraic cycles by 

$${\bf c}_1({\bf H})^{\cap \{{C^2-c_1(K_M)\cdot C\over 2}+fbd(C, M)\}}\cap 
{\bf c}_{top}({\bf H}\otimes ({\bf W}_1\oplus {\bf V}_2))
\in {\cal A}_0({\bf P}_{T(M)}({\bf V}_1\oplus {\bf V}_2))$$

 and 

$${\bf c}_1({\bf H})^{\cap \{{C^2-c_1(K_M)\cdot C\over 2}+fbd(C, M)\}}\cap 
{\bf c}_{top}({\bf H}\otimes ({\bf W}_2\oplus {\bf V}_1))\in 
{\cal A}_0({\bf P}_{T(M)}({\bf V}_1\oplus {\bf V}_2)).$$

 As before, it suffices to show the following equality in the $K$ group of
 coherent sheaves $K_0(T(M))$,

$$[{\cal V}_2] + [{\cal W}_1]=[{\cal V}_1] + [{\cal W}_2]$$

  and it is equivalent to

$$\hskip -.9in [{\cal V}_2]+[\tilde{\cal W}_1]+[{\cal O}_{T(M)}^{p_g-fbd(C, M)}]-
[\tilde{\cal V}_1] = [{\cal V}_1]+[\tilde{\cal W}_2]+
[{\cal O}_{T(M)}^{p_g-fbd(C, M)}]-[\tilde{\cal V}_2]$$

 and is then equivalent to 

$$[\tilde{\cal W}_1]-[{\cal V}_1]-[\tilde{\cal V}_1] = 
[\tilde{\cal W}_2]-[{\cal V}_2]-[\tilde{\cal V}_2].$$

 By using the exact sequence in the proof of the proposition \ref{prop; exist},
 both sides of the yet to be proved equality in $K_0(T(M))$ can be identified 
 with 

$$-{\cal R}^0\pi_{\ast}\bigl({\cal E}\bigr)+
{\cal R}^1\pi_{\ast}\bigl({\cal E}\bigr)-
{\cal R}^2\pi_{\ast}\bigl({\cal E}\bigr),$$ 
therefore, both sides must be equal. $\Box$

\bigskip

\subsubsection{\bf The Relationship of ${\cal ASW}$ with Family Seiberg-Witten
 Invariants} \label{subsection; fsw}

\bigskip

  In this subsection, we outline the relationship between 
${\cal ASW}$ and the family
 Seiberg-Witten invariant.

  For $p_g(M)=0 (b^+_2=1)$ algebraic surfaces, all the classes 
$C\in H^2(M, {\bf Z})$ automatically become
 $(1, 1)$ classes under the Hodge decomposition. And the algebraic
 Seiberg-Witten invariants of $C$
 are equal to the Seiberg-Witten invariants of the $spin^c$ class 
$2C-K_M$ (in additive notation) in the
 chambers of deformations by large Kahler forms. On the other hand,
 for algebraic
 surfaces with positive geometric genera,
 the algebraic Seiberg-Witten invariant of $C$ is defined only for
 $C\in H^2(M, {\bf Z})\cap 
 H^{1, 1}(M, {\bf C})$ and its dimension formula depends on the 
holomorphic invariant 
$0\leq fbd(C, M)\leq p_g$.

Suppose that $D_C$ is an effective curve on $M$ representing the
 cohomology class $C$
 and $s_{D_C}$ is a defining global section of $D_C$. Then 
we 
 analyze the long exact sequence associated to 

$$0\mapsto {\cal O}_M\stackrel{\otimes s_{D_C}}{\longrightarrow} 
{\cal O}_M(D_C)\mapsto {\cal O}_{D_C}(D_C)\mapsto 0,$$

$$0\mapsto H^0(M, {\cal O}_M)\mapsto H^0(M, {\cal O}_M(D_C))\mapsto 
H^0(M, {\cal O}_{D_C}(D_C))$$
$$\mapsto H^1(M, {\cal O}_M)\mapsto H^1(M, {\cal O}_M(D_C))\mapsto 
H^1(M, {\cal O}_{D_C}(D_C))$$
$$\mapsto H^2(M, {\cal O}_M)\mapsto H^2(M, {\cal O}_M(D_C))\mapsto 0.$$

 Let ${\cal T}_M$ and ${\cal T}^{\ast}_M$ denote the tangent and 
cotangent sheaves of $M$.
  As $C\in H^{1,1}(M, {\bf C})\cong H^1(M, {\cal T}^{\ast}_M)=
H^1(M, \Omega^1_M$, 
the restriction of the cup product  pairing
 
$$\cup: H^1(M, {\cal T}_M)\otimes H^1(M, \Omega_M^1)\longrightarrow 
H^2(M, {\cal O}_M)$$

  gives rises to the linear map

 $$H^1(M, {\cal T}_M)\stackrel{\cup [C]}{\longrightarrow} 
H^2(M, {\cal O}_M),$$

 sending the infinitesimal complex deformations of $M$ to the 
infinitesimal deformation of
 Hodge structures.   In the fixed complex vector space 
$H^2(M,{\bf C})$,
 the deformation of the decomposition 
$$H^2(M, {\bf C})=H^{2, 0}(M, {\bf C})\oplus H^{1, 1}(M, {\bf C})
\oplus H^{0, 2}(M, {\bf C})$$
 is a so-called deformation of the Hodge structures. The tangent 
space of infinitesimal deformations of
 $H^{1, 1}(M, {\bf C})\oplus H^{2, 0}(M, {\bf C})\subset H^2(M, {\bf C})$ can be 
 identified with 

$$\hskip -.3in Hom_{\bf C}( H^{1, 1}(M, {\bf C})\oplus H^{2, 0}(M, {\bf C}), 
H^2(M, {\bf C})/H^{1, 1}(M, {\bf C})\oplus
 H^{2, 0}(M, {\bf C}))$$
$$\cong Hom_{\bf C}(H^{1, 1}(M, {\bf C})\oplus H^{2, 0}(M, {\bf C}), 
H^{0, 2}(M, {\bf C})).$$ 

  Then one may interpret the map 
$H^1(M, {\cal T}_M)\stackrel{\cup [C]}{\longrightarrow} H^2(M, {\cal O}_M)$ 
as the 
composition of the infinitesimal period map $$H^1(M, {\cal T}_M)\mapsto 
Hom_{\bf C}(H^{1, 1}(M, {\bf C})\oplus H^{2, 0}(M, {\bf C}), H^{0, 2}(M, 
{\bf C}))$$

 and 

$$\hskip -.3in Hom_{\bf C}(H^{1, 1}(M, {\bf C})\oplus H^{2, 0}(M, {\bf C}),
 H^{0, 2}(M, {\bf C}))
   \longrightarrow Hom_{\bf C}({\bf C}C, H^{0, 2}(M, {\bf C}))\cong H^2(M, 
{\cal O}_M)$$ through the
 embedding ${\bf C}C\subset H^{1, 1}(M, {\bf C})$.  

   In general, the linear map $\cup [C]: H^1(M, {\cal T}_M)\mapsto H^2(M, 
{\cal O}_M)$ may not be
  surjective. For example, when $C=\alpha\cdot c_1(K_M)$ and when $H^1(M, 
{\cal T}_M)$ is un-obstructed, 
 the map $\cup [C]$ is always  trivial as the canonical class $c_1(K_M)$ 
persists to be of type $(1, 1)$ under
 complex deformations of $M$.

  Given an invertible sheaf ${\cal L}$, a connection on ${\cal L}$ is a 
compatible family 
 of 1st order differential operators $\nabla_{\cdot}$ 
 such that for $U$ open, $X\in \Gamma(U, {\cal T}_M)$ and $s\in 
\Gamma(U, {\cal L})$, $\nabla_X(s)\in
 \Gamma(U, {\cal L})$
satisfies $\nabla_X(f\cdot s)=X(f)\cdot s+f\cdot \nabla_X(s)$ for 
$f\in \Gamma(U, {\cal O}_M)$.

  Suppose $D_C$ is an effective divisor and $s_{D_C}$ is a defining
 global section of $D_C$ in
  $\Gamma(M, {\cal O}(D_C))$. Then $\nabla_{\cdot}(s)|_{D_C}$ 
establishes a morphism
 from $H^1(M, {\cal T}_M)$ to $H^1(D_C, {\cal O}_{D_C}(D_C))$.
  
  Then we have the following commutative diagram

\[
\begin{array}{ccc}
 H^1(M, {\cal T}_M) & \longrightarrow & H^1(D_C, {\cal O}_{D_C}(D_C))  \\
 \Big\downarrow &   & \Big\downarrow\vcenter{%
  \rlap{$\delta$}}   \\
 H^2(M, {\cal O}_M) &  = & H^2(M, {\cal O}_M)
\end{array}
\]

where $\delta$ is the connection homomorphism in the long exact sequence
 and the left vertical
 arrow is the infinitesimal period map.

We can make the following conclusion on the comparison between
 ${\cal ASW}_{pt}(1, C)$ and
 the family Seiberg-Witten theory.

\medskip
\noindent (i). Suppose that $fbd(C, M)=p_g(M)$, then the expected
 dimension of the moduli space is
 ${C^2-c_1(K_M)\cdot C\over 2}+p_g$. Let $B^{p_g}(\epsilon)$ denote
 the radius $\epsilon$ ball
 in the complex space ${\bf C}^{p_g}$.
 
Suppose that there exists a germ of deformation of complex
 structures of $M$, $\pi: {\cal X}\mapsto B^{p_g}(\epsilon)\subset
 {\bf C}^{p_g}$ such that

\noindent (a). $\pi^{-1}(0)$ is bi-holomorphic to $M$. 

and 

\medskip

\noindent (b). the infinitesimal composite period map 

$${\bf T}_0B^{p_g}\mapsto  H^1(M, {\cal T}_M)\stackrel{\cup [C]}{
\longrightarrow} H^2(M, {\cal O}_M)$$ 

 is an isomorphism.

 Then $C$ fails to be a $(1, 1)$ class in the nearby fibers.

 In the given local family of complex manifolds, the family moduli
 space of curves in $C$ localizes to be above
 $0\in B^{p_g}(\epsilon)$. In such situations, the dimension 
formula of ordinary $SW$ theory and the expected 
virtual dimension differ by $p_g$. But the latter matches up 
with the family $SW$ theory dimension formula of a
 $p_g$ family. 
Then ${\cal ASW}(C)$ is actually equal to the (local) family
 Seiberg-Witten invariant of the $spin^c$ class $2C-K_M$ (in 
additive notation) with  the family
 Seiberg-Witten equation deformed by the 
 large fiberwise Kahler forms.

\begin{rem}\label{rem; family}
The local condition (a). and (b). hold when there are $p_g$ 
dimensional 
un-obstructed infinitesimal complex deformation
 of $M$ map injectively into the infinitesimal deformation of Hodge structures.

 The global version asks for the existence of $\pi:{\cal X}\mapsto B$, 
where $B$ is a compact
 smooth $p_g$ dimensional variety and the total space ${\cal X}$ is Kahler. 

 (a).  $\pi^{-1}(b)$ is bi-holomorphic to $M$  for some $b\in B$. 
\label{condition}

 (a)'. $C\in H^2({\cal X}, {\bf Z})$ and $i_b^{\ast}C\in H^{1, 1}(\pi^{-1}(b),
 {\bf C})$ for the
 inclusion map 
 $i_b: \pi^{-1}(b)\mapsto {\cal X}$.
 
 (b).  If $C\in H^{1, 1}(\pi^{-1}(b), {\bf C})$ for any $b\in B$, then the
 infinitesimal composite period map 

$${\bf T}_bB\mapsto H^1(\pi^{-1}(b), {\cal T}_{\pi^{-1}(b)}) 
\mapsto H^2(\pi^{-1}(b),
 {\cal O}_{\pi^{-1}(b)})$$ is an isomorphism.

 It turns out that the algebraic surfaces with trivial canonical bundles 
${\cal K}_M\cong {\cal O}_M$, i.e. K3 surfaces or abelian
 surfaces, carry the so-called twistor families of complex structures 
 which have the desired global properties 
 (a), (a)', (b). In these situations, the ${\cal ASW}$ can be interpreted
 as the family
 Seiberg-Witten invariants in the chambers deformed by large ${\bf S}^2$ 
family of hyperkahler forms [LL1].
  
  Besides hyperkahler ${\bf S}^2$ families of complex structures, it 
is probably too strong
 to find the germs of  smooth families ${\cal X}\mapsto B$ with smooth 
fibers which satisfy all (a), (a)' and
 (b).
\end{rem}
\medskip
  Suppose that $c_1(M)=0$, then all the $C\in H^2(M, 
{\bf Z})$, $C\cdot \omega_M>0$ satisfy
 $fbd(C, M)=p_g$. On the other hand let ${\bf E}$ be sufficiently
 relative very ample on $M\times T(M)\mapsto T(M)$, 
Kodaira vanishing theorem
 also implies that $fbd(C, M)=p_g$ for such classes $c_1({\bf E})=C$. 
In particular, the
 algebraic Seiberg-Witten invariant of $C$ for such ${\bf E}$ can be thought 
formally as the
 family invariant on a germ of $p_g$ dimensional infinitesimal family.

\medskip
\noindent (ii). Suppose that $fbd(C, M)=0$, e.g. $C=\sum E_i$, then the dimension 
formula 
 ${C^2-C\cdot c_1(K_M)\over 2}$ is identical to the $p_g=0$ case.

\medskip

 By combining with the family blowing up formula, the independence of dimension 
formula on $p_g$ 
explains why type $I$ exceptional curves
 in the universal family [L1] of $p_g>0$ algebraic surfaces obey the same 
dimension formula in the $p_g=0$ 
algebraic surfaces.

\medskip
\noindent (iii). Suppose that $0<fbd(C, M)<p_g$, then the dimension formula 
of $C$ is shifted to
 ${C^2-c_1(K_M)\cdot C\over 2}+fbd(C, M)$.
\medskip
 We discuss $(ii)$ and $(iii)$ together.

  By the definition of the number $fbd(C, M)$, there exists a subspace ${\bf F}$ 
of $H^2(M, {\cal O}_M)$ of dimension 
 $fbd(C, M)$ which is the intersection of all the kernels of $H^2(M, 
{\cal O}_M)\stackrel{\otimes s_t}
{\longrightarrow}
 H^2(M, {\cal E}_{M\times \{t\}})$.  Then the subspace ${\bf F}$ defines a 
$fbd(C, M)$ dimensional 
subspace of infinitesimal
 deformation of Hodge structures deforming $C$ away from being 
 a $(1, 1)$ class.
  For positive $fbd(C,M)$, the ${\bf F}$ defines a trivial factor on the 
obstruction bundle of the
 Kahler Seiberg-Witten theory, which causes the usual $SW$ invariant to 
vanish.
 By removing the trivial factor ${\bf F}$, the intersection number defining 
${\cal ASW}(C)$ are generically nonzero.
 Only when certain conditions like $(a)$, $(a)'$, $(b)$ (on page 
\pageref{condition}) hold, the operation of removing 
 a trivial factor in the obstruction bundle can be interpreted as extending 
$M$ into a local family.

\medskip

\bigskip

\section{\bf The Algebraic Proof of the Family Blowup Formula}\label{section; ap}

\medskip

  Having defined the algebraic Seiberg-Witten invariants in section
 \ref{section; define},
 we offer an algebraic proof of the family blowup 
formula which includes the important special cases of the universal
 families $M_{n+1}\mapsto M_n$, $n\in {\bf N}$
 and
 its sub-families $M_{n+1}\times_{M_n}Y(\Gamma)\mapsto Y(\Gamma)$ on 
the closures of admissible strata See [Liu1]).

 Suppose that $\pi:{\cal X}\mapsto B$ is an algebraic family of smooth 
algebraic surfaces over a smooth 
algebraic base manifold $B$ and let $s:B\mapsto {\cal X}$
 denote an algebraic cross section of the fibration.
 Let $C$ be an element in $H^{1,1}({\cal X}, {\bf Z})$ which restricts 
to a monodromy invariant class on the fibers.
  The inclusion $s(B)\subset {\cal X}$ is a codimension two smooth sub-variety 
of ${\cal X}$ with
 normal bundle ${\bf N}_{s(B)}{\cal X}$. Blowing up $s(B)\subset {\cal X}$ 
produces the blown up variety
 ${\cal X}'$ with an exceptional divisor $E_B \cong {\bf P}_B({\bf N})$. We 
denote the
 fiberwise exceptional curve (and the corresponding cohomology class) of 
${\cal X}'\mapsto B$ by the
 notation $E$. Then the blowup formula is expected to relate the algebraic
 family invariants of $C$ over
 ${\cal X}\mapsto B$ and $C'=C+mE, m\in {\bf Z}$ over ${\cal X}'\mapsto B$.

\medskip
  Given the family $\pi:{\cal X}\mapsto B$,  the family of the relative 
Picard group $Pic^0$ associated to 
 the fiber algebraic surfaces as in proposition \ref{prop; torus} forms a 
fiber bundle of complex tori over $B$,
 denoted by ${\cal T}_B({\cal X})$. The fibration ${\cal T}_B({\cal X})\mapsto 
B$ can be identified with the
 quotient
 ${\cal R}^1\pi_{\ast}\bigl({\cal O}_{\cal X}\bigr)/{\cal R}^1\pi_{\ast}
\bigl({\bf Z}\bigr)$.

 As in the $B=pt$ case, there is a holomorphic line bundle ${\bf E}$ over
 ${\cal T}_B({\cal X})$ 
with first Chern class (the pull-back of) $C$. When ${\cal X}\mapsto B$ 
is the universal family 
 $M_{n+1}\mapsto M_n$ or its sub-families, it is easy to use the following
 proposition, prop. \ref{prop; same} repeatedly
 to show that ${\cal T}_B({\cal X})$ is isomorphic to
the trivial product $T(M)\times B$.

\medskip
\begin{prop}\label{prop; same}
Let ${\cal X}'\mapsto B$ denote the blown up fibration from ${\cal X}
\mapsto B$ along $s:B\mapsto {\cal X}$.
Then ${\cal T}_B({\cal X}')={\cal T}_B({\cal X})$. Namely, the torus 
fibration is invariant under the
 blowing up along a cross section.
\end{prop}

\medskip

\noindent Proof: The fibers of ${\cal T}_B({\cal X})$ are constructed from 
the Hodge group $H^{0,1}$ of the
 fibers of ${\cal X}\mapsto B$. The $(0,1)$ component of the Hodge 
decomposition of an algebraic surface 
is invariant under
 blowing ups. By applying the observation to a family 
that is why ${\cal T}_B({\cal X})$ is invariant under blowing ups 
along cross sections.
$\Box$

\bigskip

We introduce some technical conditions on the fibration which 
guarantees the existence of 
algebraic Kuranishi models.

\medskip

\begin{defin}\label{defin; relatively}
 An algebraic fibration $X\mapsto B$ of relative dimension two is 
said to be relatively good if there exists an
 effective very ample divisor $D$ in $X$ which is of relative 
dimension one under $X\mapsto B$.
 The algebraic fibration $X\mapsto B$ is said to be two-relatively 
good if there exist two
 effective very ample divisors $D_1, D_2$ such that 

\noindent (i). $D_1$ and $D_2$ are both smooth and are of relative 
dimension one over $B$.

\noindent (ii). $D_1\cap D_2\subset D_i, i=1, 2$ is a smooth divisor
 in $D_1$ and $D_2$ 
 and is of relative dimension 
 zero over $B$.
\end{defin}

\medskip

We introduce the following definition of formal excess base dimension 
$febd(C, {\cal X}/B)$, 
extending definition \ref{defin; r} for $B=pt$,

\medskip

\begin{defin} \label{defin; extend r}
 Let $\pi:{\cal X}\mapsto B$ be an algebraic fiber bundle of algebraic 
surfaces and 
let $C\in H^{1,1}({\cal X}, {\bf C})\cap H^2({\cal X}, {\bf Z})$
 be a monodromy invariant fiberwise cohomology class. Let ${\cal E}$ 
be the invertible 
sheaf over ${\cal T}_B({\cal X})$
  and consider the following natural pairing

 $${\cal R}^0\pi_{\ast}({\cal E})\otimes 
{\cal R}^2\pi_{\ast}({\cal O}_{\cal X})\mapsto
  {\cal R}^2\pi_{\ast}({\cal E}),$$

  Define $Ann({\cal R}^0\pi_{\ast}({\cal E}))\subset 
{\cal R}^2\pi_{\ast}({\cal O}_{\cal X})$ 
to be the annihilator
 of ${\cal R}^0\pi_{\ast}({\cal E})$ under the pairing.

  If ${\cal R}^2\pi_{\ast}({\cal E})=0$, define $febd(C,
 {\cal X}/B)=p_g$ in this
 trivial case.

 If ${\cal R}^2\pi_{\ast}({\cal E})\not=0$,
 then the formal excess base dimension $febd(C, {\cal X}/B)$ 
is defined to be the rank of the maximal trivial locally free 
 subsheaves of $Ann({\cal R}^0\pi_{\ast}({\cal E}))$.
\end{defin}

\medskip

 To simplify our notations, in the following theorem and its proof 
we do not write explicitly the pull-back maps on line bundles and 
the cohomologies
 $H^{\ast}({\cal X}, {\bf Z})\mapsto H^{\ast}({\cal X}', {\bf Z})$ 
induced by the blowing down map
 ${\cal X}'\mapsto {\cal X}$. Thus, we use the same symbol $C$ to 
denote the fiberwise class on ${\cal X}$
 and its pull-back to ${\cal X}'$. The readers should be able to 
determine from the context whether we refer to
 the class $C$ on ${\cal X}$ or its pull-back to ${\cal X}'$.

 Let ${\cal E}$ be an invertible sheaf on 
${\cal X}'\times_B{\cal T}_B({\cal X}')$,
 then $c_1({\cal E})$ pulls back to a class in 
$H^2({\cal X}', {\bf Z})$. We abuse the
 notation and denote it by the same symbol
 $c_1({\cal E})$. In the following, we assume 
  that the pull-back first Chern class 
 $c_1({\cal E})=C$ or $C+mE$ in $H^2({\cal X}', {\bf Z})$.

\medskip

If the geometric genus $p_g$ of the fiber algebraic surfaces is zero or if 
$c_1(K_{{\cal X}/B})-C$ is of 
non-positive degree (with respect to an ample polarization), we assume that 
fiber bundle 
$\pi':{\cal X}'\mapsto B$ is relatively good. If $p_g$ of the fiber surfaces 
is greater than zero and 
 the relative degree of the class $c_1(K_{{\cal X}/B})-C$ is positive, then 
we assume that 
 $\pi':{\cal X}'\mapsto B$ is two-relatively good.

 Once we make such assumptions on the fiber bundle ${\cal X}'\mapsto B$, 
we can mimic the construction in
 section \ref{section; define} and construct the algebraic family Kuranishi
 models of the family 
moduli spaces of curves dual to $C$ and $C+mE$. 
 Because the current construction is basically the family extension of
 the construction for 
$B=pt$, we omit much of the details.
 Nevertheless, we address the analogue of proposition \ref{prop; exist}, 
which constructs the
 algebraic family Kuranishi model $({\bf V}, {\bf W}, \Phi_{{\bf V}{\bf W}})$ 
of the case $p_g>0$ and 
${\cal R}^2\pi_{\ast}\bigl({\cal E}\bigr)\not=0$.

\medskip

By definition \ref{defin; extend r} of $febd(c_1({\cal E}), 
{\cal X}'/B)$, there exists an 
inclusion
${\cal O}^{febd(c_1({\cal E}), {\cal X}'/B)}_{{\cal T}_({\cal X}')}\subset 
{\cal R}^2\pi'_{\ast}\bigl({\cal O}_{{\cal X}'}\bigr)$.
Define the quotient locally free sheaf 
${\cal R}^2\pi'_{\ast}\bigl({\cal O}_{{\cal X}'}\bigr)/
{\cal O}^{febd(c_1({\cal E}),
 {\cal X}'/B)}_{{\cal T}_B({\cal X}')}$
 to be ${\cal F}$. Denote $i_{\cal F}$ to be the surjective
 morphism ${\cal F}\mapsto {\cal R}^2\pi'_{\ast}\bigl({\cal E}\bigr)$.

\medskip
 
\begin{prop}\label{prop; analogue}
 Suppose that the algebraic fiber bundle $\pi: {\cal X}'\mapsto B$ 
is two-relatively good, 
 then for the given invertible ${\cal E}$ over ${\cal X}'\times_B 
 {\cal T}_B({\cal X}')$,

\noindent (i). There exists a pair of locally free sheaves 
 $\tilde{\cal W}, \tilde{\cal V}$ and a sheaf morphism 
 $\Psi_{\tilde{\cal W}\tilde{\cal V}}:\tilde{\cal W}\mapsto 
\tilde{\cal V}$
 such that $Coker(\Psi_{\tilde{\cal W}\tilde{\cal V}})\cong 
{\cal R}^2\pi'_{\ast}\bigl({\cal E}\bigr)$.

\medskip

\noindent (ii). The locally free sheaf ${\cal W}$ associated to ${\bf W}$ 
is taken to be
$Ker(\Psi_{\tilde{\cal W}\tilde{\cal V}})\oplus Ker(i_{\cal F})$ 
and there exists the algebraic family
 Kuranishi morphism $\Phi_{{\cal V}{\cal W}}:{\cal V}\mapsto {\cal W}$ such
 that
 $Ker(\Phi_{{\cal V}{\cal W}})\cong {\cal R}^0\pi'_{\ast}\bigl({\cal E}\bigr)$
 and $Coker(\Phi_{{\cal V}{\cal W}})$
 contains ${\cal R}^0\pi'_{\ast}\bigl({\cal E}\bigr)$ as a sub-sheaf.

\medskip

\noindent (iii). There is the following formula on the virtual rank,

 $$rank_{\bf C}{\bf V}-rank_{\bf C}{\bf W}=
{c_1^2({\cal E})-c_1({\cal E})\cdot c_1({\bf K}_{{\cal X}'/B})\over 2}
-q+febd(c_1({\cal E}), {\cal X}'/B)+1.$$
\end{prop}

\medskip

\noindent Proof of the proposition: 

\noindent Step One: By using the condition of $\pi':{\cal X}'\mapsto B$ being 
``2-relatively good'', the proof of the statement (i) is parallel to the
 $B=pt$ case.
 
As before, we use bold character like ${\bf W}$ to denote the 
algebraic vector bundle and
 use the calligraphic character like ${\cal W}$ to denote the 
corresponding sheaf of sections.

  We set $D_1=D$ and $lD_2=\Delta$ and take $n, l$ to be 
sufficiently large integers.
  We take ${\cal V}$, $\tilde{\cal V}$, $\tilde{\cal W}$ to be 
  ${\cal R}^0\pi'_{\ast}\bigl({\cal O}_{{\cal X}'}(nD)\otimes {\cal E}\bigr)$,
 ${\cal R}^0\pi'_{\ast}\bigl( {\cal O}_{nD\cap \Delta}(nD+\Delta)\otimes 
{\cal E}\bigr)$
 and ${\cal R}^0\pi'_{\ast}\bigl( {\cal O}_{nD}(nD+\Delta)\otimes 
{\cal E}\bigr)$, respectively.
 
  Then we have for sufficiently large $n$ and $l$, 

 $$\hskip -.9in 0\mapsto {\cal R}^0\pi'_{\ast}\bigl({\cal E}\bigr)\mapsto 
{\cal V}\mapsto 
 {\cal R}^0\pi'_{\ast}\bigl({\cal O}_{nD}(nD)\otimes{\cal E}\bigr)
\mapsto {\cal R}^1\pi_{\ast}\bigl({\cal E}\bigr)\mapsto 0,$$

 $$\hskip -.9in 0\mapsto {\cal R}^0\pi'_{\ast}\bigl({\cal O}_{nD}(nD)
\otimes{\cal E}\bigr)
\mapsto \tilde{\cal W}\mapsto \tilde{\cal V}\mapsto 
{\cal R}^1\pi'_{\ast}\bigl({\cal O}_{nD}(nD)\otimes{\cal E}\bigr)\mapsto 0,$$

 and $${\cal R}^1\pi'_{\ast}\bigl({\cal O}_{nD}(nD)\otimes{\cal E}\bigr)
\cong {\cal R}^2\pi'_{\ast}\bigl(
{\cal E}\bigr),$$

 similar to the $B=pt$ case.

\medskip

\noindent Step Two: For statement (ii)., we focus on the locally free-ness of 
${\cal W}=Ker(\Psi_{\tilde{\cal W}\tilde{\cal V}})\oplus Ker(i_{\cal F})$. 
 According to [Ha] page 174, chapter II lemma 8.9, 
it suffices to check that for all $y\in {\cal T}_B({\cal X}')$,
 $dim_{k(y)}{\cal W}\otimes k(y)$ is a constant independent of $y$.

We demonstrate this by showing (using the exact sequences listed 
in step one) 

$$dim_{k(y)}{\cal W}\otimes k(y)=dim_{k(y)}Ker(\Psi_{\tilde{\cal W}
\tilde{\cal V}})\otimes k(y)+
dim_{k(y)}Ker(i_{\cal F})\otimes k(y)$$

$$=\bigl( dim_{k(y)}\tilde{\cal W}\otimes k(y)-dim_{k(y)}\tilde{\cal V}
\otimes k(y)+
dim_{k(y)}{\cal R}^2\pi_{\ast}\bigl({\cal E}\bigr)\otimes k(y)\bigr)+$$

$$\bigl(dim_{k(y)}{\cal F}\otimes k(y)-dim_{k(y)}{\cal R}^2\pi_{\ast}
\bigl({\cal E}\bigr)\otimes k(y)\bigr)=
  rank_{\bf C}\tilde{\bf W}-rank_{\bf C}\tilde{\bf V}+rank_{\bf C}{\bf F}.$$

 Thus, the sheaf ${\cal W}$ is locally free.
 The remaining conclusions about $Ker$ or $Coker$ are
 identical to the original argument in prop. \ref{prop; exist}
 and we omit it here.

\medskip

\noindent Step Three: 
For statement (iii). in the proposition, we notice that 
 the construction of $\tilde{\cal V}, \tilde{\cal W}$ and ${\cal V}$ are
 the
 family extension of the construction in prop \ref{prop; exist} and their 
ranks are independent to 
 this extension, we still have as before

  $$rank_{\bf C}{\bf V}+rank_{\bf C}\tilde{\bf V}-rank_{\bf C}\tilde{\bf W}=
{c_1({\cal E})^2-c_1({\cal E})\cdot c_1({\bf K}_{{\cal X}'/B})\over 2}-q+p_g+1$$ 

 by surface Riemann-Roch formula.

 Then the formula in (iii). can be derived by using the calculation in the step two and the
 formula $rank_{\bf C}{\bf F}=p_g-febd(c_1({\cal E}), {\cal X}'/B)$.
$\Box$

\medskip

\begin{mem}\label{mem; K}
 The rank calculation in the step two of the proof can be strengthened to
 imply that
 $$[{\cal W}]=[\tilde{\cal W}]-[\tilde{\cal V}]+[{\cal F}]$$
in the $K$ group of coherent sheaves on ${\cal T}_B({\cal X}')$,
 $K_0({\cal T}_B({\cal X}'))$.
 This identity will be used implicitly in the derivation of the 
 algebraic family blowup formula.
\end{mem}

\medskip

 By using the data $({\bf V}, {\bf W}, \Phi_{{\bf V}{\bf W}})$ of an 
algebraic family Kuranishi model, we 
define the algebraic family Seiberg-Witten invariants similar to the $B=pt$ case.
  The corresponding
 mixed invariants for $C+mE=c_1({\cal E})$ will be defined to be  
$$\hskip -.9in {\cal AFSW}_{{\cal X}'\mapsto B}(c, C+mE)=
c_1^{{c_1^2({\cal E})-c_1({\cal E})
\cdot c_1({\bf K}_{{\cal X}'/B})\over 2}+febd(c_1({\cal E}), 
{\cal X}'/B)}({\bf H})\cap
c_{top}({\bf H}\otimes {\bf W})\in 
{\cal A}_0({\bf P}_{{\cal T}_B({\cal X}')}({\bf V})),$$ 
 
 for $c\in {\cal A}_{\ast}(B)$.
 
\medskip
\begin{theo}\label{theo; al-bl} 
 Let ${\cal X}'\mapsto B$ and ${\cal X}\mapsto B$ and $C$ be as 
described above with the
 appropriated relatively good conditions, and let $m$ be an integer such that
 family dimension 
$\int_{{\cal X}/B}{(C^2-c_1(K_1{{\cal X}/B}))\cdot C\over 2}+febd(C, 
{\cal X}/B)-
{m^2-m\over 2}+dim_{\bf C}B\geq 0$.  
Then for any $c\in H^{\ast}(B, {\bf Z})$ (or any class $c$ of algebraic
 cycles in $B$)
the mixed algebraic family Seiberg-Witten invariant of the class $C+mE$ 
over the blown up fibration ${\cal X}'\mapsto B$ is 
 related to the algebraic family Seiberg-Witten invariant of $C$ by the 
following formula, 

$${\cal AFSW}_{{\cal X}'\mapsto B}(c, C+mE)=\sum_{i\geq 0}
{\cal ASW}_{{\cal X}\mapsto B}(
c\cap c_i({\bf E}\otimes ({\bf S}^{m-2}({\bf C}_B\oplus 
{\bf N}_{s(B)}{\cal X}))), C) $$

 for $m\geq 1$ and 

$${\cal AFSW}_{{\cal X}'\mapsto B}(c, C+mE)=
\sum_{i\geq 0}{\cal ASW}_{{\cal X}\mapsto B}(
c\cap c_i({\bf E}\otimes ({\bf S}^{-m-1}({\bf C}_B\oplus
 {\bf N}^{\ast}_{s(B)}{\cal X}))), C) $$

 for $m\leq 0$.

\end{theo}
\medskip

\begin{rem}\label{rem; section}
  One may assume additionally that ${\cal X}'\mapsto B$ carries 
 cross sections. Then the image of each cross section
  to ${\cal X}'$ under the blowing down map ${\cal X}'\mapsto {\cal X}$ 
induces a cross section on ${\cal X}\mapsto B$.
  By using the cross 
sections, one may interpret formally the algebraic family 
Seiberg-Witten invariant as a virtual count of
 holomorphic curves through generic cross sections of ${\cal X}'
\mapsto B$, ${\cal X}\mapsto B$ and the
 family blowup formula relates these invariants.
\end{rem}

\medskip
\noindent Proof of the theorem:
 Suppose that abstractly
 we are given an algebraic Kuranishi model of the fiberwise class $C$ 
over ${\cal X}\mapsto B$. 
Namely, there is a bundle map $\Phi_{\bf VW}: {\bf V}\longrightarrow 
{\bf W}$ over ${\cal T}_B({\cal X})$ 
and the 
 algebraic family invariant is defined either as an intersection number 
involving 
 $c_{top}({\bf H}\otimes {\bf W})$ over the whole projective space bundle
${\bf P}_{{\cal T}_B({\cal X})}$ or over an open subset of it. As usual, 
we use ${\cal V}$ and ${\cal W}$ to
denote the corresponding locally free sheaves associated to ${\bf V}, 
{\bf W}$.

  Let ${\cal E}$ and ${\cal E}_m$ denote the invertible sheaves associated 
with $C$ and $C+mE$ on ${\cal X}'$, respectively.
  By proposition \ref{prop; same}, we can identify their base spaces to be 
${\cal T}_B({\cal X}')={\cal T}_B({\cal X})$. 

\medskip

 \noindent Case One:
 We assume that $m>0$, then there is an short exact sequence relating 
${\cal E}_m$ and ${\cal E}$,

$$0\mapsto {\cal E}\mapsto {\cal E}_m\mapsto {\cal O}_{mE_B}({mE_B})
\otimes {\cal E}\mapsto 0,$$

 where $E_B={\bf P}_B({\bf N}_{s(B)}{\cal X})$ 
denotes the exceptional divisor of the blowing down map ${\cal X}'
\mapsto {\cal X}$.

 Firstly, we establish the following lemma,

\medskip

\begin{lemm}\label{lemm; projective}
Let $E_B={\bf P}_B({\bf N}_{s(B)}{\cal X})$ denote the 
exceptional divisor of 
 ${\cal X}'\mapsto {\cal X}$. Then for all $m>0$, the
 sheaf ${\cal R}^0\pi'_{\ast}\bigl({\cal O}_{mE_B}(mE_B)
\otimes {\cal E}\bigr)$ over ${\cal T}_B({\cal X})$  
is the zero sheaf and
 the sheaf ${\cal R}^0\pi'_{\ast}\bigl({\cal O}_{mE_B}(mE_B)
\otimes{\cal E}\bigr)$  is
 locally free.
\end{lemm}

\medskip

\noindent Proof: 
\medskip
  By induction we assume that for $k=m-1$ the sheaf   
${\cal R}^0\pi'_{\ast}\bigl({\cal O}_{kE_B}(kE_B)\bigr)$ over has
 been proved to be the zero sheaf. We prove that
 for $k+1=m$, ${\cal R}^0\pi'_{\ast}\bigl({\cal O}_{(k+1)B}((k+1)E_B)\bigr)$ 
is also trivial.

Recall that for two effective divisors $A$ and $B$ on an algebraic manifold, 
we have the
 following short exact sheaf sequence (see for example [Fr], pp 15, equation 
(1.9)),

$$0\mapsto {\cal O}_A(A)\mapsto {\cal O}_{A+B}(A+B)\mapsto {\cal O}_B(A+B)
\mapsto 0.$$ 

  In our case we take $A=kE_B$ and $B=E_B$, then 

$$0\mapsto {\cal O}_{kE_B}(kE_B)\mapsto {\cal O}_{(k+1)E_B}((k+1)E_B)
\mapsto {\cal O}_{E_B}((k+1)E_B)\mapsto 0.$$

Then we have the following (portion of) 
derived long exact sequence by pushing forward along $\pi':{\cal X}'\mapsto B$,

$$0\mapsto {\cal R}^0\pi'_{\ast}\bigl({\cal O}_{kE_B}(kE_B)\bigr)\mapsto 
{\cal R}^0\pi'_{\ast}\bigl({\cal O}_{(k+1)E_B}((k+1)E_B)\bigr)\mapsto 
{\cal R}^0\pi'_{\ast}\bigl({\cal O}_{E_B}((k+1)E_B)\bigr).$$

  It suffices to 
argue that ${\cal R}^0\pi'_{\ast}\bigl({\cal O}_{E_B}((k+1)E_B)\bigr)$ 
vanishes.
  By base change to all closed points $b\in B$, it suffices to check that 
$h^0({\bf P}^1, {\cal O}_{{\bf P}^1}(-k-1))=0$ for $k\geq 0$. It 
follows from the negativity of the degree
 $deg {\cal O}_{{\bf P}^1}(-k-1)$ on  ${\bf P}^1$. Then by 
curve-Riemann-Roch and Grauert
criterion (see [Ha] page 288 cor. 12.9) the sheaf 
${\cal R}^1\pi'_{\ast}\bigl({\cal O}_{(k+1)E_B}((k+1)E_B)\bigr)$
 is locally free.

 By tensoring with ${\cal E}$ from the base ${\cal T}_B({\cal X})$ 
and by moving ${\cal E}$ into the right derived image functor we 
find that
 $ {\cal R}^1\pi'_{\ast}\bigl({\cal O}_{(k+1)E_B}((k+1)E_B)\otimes 
{\cal E}\bigr)$ is locally free and
 ${\cal R}^0\pi'_{\ast}\bigl({\cal O}_{(k+1)E_B}((k+1)E_B)\otimes 
{\cal E}\bigr)$ is the zero sheaf.
 \medskip

 Moreover, we have the following short exact sequence of right 
derived sheaves, which will be used in
 identifying the relative obstruction bundles appearing in the 
statement of the theorem,

$$\hskip -.9in 
0\mapsto {\cal R}^1\pi'_{\ast}\bigl({\cal O}_{(m-1)E_B}((m-1)E_B)\bigr)\mapsto 
 {\cal R}^1\pi'_{\ast}\bigl({\cal O}_{mE_B}(mE_B)\bigr)\mapsto
 {\cal R}^1\pi'_{\ast}\bigl({\cal O}_{E_B}(mE_B)\bigr)\mapsto 0.$$

   By using the fact ${\cal R}^0\pi'_{\ast}
\bigl({\cal O}_{{\bf P}_B({\bf N}_{s(B)}({\cal X}))}(-E_B)\bigr)=
{\cal N}^{\ast}_{s(B)}{\cal X}$, the conormal 
sheaf of $s(B)\subset {\cal X}$ and the relative Serre duality
 along ${\bf P}_B({\bf N}_{s(B)}{\cal X})\mapsto B$, one can identify 
${\cal R}^1\pi'_{\ast}\bigl({\cal O}_{E_B}(mE_B)\bigr)$ with 
the $m-2$ symmetric power
  ${\bf S}^{m-2}({\cal N}_{s(B)}{\cal X})$.

Thus, there is an identity among total Chern classes of these derived sheaves,
\hskip -1in

$$\hskip -1in
c_{total}({\cal R}^1\pi'_{\ast}\bigl({\cal O}_{mE_B}(mE_B)\bigr)\otimes {\cal E})=
c_{total}({\cal R}^1\pi'_{\ast}\bigl({\cal O}_{(m-1)E_B}((m-1)E_B)\bigr)
 \otimes {\cal E})\cdot
 c_{total}({\bf S}^{m-2}({\cal N}_{s(B)}{\cal X})\otimes {\cal E}).$$

 This identity will be used in the proof of the theorem \ref{theo; al-bl}.
 This finishes the proof of lemma 
\ref{lemm; projective}.
$\Box$

  Firstly, we deal with the case when
 ${\cal R}^2\pi'_{\ast}\bigl({\cal E}\bigr)$, ${\cal R}^2\pi'_{\ast}\bigl(
 {\cal E}_m\bigr)$ vanish. It happens when the fiber 
surfaces have a trivial geometric genus
 or when $c_1(K_{{\cal X}/B}-C)$ is of negative relative degree.

  By using lemma \ref{lemm; projective}, we have 

$${\cal R}^0\pi'_{\ast}\bigl({\cal E}\bigr)\cong 
{\cal R}^0\pi'_{\ast}\bigl({\cal E}_m\bigr)$$

and

$$0\mapsto {\cal R}^1\pi'_{\ast}\bigl({\cal E}\bigr)\mapsto 
{\cal R}^1\pi'_{\ast}\bigl({\cal E}_m\bigr)
\mapsto {\cal R}^1\pi'_{\ast}\bigl({\cal O}_{mE_B}(mE_B)
\otimes{\cal E}\bigr)\mapsto 0.$$

 This suggests that one can build up algebraic Kuranishi models 
of $C, C+mE$ by
 $\Phi_{\cal VW}:{\cal V}\longrightarrow {\cal W}$,
 $\Phi_{\cal V'W'}:{\cal V'}\longrightarrow {\cal W}'$, ${\cal V'}={\cal V}$,
  with an short exact sequence relating the obstruction
 sheaves,

 $$0\mapsto {\cal W}\mapsto {\cal W}'\mapsto 
{\cal R}^1\pi'_{\ast}\bigl({\cal O}_{mE_B}(mE_B)\otimes{\cal E}\bigr)\mapsto 0.$$

  To achieve this, we choose a large $n$ and take 
 ${\cal R}^0\pi'_{\ast}\bigl({\cal O}(nD)\otimes {\cal E}_m\bigr)={\cal V}'$ and
 ${\cal R}^0\pi'_{\ast}\bigl({\cal O}_{nD}(nD)\otimes 
{\cal E}_m\bigr)={\cal W}'$, similar to 
subsection \ref{subsection; q>0}.

 Let ${\cal R}^0\pi'_{\ast}\bigl({\cal O}_{nD}(nD)\otimes {\cal E}_m\bigr)\mapsto 
 {\cal R}^1\pi'_{\ast}\bigl({\cal O}_{mE_B}(mE_B)\otimes {\cal E}\bigr)$ be the
 composition of two surjective maps

 $${\cal R}^0\pi'_{\ast}\bigl({\cal O}_{nD}(nD)\otimes {\cal E}_m\bigr)\mapsto 
{\cal R}^1\pi'_{\ast}\bigl({\cal E}_m\bigr)$$

 and $${\cal R}^1\pi'_{\ast}\bigl({\cal E}_m\bigr)\mapsto 
{\cal R}^1\pi'_{\ast}\bigl({\cal O}_{mE_B}(mE_B)\otimes {\cal E}\bigr).$$

\medskip
\begin{prop}
 Define ${\cal W}$ to be the kernel of 

$${\cal R}^0\pi'_{\ast}\bigl({\cal O}_{nD}(nD)\otimes {\cal E}_m\bigr)\mapsto 
 {\cal R}^1\pi'_{\ast}\bigl({\cal O}_{mE_B}(mE_B)\otimes {\cal E}\bigr) 
\mapsto 0.$$

 Then ${\cal W}$ is locally free and there is a four term exact sheaf sequence

$$\hskip -.3in 0\mapsto {\cal R}^0\pi'_{\ast}\bigl({\cal E}\bigr)\mapsto 
{\cal R}^0\pi'_{\ast}\bigl({\cal O}(nD)\otimes {\cal E}_m\bigr)\mapsto
 {\cal W}\mapsto {\cal R}^1\pi'_{\ast}\bigl({\cal E}\bigr)\mapsto 0.$$

\end{prop}

\medskip

Because of this exact sequence, one may take 
${\cal V}={\cal V}', {\cal W}$ to define a Kuranishi model of
 $C$.

\medskip

\noindent Proof:

 The sheaf ${\cal W}$ is locally free if and only if 
$\phi(x)=dim_{k(x)}{\cal W}_x\otimes k(x)$ is a 
 constant, which follows from the locally freeness of both 
${\cal R}^0\pi'_{\ast}\bigl({\cal O}_{nD}(nD)\otimes 
{\cal E}_m\bigr)$ and 
$ {\cal R}^1\pi'_{\ast}\bigl({\cal O}_{mE_B}(mE_B)\otimes 
{\cal E}\bigr)$ (see lemma
 \ref{lemm; projective}).

Because ${\cal R}^0\pi'_{\ast}\bigl({\cal E}\bigr)\cong 
{\cal R}^0\pi'_{\ast}\bigl({\cal E}_m\bigr)$, the
 injectivity of ${\cal R}^0\pi'_{\ast}\bigl({\cal E}\bigr)\mapsto 
{\cal R}^0\pi'_{\ast}\bigl({\cal O}(nD)\otimes {\cal E}_m\bigr)$ follows.

On the other hand, the exactness of 

$$\hskip -.3in 0\mapsto {\cal R}^0\pi'_{\ast}\bigl({\cal E}_m\bigr)\mapsto 
{\cal V}'\mapsto
 {\cal W}'\mapsto {\cal R}^1\pi'_{\ast}\bigl({\cal E}_m\bigr)\mapsto 0$$

 implies the composition ${\cal V}'\mapsto {\cal R}^1\pi'_{\ast}
\bigl({\cal E}_m\bigr)$ is trivial.
Thus ${\cal V}'\mapsto {\cal R}^1\pi'_{\ast}\bigl({\cal O}_{mE_B}(mE_B)
\otimes {\cal E}\bigr)$ is
 also trivial and then ${\cal V}'\mapsto {\cal W}'$ must factor through
 ${\cal W}\mapsto {\cal W}'$, the kernel of ${\cal W}'\mapsto 
{\cal R}^1\pi'_{\ast}\bigl({\cal O}_{mE_B}(mE_B)\otimes {\cal E}\bigr)$.

 This implies that there is a four term exact sequence

$$0\mapsto {\cal R}^0\pi'_{\ast}\bigl({\cal E}\bigr)\mapsto {\cal V}'
\stackrel{\Phi_{{\cal V}'{\cal W}}}{\longrightarrow} 
{\cal W}\mapsto Cokernel( \Phi_{{\cal V}'{\cal W}})\mapsto 0.$$

 It suffices to show that 
$Coker(\Phi_{{\cal V}'{\cal W}})\cong {\cal R}^1\pi'_{\ast}\bigl({\cal E}\bigr)$
 is a natural isomorphism. This follows from the following commutative diagram,

\[
\begin{array}{ccccc}
\hskip -.5in

  Coker(\Phi_{{\cal V}'{\cal W}}) & \mapsto & \Phi_{{\cal V}'{\cal W}'} & 
\mapsto  & 
{\cal R}^1\pi'_{\ast}\bigl({\cal O}_{mE_B}(mE_B)\otimes {\cal E}\bigr) \\
   &   & \Big\downarrow  &  & \Big\downarrow \\
{\cal R}^1\pi'_{\ast}\bigl({\cal E}\bigr)  & \mapsto  & 
{\cal R}^1\pi'_{\ast}\bigl({\cal E}_m\bigr)
  & \mapsto & 
{\cal R}^1\pi'_{\ast}\bigl({\cal O}_{mE_B}(mE_B)\otimes {\cal E}\bigr)
\end{array}
\]

\medskip
  Both rows are short exact sequences and the vertical arrows are natural 
isomorphisms.  Then
 $\Phi_{{\cal V}'{\cal W}'}\cong {\cal R}^1\pi'_{\ast}\bigl({\cal E}_m\bigr)$ 
induces the
 isomorphism $\Phi_{{\cal V}'{\cal W}}
\cong {\cal R}^1\pi'_{\ast}\bigl({\cal E}\bigr)$
$\Box$

 When one pulls back ${\bf W}'$, ${\bf W}$ to 
${\bf P}_{{\cal T}_B({\cal X})}({\bf V})$, one has
 the following identity on the top Chern class

 $$c_{top}({\cal H}\otimes {\cal W}')=c_{top}({\cal H}\otimes {\cal W})\cdot
 c_{top}({\cal H}\otimes 
({\cal R}^1\pi'_{\ast}\bigl({\cal O}_{mE_B}(mE_B)\otimes{\cal E}\bigr))=$$

$$c_{top}({\cal H}\otimes {\cal W})\cdot c_{top}({\cal H}\otimes 
(\oplus_{0\leq j\leq m-2} 
{\bf S}^j({\cal N}_{s(B)}{\cal X}))\otimes {\cal E}).$$

 We have inductively used the Chern class identity at the end of the proof of 
lemma \ref{lemm; projective}.

\medskip

 After finishing the vanishing 
${\cal R}^2\pi'_{\ast}\bigl({\cal E}\bigr)$ case, we have to
deal with ${\cal R}^2\pi'_{\ast}\bigl({\cal E}\bigr)$ non-vanishing case. 
\label{dsimilar}

  Because of the relative Serre duality 

$${\cal E}, {\cal E}_m={\cal E}\otimes {\cal O}(mE_B)\Longrightarrow 
{\cal K}_{{\cal X}'/B}\otimes {\cal E}^{\ast},
 {\cal K}_{{\cal X}'/B}\otimes ({\cal E}_m)^{\ast}=
{\cal K}_{{\cal X}'/B}\otimes {\cal E}^{\ast}\otimes {\cal O}(-mE_B),$$

  the case can be reduced to a parallel discussion with $m<0$. We will 
 address it at the end of the subsection on page
 \pageref{similar}.

\medskip

\noindent Case Two: We assume that $m\leq 0$. The $m=0$ case is trivial 
and we consider negative $m=-k, k\in
{\bf N}$.
We consider the following short exact sheaf sequence,

$$0\mapsto {\cal E}_m\mapsto {\cal E}\mapsto {\cal O}_{kE_B}\otimes {\cal E}
\mapsto 0.$$

  Its derived long exact sequence relates the derived images of 
 ${\cal E}_m$ and ${\cal E}$,

$$0\mapsto {\cal R}^0\pi'_{\ast}\bigl({\cal E}_m\bigr)\mapsto 
{\cal R}^0\pi'_{\ast}\bigl({\cal E}\bigr)\mapsto
{\cal R}^0\pi'_{\ast}\bigl( {\cal O}_{kE_B}\otimes {\cal E}\bigr)\mapsto $$
$${\cal R}^1\pi'_{\ast}\bigl({\cal E}_m\bigr)\mapsto 
{\cal R}^1\pi'_{\ast}\bigl({\cal E}\bigr)\mapsto
{\cal R}^1\pi'_{\ast}\bigl( {\cal O}_{kE_B}\otimes {\cal E}\bigr).$$

\medskip
 Similar to lemma \ref{lemm; projective}, we have the following lemma regarding 
${\cal R}^i\pi'_{\ast}\bigl( {\cal O}_{kE_B}\otimes {\cal E}\bigr)=
{\cal R}^i\pi'_{\ast}\bigl({\cal O}_{kE_B}\bigr)\otimes {\cal E}, i=0, 1$.

\medskip
\begin{lemm}\label{lemm; projective2}
For $k\in {\bf N}$, the first derived image sheaf
${\cal R}^1\pi'_{\ast}\bigl({\cal O}_{kE_B}\bigr)$ is the zero sheaf and 
${\cal R}^0\pi'_{\ast}\bigl({\cal O}_{kE_B}\bigr)$ is locally free. Moreover, 
 $c_{total}({\cal R}^0\pi'_{\ast}\bigl({\cal O}_{kE_B}\bigr))$ is equal to
 $c_{total}({\bf S}^{k-1}({\cal O}_B\oplus {\cal N}^{\ast}_{s(B)}{\cal X}))$.
\end{lemm}
\medskip
\noindent Proof:
 The argument of this lemma is very similar to lemma \ref{lemm; projective}.
 We consider the following short exact sequence 
 $$0\mapsto {\cal O}_{E_B}(-(k-1)E_B)\mapsto {\cal O}_{kE_B}\mapsto 
{\cal O}_{(k-1)E_B}\mapsto 0,$$
 for two divisors $A=(k-1)E_B$ and $B=E_B$ and their sum $A+B=kE_B$. 
Then we have to take the derived sequence
 and make an induction on $k$. 

 The key to identify the total Chern classes is to show that for 
$E_B={\bf P}_B({\bf N}_{s(B)}{\cal X})$,
 ${\cal R}^0\pi'_{\ast}\bigl({\cal O}_{E_B}(-(r-1)E_B)\bigr)=
{\bf S}^{r-1}({\cal N}^{\ast}_{s(B)}{\cal X})$,
 $r\in {\bf N}$, which follows from the projection formula 
(see [Ha] page 253 exercise 8.4.) 
of the ${\bf P}^1$ bundle. \label{c=}
We leave the remaining detail as an exercise to the reader, 
based on the argument of
 lemma \ref{lemm; projective}.
 $\Box$

If one works in the $C^{\infty}$ category, the derived long exact sequence 
of ${\cal E}_m, {\cal E}$ and
 lemma \ref{lemm; projective2} is enough for us to derive the desired family 
blowup formula stated in the theorem \ref{theo; al-bl}. To give an algebraic
 proof, we have to go through
 a detour.

 To simplify the discussion, firstly we assume (i). $p_g$ of the fiber 
surface to be $0$ or (ii). 
 $c_1(K_{{\cal X}/B})-C$ to be of non-positive degree with respect to an ample 
relative polarization. 
Under these assumptions, 
 the second derived sheave ${\cal R}^2\pi'_{\ast}\bigl({\cal E}\bigr)$ vanishes.

 By assumption, the fibration ${\cal X}'\mapsto B$ is relatively good. 
 Consider an effective ample divisor $D$ on the total space ${\cal X}'$, 
relative dimension one over $B$,
 and a large multiple 
$nD$ for $n\gg 0$. The same argument as in subsection \ref{subsection; q>0} 
allows us to construct explicit algebraic Kuranishi models
 for $C-kE$ and $C$ over ${\cal T}_B({\cal X})$. The key idea is to relate
 the algebraic Kuranishi models
 of $C-kE$ and $C$.

 Consider the short exact sequence 

 $$0\mapsto {\cal E}_m\otimes {\cal O}(nD)\mapsto {\cal E}\otimes {\cal O}(nD)
\mapsto {\cal O}_{kE_B}(nD)\otimes
 {\cal E}\mapsto 0.$$

  By choosing $n$ large enough one can make the first derived image sheaves 
 ${\cal R}^1\pi'_{\ast}\bigl({\cal O}(nD)\otimes {\cal E}_m\bigr)$,
 ${\cal R}^1\pi'_{\ast}\bigl({\cal O}(nD)\otimes 
{\cal E}\bigr)$ vanish.

 Thus, by lemma \ref{lemm; projective2} one gets a short exact sequence of
 locally free sheaves,

$$0\mapsto{\cal R}^0\pi'_{\ast}\bigl({\cal E}_m\otimes {\cal O}(nD)\bigr)\mapsto  
{\cal R}^0\pi'_{\ast}\bigl({\cal E} \otimes {\cal O}(nD)\bigr)
\mapsto {\cal R}^0\pi'_{\ast}\bigl({\cal O}_{kE_B}(nD))\otimes 
{\cal E}\mapsto 0. $$

On the other hand, one may get a short exact sequence on the obstruction
 sheaves of the algebraic 
Kuranishi models, 

$$\hskip -.5in 
0\mapsto{\cal R}^0\pi'_{\ast} \bigl({\cal O}_{nD}(nD)\otimes {\cal E}_m\bigr)
\mapsto {\cal R}^0\pi'_{\ast}\bigl({\cal O}_{nD}(nD)\otimes {\cal E}\bigr)\mapsto 
{\cal R}^0\pi'_{\ast}\bigl({\cal 0}_{kE_B\cap nD}(nD)\bigr)\otimes {\cal E}
\mapsto 0,$$

 due to the vanishing of higher derived image sheaves 
 ${\cal R}^i\pi'_{\ast}\bigl( {\cal O}_{nD}(nD)\otimes {\cal E}_m \bigr)$, 
${\cal R}^i\pi'_{\ast}\bigl( {\cal O}_{nD}(nD)\otimes {\cal E} \bigr)$, $i\geq 1$.
\medskip

 The following short exact sequence on the locally free sheaves 
is vital to the remaining discussion.

$$\hskip -.5in 
0\mapsto {\cal R}^0\pi'_{\ast}\bigl({\cal O}_{kE_B}\bigr)\otimes {\cal E}\mapsto
 {\cal R}^0\pi'_{\ast}\bigl( {\cal O}_{kE_B}(nD)\bigr)\otimes {\cal E}
\mapsto {\cal R}^0\pi'_{\ast}\bigl( {\cal 0}_{kE_B\cap nD}(nD)\bigr)\otimes
 {\cal E}\mapsto  0.$$

 Consider the following lemma,

\medskip

\begin{lemm}\label{lemm; stablization}
  Let ${\bf V}$ be an algebraic vector bundle over ${\cal T}_B({\cal X})$ 
and $0\mapsto {\bf V}\mapsto {\bf V}'\mapsto {\bf U}\mapsto 0$ be a bundle
 extension 
 of ${\bf V}$ by ${\bf U}$. 
Denote ${\bf H}^{\ast}$ to be the tautological line bundle over 
${\bf P}_{{\cal T}_B({\cal X})}({\bf V}')$.
 Then the projective space bundle 
${\bf P}_{{\cal T}_B({\cal X})}({\bf V})$ can be identified with the subspace of
 ${\bf P}_{{\cal T}_B({\cal X})}({\bf V}')$ defined by the zero locus of a
 canonical section
 of ${\bf H}\otimes {\bf U}$.
\end{lemm}
\medskip
\noindent Proof:
 Consider the bundle surjection ${\bf V}'\mapsto {\bf U}$, it induces a map
 from the
 tautological line bundle ${\bf H}^{\ast}$ to ${\bf U}$ and therefore a 
canonical 
section of ${\bf H}\otimes {\bf U}$
 over ${\bf P}_{{\cal T}_B({\cal X})}({\bf V}')$.
 On the other hand, the fibers of ${\bf H}^{\ast}$ can be identified with
 the rays of ${\bf V}'$, which
 map trivially to ${\bf U}$ if and only if the rays are from the sub-bundle 
${\bf V}$. A direct
 investigation shows that the ${\bf P}_{{\cal T}_B({\cal X})}({\bf V})$ is 
the transversal zero locus of
 the canonical section of ${\bf H}\otimes {\bf U}$. $\Box$

\medskip

  This lemma implies that one may thicken the projective space bundle
 by adding ${\bf H}\otimes {\bf U}$ to the obstruction bundle.

 Now one compares the algebraic Kuranishi models of $C-kE$ and $C$. Let 
 ${\bf V}$, ${\bf V}'$ and ${\bf U}$ be the algebraic vector 
bundles associated with the
 locally free sheaves ${\cal R}^0\pi'_{\ast}\bigl({\cal E}_m\otimes
 {\cal O}(nD)\bigr)$,
 ${\cal R}^0\pi'_{\ast}\bigl({\cal E}\otimes {\cal O}(nD)\bigr)$ and
 ${\cal R}^0\pi'_{\ast}\bigl({\cal O}_{kE_B}(nD))\otimes 
{\cal E}$, respectively. 

  Then there is a bundle short exact sequence $0\mapsto {\bf V}\mapsto
 {\bf V}'\mapsto {\bf U}_V\mapsto 0$.
 The original algebraic Kuranishi model of $C-kE$ and $C$ are realized 
as intersection numbers in 
 (open sets of) ${\bf P}_{{\cal T}_B({\cal X})}({\bf V})$ and 
${\bf P}_{{\cal T}_B({\cal X})}({\bf V}')$.
By lemma \ref{lemm; stablization}, One may embed 
${\bf P}_{{\cal T}_B({\cal X})}({\bf V})$ in 
 ${\bf P}_{{\cal T}_B({\cal X})}({\bf V}')$ as the transversal
 zero locus of ${\bf H}\otimes {\bf U}$.

 Thus, to compare the algebraic family invariants of $C-kE$ and $C$, 
it suffices to compare the 
 obstruction bundles ${\bf W}'$ and ${\bf W}$ and their canonical sections.
 Let ${\bf W}, {\bf W}'$ and $\tilde{\bf U}$ be the algebraic vector bundle 
over ${\cal T}_B({\cal X})$
 associated with ${\cal R}^0\pi'_{\ast} \bigl({\cal O}_{nD}(nD)\otimes 
{\cal E}_m\bigr)$,
 ${\cal R}^0\pi'_{\ast} \bigl({\cal O}_{nD}(nD)\otimes {\cal E}\bigr)$ and 
${\cal R}^0\pi'_{\ast}\bigl({\cal 0}_{kE_B\cap nD}(nD)\bigr)\otimes {\cal E}$.
 Then there is also a bundle short exact sequence 

 $$0\mapsto {\bf W}\mapsto {\bf W}'\mapsto {\bf U}_W\mapsto 0$$
 and the following diagram is commutative

 \[
  \begin{array}{cccc}
  \Phi_{\bf VW}: & {\bf V}  &  \longrightarrow & {\bf W} \\
   & \Big\downarrow &  & \Big\downarrow \\
  \Phi_{\bf V'W'}: & {\bf V}'  & \longrightarrow  &  {\bf W}'  
  \end{array}
\]

  While the family moduli space of curves in $C-kE$ and in $C$ 
 are both viewed as subspaces in ${\bf P}_{{\cal T}_B({\cal X})}({\bf V}')$,
 the obstruction bundles of the class $C-kE$ and $C$ are ${\bf H}
\otimes ({\bf W}\oplus {\bf U}_V)$, 
${\bf H}\otimes {\bf W}'$, respectively.

  One see easily that 
 $$[{\bf W}] + [{\bf U}] = [{\bf W}']+[{\bf U}_V]-[{\bf U}_W]$$

 in the $K$ group of ${\cal T}_B({\cal X})$. By using 

 $$0\mapsto  {\bf E}\otimes ({\bf S}^{k-1}({\bf C}_B\oplus 
{\bf N}^{\ast}_{s(B)}{\cal X}))
  \mapsto {\bf U}\mapsto \tilde {\bf U}\mapsto 0,$$

  one may conclude 

$$[{\bf W}] + [{\bf U}_V] = [{\bf W}']+
[{\bf E}\otimes ({\bf S}^{k-1}({\bf C}_B\oplus 
{\bf N}^{\ast}_{s(B)}{\cal X})))],$$

$$rank_{\bf C}{\bf W}+rank_{\bf C} {\bf U} = rank_{\bf C} {\bf W}' +
rank_{\bf C} {\bf S}^{k-1}({\bf C}_B\oplus {\bf N}^{\ast}_{s(B)}{\cal X}).$$

 and the following identity on top Chern classes

$$c_{top}({\bf H}\otimes ({\bf W}\oplus {\bf U}_V))=c_{top}({\bf H}\otimes 
{\bf W}')\cdot
 c_{top}({\bf H}\otimes 
{\bf E}\otimes ({\bf S}^{k-1}({\bf C}_B\oplus {\bf N}^{\ast}_{s(B)}{\cal X}))).$$

 The identity on ${\cal AFSW}$  of $C$ and $C-kE$ follows from the equality of 
top Chern classes.                            

\medskip

  Finally let us derive the family blowup formula for $p_g>0$ fiber surfaces 
with a fiberwise invariant 
 class $C$ with positive relative degree on $c_1(K_{{\cal X}/B})-C$. The 
derivation is much more evolved than the
 previous case as several derived image sheaves fail to be locally free in 
this case.

 By the assumption of the theorem, the fibration ${\cal X}'\mapsto B$ is 
assumed to be two-relatively good. Namely, there 
 are two ample divisors $D_1, D_2$ which are (1). of relative dimension 
one over $B$.   (2). The intersection 
$D_1\cap D_2$ is of relative dimension zero over $B$.

  Following subsection \ref{subsection; pg>0},
 we adopt the asymmetric notation $D=D_1$ and $\Delta=lD_2$. 

\medskip

\begin{rem}\label{rem; difference}
 Basically we have to extend the construction in subsection 
\ref{subsection; pg>0} to a relative version.
 The essential difference is that in the $B=pt$ case, we may choose a 
single $\Sigma\in |nD|$ to be smooth. 
 Then the vanishing of certain derived image sheaves on $\Sigma\times 
T(M)\mapsto T(M)$ follows from Kodaira vanishing
 theorem on $\Sigma$.  
When
 ${\cal X}'\mapsto B$ is a non-trivial fiber bundle, it may be hard to 
find an effective
 divisor in $|nD|$ which is relatively smooth of
 dimension one over $B$ unless one makes some additional assumption on 
${\cal X}'\mapsto B$. 
 But the vanishing result on the first derived image sheaves of 
sufficiently very ample invertible sheaves 
 on a non-smooth (or even non-reduced) relative dimension one divisor
 $nD$ can still
 be derived by using the exact sequence $0\mapsto 
{\cal O}\mapsto {\cal O}(nD)\mapsto {\cal O}_{nD}(nD)\mapsto 0$
 suitably.
  
\end{rem}

\medskip

We have the following commutative
 diagram of coherent sheaves for a sufficiently large $n$, and $l$.

 \[
\hskip -1.5in  
\begin{array}{ccccc}
 {\cal R}^0\pi'_{\ast}\bigl( {\cal O}_{nD}(nD)\otimes {\cal E}_m \bigr)  &\mapsto& 
{\cal R}^0\pi'_{\ast}\bigl( {\cal O}_{nD}(nD)\otimes {\cal E} \bigr) &
\mapsto
 &{\cal R}^0\pi'_{\ast}\bigl( {\cal O}_{kE_B\cap nD}(nD)\otimes {\cal E}\bigr) \\
  \Big\downarrow & & \Big\downarrow &  &\Big\downarrow \\
 {\cal R}^0\pi'_{\ast}\bigl( {\cal O}_{nD}(nD+\Delta)\otimes {\cal E}_m \bigr) & 
\mapsto& {\cal R}^0\pi'_{\ast}\bigl( {\cal O}_{nD}(nD+\Delta)\otimes
 {\cal E} \bigr) & \mapsto
  & {\cal R}^0\pi'_{\ast}\bigl( {\cal O}_{kE_B\cap nD}(nD+\Delta)\otimes
 {\cal E} \bigr)\\
 \Big\downarrow & &\Big\downarrow & & \Big\downarrow\\
 {\cal R}^0\pi'_{\ast}\bigl( {\cal O}_{nD\cap \Delta}(nD+\Delta)\otimes 
{\cal E}_m\bigr) & 
\mapsto&{\cal R}^0\pi'_{\ast}\bigl( {\cal O}_{nD\cap \Delta}(nD+\Delta)
\otimes {\cal E}\bigr) & \mapsto
 & {\cal R}^0\pi'_{\ast}\bigl( {\cal O}_{nD\cap kE_B\cap \Delta}(nD+\Delta)
\otimes {\cal E}\bigr) \\
 \Big\downarrow & & \Big\downarrow &  &\Big\downarrow \\
{\cal R}^1\pi'_{\ast}\bigl( {\cal O}_{nD}(nD)\otimes {\cal E}_m  \bigr)  & \mapsto&
{\cal R}^1\pi'_{\ast}\bigl( {\cal O}_{nD}(nD)\otimes {\cal E} \bigr) & 
\mapsto& {\cal R}^1\pi'_{\ast}\bigl( {\cal O}_{kE_B\cap nD}(nD)\otimes 
{\cal E} \bigr)
  \end{array}
\]

\label{comm}

In the commutative 
diagram there are twelve derived image sheaves and form four rows and three 
columns. 
The first row and the fourth
 row are connected by the connecting homomorphism, which is not 
included in the diagram. Thus, the first row is left 
exact but generally not right exact.
 Because $D\mapsto B$ is
 of relative dimension one, by base change one sees that the the 
fourth row is right exact. On the other hand, the 
 three columns are parts of derived long exact sequences from sheaves 
short exact sequences. If we choose $l$
 and $n$ to be large enough, the three columns are four terms exact 
sequences.
 This also implies that the second row is a short exact sequence.  
Because $nD\cap \Delta$ is relative dimension one
(non-reduced) over the base, the third row is also short exact.
 
  Similar to the arguments of lemma \ref{lemm; projective} and lemma 
\ref{lemm; projective2}, 
the derived image sheaves in the second and the third rows of the 
diagram are
 locally free sheaves. The $(2,1)$th, $(3, 1)$th, $(2, 2)$th, and 
$(3, 2)$th entries
 are used in constructing the algebraic Kuranishi models of $C-kE$ and $C$.

\medskip

\noindent Proof of the Case with non-zero 
${\cal R}^2\pi'_{\ast}\bigl({\cal E}\bigr)$: 

Denote the vector bundle associated with the $(2, 1)$th, $(2, 2)$th,
 $(3, 1)$, $(3, 2)$th and $(2, 3)$th entries of the sheaves commutative 
diagram on page \pageref{comm}
 by $\tilde{\bf W}, \tilde{\bf W}'$, 
$\tilde{\bf V}$, $\tilde{\bf V}'$ and ${\bf R}$, respectively. Then one 
may construct
 ${\bf W}$ and ${\bf W}'$ by using the recipe of proposition \ref{prop; exist}.
  Denote the algebraic vector bundles 
associated with the locally free ${\cal R}^0\pi'_{\ast}\bigl({\cal O}(nD)
\otimes {\cal E}_m\bigr)$, 
 ${\cal R}^0\pi'_{\ast}\bigl({\cal O}(nD)\otimes {\cal E}\bigr)$ by 
${\bf V}$, ${\bf V}'$.  Then 
 a parallel discussion similar to proposition \ref{prop; exist}
 implies that the algebraic family 
Kuranishi models of $C-kE$ and $C$ can be 
built from the algebraic bundle maps 

$$\Phi_{{\bf V}{\bf W}}:{\bf V}\longrightarrow {\bf W},$$ 

 and 

$$\Phi_{{\bf V}'{\bf W}'}:{\bf V}'\longrightarrow {\bf W}'.$$
\label{mark}

  We concentrate on how does the relative obstruction bundle of the 
family blowup formula appear in the current picture. 

 Firstly, notice that there are short exact sequences 

$$0\mapsto {\bf V}\mapsto {\bf V}'\mapsto {\bf U}_V\mapsto 0,$$ 

$$0\mapsto \tilde{\bf V}\mapsto \tilde{\bf V}'\mapsto \tilde{\bf U}_V\mapsto 0,$$

$$0\mapsto \tilde{\bf W}\mapsto \tilde{\bf W}'\mapsto \tilde{\bf U}_W\mapsto 0,$$

 with ${\bf U}_V$, $\tilde{\bf U}_V$
 being the algebraic vector bundle associated with the locally free sheaves
${\cal R}^0\pi'_{\ast}\bigl( {\cal O}_{kE_B}(nD)\otimes
 {\cal E} \bigr)$ and 
${\cal R}^0\pi'_{\ast}\bigl({\cal O}_{nD\cap kE_B\cap 
\Delta}(nD+\Delta)\otimes {\cal E}\bigr)$. 

Recall that $[{\bf W}]=[\tilde{\bf W}]-[\tilde{\bf V}]+[\tilde{\bf F}]$ and
    $[{\bf W}']=[\tilde{\bf W}']-[\tilde{\bf V}']+[\tilde{\bf F}']$ in the
 $K$ group of algebraic vector bundles on 
${\cal T}_B({\cal X}')$, by remark \ref{mem; K}.

It is easy to see that for the formal excess base
 dimensions, we have 
$febd(C, {\cal X}'/B)=febd(C-mE, {\cal X}'/B)$ by using the isomorphism 
 ${\cal R}^2\pi'_{\ast}\bigl({\cal E}_m\bigr)\cong 
{\cal R}^2\pi'_{\ast}\bigl({\cal E}\bigr)$ and then
 ${\bf F}\cong {\bf F}'$. (read the paragraph before prop. 
\ref{prop; analogue} for the definition of 
 the corresponding locally free sheaf ${\cal F}$)

  Again by lemma \ref{lemm; stablization} and the same 
argument in the previous case, it suffices to show the following
 identity on
 the virtual bundles $$[{\bf U}_V\oplus \tilde{\bf U}_V-\tilde{\bf U}_W]
=[{\bf E}\otimes_{\bf C} ({\bf S}^{k-1}({\bf C}_B\oplus 
{\bf N}^{\ast}_{s(B)}{\cal X}))]$$ 
in the $K$ group of algebraic vector bundles on ${\cal T}_B({\cal X})$ and
 the corresponding equality on their virtual ranks.

The follows from the following calculation in the $K$ group of coherent sheaves
 $K_0({\cal T}_B({\cal X}))$ and the short exact sequences,

\label{reduce}

 $$\hskip -.9in -{\cal R}^0\pi'_{\ast}\bigl({\cal O}_{nD\cap kE_B}(nD+\Delta)
\otimes {\cal E}\bigr)+
{\cal R}^0\pi'_{\ast}\bigl( {\cal O}_{kE_B}(nD)\otimes
 {\cal E} \bigr)\oplus 
  {\cal R}^0\pi'_{\ast}\bigl({\cal O}_{nD\cap kE_B\cap \Delta}(nD+\Delta)
\otimes {\cal E}\bigr)$$
$$\hskip -.9in 
=-\bigl(  {\cal R}^0\pi'_{\ast}\bigl( {\cal O}_{nD\cap kE_B}(nD+\Delta)
\otimes {\cal E}\bigr)-
 {\cal R}^0\pi'_{\ast}\bigl({\cal O}_{nD\cap kE_B\cap \Delta}(nD+\Delta)
\otimes {\cal E}\bigr)\bigr)
+{\cal R}^0\pi'_{\ast}\bigl( {\cal O}_{kE_B}(nD)\otimes
 {\cal E} \bigr)$$
$$=-{\cal R}^0\pi'_{\ast}\bigl( {\cal O}_{nD\cap kE_B}(nD)\otimes {\cal E} \bigr)+
{\cal R}^0\pi'_{\ast}\bigl( {\cal O}_{kE_B}(nD)\otimes {\cal E} \bigr)$$
$$={\cal R}^0\pi'_{\ast}\bigl( {\cal O}_{kE_B}\otimes {\cal E} \bigr)-
{\cal R}^1\pi'_{\ast}\bigl( {\cal O}_{kE_B}\otimes {\cal E} \bigr) $$
$$={\cal R}^0\pi'_{\ast}\bigl( {\cal O}_{kE_B}\otimes {\cal E} \bigr)=
{\cal R}^0\pi'_{\ast}\bigl({\cal O}_{kE_B} \bigr) \otimes {\cal E}.$$
 
 On page \pageref{c=}, lemma \ref{lemm; projective2}, we have
already identified the 
the total Chern classes of the two coherent sheaves 
${\cal R}^0\pi'_{\ast}\bigl( 
{\cal O}_{kE_B} \bigr)$ and 
the symmetric power ${\bf S}^{k-1}({\cal O}_B\oplus 
{\cal N}^{\ast}_{s(B)}{\cal X})$, so
 the proof of $m<0$, ${\cal R}^2\pi'_{\ast}\bigl({\cal E}\bigr)\not=0$
 case is done.
\label{similar}

 At the end of the proof, 
let us address the $m>0$, ${\cal R}^2\pi'_{\ast}\bigl({\cal E}\bigr)\not=0$
 case leftover on
 page \pageref{dsimilar}. As in the above discussion, ${\cal X}'\mapsto B
$ is assumed to 
be two relatively good.
 We have chosen $nD$ and $\Delta=lD_2|_{nD}$ to construct the algebraic 
family Kuranishi models.

 By interchanging the roles of ${\cal E}$ and ${\cal E}_m$ symbolically and by 
replacing $k$ by $m$, 
there is a corresponding twelve-term
commutative exact diagram similar to the one on page \pageref{comm}.
 Following the convention on page \pageref{mark}, we take
 ${\bf V}$, $\tilde{\bf W}$, $\tilde{\bf V}$ and ${\bf V}'$, $\tilde{\bf W}'$,
 $\tilde{\bf V}'$ to be  
 the algebraic bundles associated to the algebraic 
family Kuranishi models of $C$ and $C+mE$, respectively. 

 Then as before there are exact sequences 

$$0\mapsto {\bf V}\mapsto {\bf V}'\mapsto {\bf U}_V\mapsto 0,$$
 
$$0\mapsto \tilde{\bf V}\mapsto \tilde{\bf V}'\mapsto \tilde{\bf U}_V\mapsto 0,$$
 
and $$0\mapsto \tilde{\bf W}\mapsto \tilde{\bf W}'\mapsto \tilde{\bf U}_W
\mapsto 0.$$

 In this case, ${\bf F}\cong {\bf F}'$ may not hold, but we still have 
 $[{\bf F}]=[{\bf F}']$ in the reduced $K$ group as both of their associated 
locally free sheaves
 ${\cal F}$ and ${\cal F}'$ are 
 constructed from  ${\cal R}^2\pi'_{\ast}\bigl({\cal O}_{{\cal X}'}\bigr)$ by
 quotienting out some trivial sub-factors ${\cal O}^{febd(C, 
{\cal X}'/B)}_{{\cal T}_B({\cal X}')}$ and
${\cal O}^{febd(C+mE, {\cal X}'/B)}_{{\cal T}_B({\cal X}')}$, respectively.

 Our goal is to prove that, 

\medskip

(i). $-[\tilde{\bf U}_W-\tilde{\bf U}_V-{\bf U}_V]=-[{\bf E}\otimes 
{\bf S}^{m-2}({\bf C}_B\oplus {\bf N}_{s(B)}{\cal X})]$ in the $K$ group.

\medskip

(ii). The virtual ranks of the virtual vector bundles in (i). match.

\medskip

 We can go through the same calculation above on page \pageref{reduce}, 
except that 
the invertible sheaf ${\cal E}_m$ which replaces ${\cal E}$
 cannot be pulled out of the derived images of $\pi'$.

 The final answer is 

$${\cal R}^0\pi'_{\ast}\bigl({\cal O}_{mE_B}(mE_B)\otimes {\cal E}\bigr)-
 {\cal R}^1\pi'_{\ast}\bigl({\cal O}_{mE_B}(mE_B)\otimes {\cal E}\bigr)$$

 by lemma \ref{lemm; projective}, 

$$=-{\cal R}^1\pi'_{\ast}\bigl({\cal O}_{mE_B}(mE_B)\otimes {\cal E}\bigr)=
-{\cal R}^1\pi'_{\ast}\bigl({\cal O}_{mE_B}(mE_B)\bigr)\otimes {\cal E}.$$

\medskip

this ends the proof of
 the theorem in the ${\cal R}^2\pi'_{\ast}\bigl({\cal E}\bigr)\not=0$, $m<0$ case
 and therefore the proof of the whole theorem. $\Box$

\medskip

\subsection{\bf The Family Blowup Formula of the Universal Family}
\label{subsection; universal}

\bigskip

 In the previous section, section \ref{section; ap},
 we have derived the family blowup formula of 
the algebraic family invariant. 
 When we prove the formula, 
we have to construct algebraic family Kuranishi models based on 
relatively very ample divisors
 on these families. The readers should be aware that the 
non-uniqueness of the 
choices of the very ample divisors and 
therefore the choices of algebraic family Kuranishi models.
In this sub-section, we focus on the universal family 
$M_{n+1}\mapsto M_n$ and make a brief remark 
 regarding the 'canonical obstruction bundles' for these special families.
 Following the notation in [Liu1], we would like to explain 
how does the blowup formula relate to the
 canonical obstruction bundles of these families.

 Let ${\cal E}$ be the invertible sheaf over $M\times T(M)$ 
whose first Chern class over
 $M\times \{t\}, t\in T(M)$ is $C\in H^{1, 1}(M, {\bf Z})$.
 One may apply the blowing down map $T(M)\times M_{n+1}\mapsto 
T(M)\times M\times M_n={\cal T}_{M_n}(M_{n+1})$
to pull ${\cal E}$  back, which defines an invertible sheaf over 
$T(M)\times M_{n+1}$. We abuse the notation and
 denoting them by the same symbol. Let $E_1, E_2, \cdots, E_n$ 
denote the exceptional divisors of
 the blowing down map $M_{n+1}\mapsto M\times M_n$. Then we have 
the following short exact sequence 

$$0\mapsto {\cal E}(-\sum m_i E_i)\mapsto {\cal E}\mapsto 
{\cal O}_{\sum m_i E_i}\otimes {\cal E}\mapsto 0.$$

  Taking the right derived images of this sequence along 
$T(M)\times M_{n+1}\mapsto T(M)\times M_n$, then we get a 
long exact 
sequence of coherent sheaves on $T(M)\times M_n$. In the long 
paper[Liu1], our main concern is study the
 case when ${\cal E}$ is sufficiently very ample. So let us 
make a simplifying assumption on ${\cal E}$,

\medskip
\begin{assum} \label{assum; vani}
 The invertible sheaf ${\cal E}\otimes {\cal K}^{-1}_M$ is 
ample.
\end{assum} 

 By Nakai criterion, this assumption is reduced to an 
assumption of the positivity of $C-c_1(K_M)$ on the
 curve cone of $M$.
 If ${\cal E}$ is a high power of an ample invertible 
sheaf, the assumption always holds.
 Under this assumption, the higher sheaf cohomologies 
of ${\cal E}$ vanish. Then the sheaf short exact sequence
 induces a derived long exact sequence

$$0\mapsto {\cal R}^0(f_n)_{\ast}\bigl({\cal E}(-\sum m_i E_i)\bigr)\mapsto 
 {\cal R}^0(f_n)_{\ast}\bigl( {\cal E}\bigr)\mapsto 
{\cal R}^0(f_n)_{\ast}\bigl( {\cal O}_{\sum_{i\leq n}
 m_i E_i}\otimes {\cal E} \bigr)$$
$$\mapsto {\cal R}^1(f_n)_{\ast}\bigl({\cal E}(-\sum_{i\leq n}
 m_i E_i)\bigr)\mapsto 0.$$

 Because ${\cal E}$ is pulled back from $M\times T(M)$ to the 
universal family, 
${\cal R}^0(f_n)_{\ast}\bigl( {\cal E}\bigr)$ is pulled back 
from $T(M)$ and is constant along the $M_n$ factor.

 By the similar argument as in lemma \ref{lemm; projective2},
 we can prove the following,
\medskip

\begin{prop}\label{prop; canonical}
 The coherent sheaf ${\cal R}^0(f_n)_{\ast}\bigl( 
{\cal O}_{\sum_{i\leq n} m_i E_i}\otimes {\cal E}\bigr)$
 is locally free of
 rank $\sum_{i\leq n} {(m_i+1)m_i\over 2}$.
\end{prop}
\medskip
\noindent Proof:  Because the sheaf is the zero-th derived
 image of a relatively one dimensional fibration, 
 it suffices to prove inductively that the first derive
 image sheaf vanishes.

 This follows from the following short exact sequence 

$$0\mapsto {\cal O}_{m_nE_n}(-\sum_{i\leq n-1}m_iE_i)\otimes 
{\cal E}\mapsto {\cal O}_{\sum_{i\leq n}m_i E_i}\otimes {\cal E}\mapsto 
 {\cal O}_{\sum_{i\leq n-1} m_i E_i}\otimes {\cal E}\mapsto 0.$$
 
  The remaining of the proof is very similar to the lemmas 
\ref{lemm; projective} and \ref{lemm; projective2}, 
we omit the details.
$\Box$

 As in [Liu1], the $n+1$th universal space $M_{n+1}$ can be constructed 
by blowing up the relative diagonal 
$M_n\mapsto M_n\times_{M_{n-1}}M_n$ from the fiber product

\[
\begin{array}{ccc}
 M_n\times_{M_{n-1}}M_n  & \longrightarrow & M_n \\
 \Big\downarrow  &  & \Big\downarrow\vcenter{%
  \rlap{$\scriptstyle{\mathrm{f_n}}\,$}}    \\
  M_n & \stackrel{\mathrm{f_n}}{\longrightarrow}& M_{n-1}
\end{array}
\]
 The space $M_n\times_{M_{n-1}}M_n$ fibers over $M_n$ and the map 
$M_n\times_{M_{n-1}}M_n\mapsto M_n$ is
 smooth of relative dimension two. 
 One may view the relative diagonal as a cross section from the base 
$M_n$ to $M_n\times_{M_{n-1}}M_n$ and the
 normal bundle of the cross section is isomorphic to the relative tangent
 bundle of $M_n\mapsto M_{n-1}$,
 ${\bf T}_{M_n/M_{n-1}}\cong {\bf T}M_n/f_n^{\ast}{\bf T}M_{n-1}$.

 Knowing that ${\cal R}^0(f_n)_{\ast}\bigl({\cal O}_{\sum_{i\leq n} m_i E_i}
\otimes {\cal E}\bigr)$
 is locally free, 
its associated vector bundle can be used to build the canonical obstruction 
bundle of the family 
Seiberg-Witten invariant
(under the additional
 assumption \ref{assum; vani}).

\medskip

\begin{defin}\label{defin;canon}
Under the assumption \ref{assum; vani} on ${\cal E}$, 
let ${\bf V}_{canon}$ and  ${\bf W}_{canon}$ denote the algebraic 
vector bundles associated 
to the locally free ${\cal R}^0(f_n)_{\ast}\bigl({\cal E}\bigr)$ and 
${\cal R}^0(f_n)_{\ast}\bigl({\cal O}_{\sum_{i\leq n} m_i E_i}\otimes 
{\cal E}\bigr)$.

 Then for a given tuple $(m_1, m_2, m_3, \cdots, m_n)$ of singular 
multiplicities,
 the tuple 
 $$({\bf V}_{canon}, {\bf W}_{canon}, 
\Phi_{{\bf V}_{canon}{\bf W}_{canon}})$$

 with the bundle map

$$\Phi_{{\bf V}_{canon}{\bf W}_{canon}}:{\bf V}_{canon}\mapsto 
{\bf W}_{canon}$$

induced from $${\cal R}^0(f_n)_{\ast}\bigl({\cal E}\bigr)\mapsto 
{\cal R}^0(f_n)_{\ast}\bigl({\cal O}_{\sum_{i\leq n} m_i E_i}\otimes 
{\cal E}\bigr)$$

 is called the canonical algebraic family Kuranishi model of 
$C-\sum_{i\leq n}m_i E_i$.
 
 The vector bundle 
${\bf Obs}_{canonical; C-\sum_{i\leq n}m_i E_i}=
{\bf H}\otimes \pi^{\ast}_{{\bf P}({\bf V}_{canon})} {\bf W}_{canon}$
 over the projective space bundle ${\bf P}({\bf V}_{canon})$
is defined to be the canonical obstruction bundle of the class $C-\sum_i m_i E_i$.

\end{defin}

\medskip

 Moreover, one see explicitly that the first term of the following
 short exact sequence relating the canonical obstruction bundles of 
$C-\sum_{i\leq n}m_i E_i$ and $C-\sum_{i\leq n-1}m_i E_i$,

$$\hskip -1.3in 0\mapsto {\bf H}\otimes \pi^{\ast}_{{\bf P}({\bf V}_{canon})}
 \bigl({\bf E}\otimes {\bf S}^{m_n-1}({\cal O}_B\oplus 
{\bf T}^{\ast}_{M_n/M_{n-1}})\bigr)
\mapsto {\bf Obs}_{canonical; C-\sum_{i\leq n}m_i E_i}\mapsto 
{\bf Obs}_{canonical; C-\sum_{i\leq n-1}m_i E_i}
\mapsto 0.$$

  This gives a transparent explanation why the mixed invariants
 ${\cal AFSW}_{M_{n+1}\mapsto M_n}(c, C-\sum_{i\leq n}m_iE_i)$ is 
equal to the combination of mixed invariants
  $\sum_{p\leq {m_n(m_n+1)\over 2} }{\cal AFSW}_{M_{n+1}\mapsto M_n}(c\cup 
c_p({\bf S}^{m_n-1}({\bf C}_B\oplus {\bf T}^{\ast}_{M_n/M_{n-1}})), 
C-\sum_{i\leq n-1} m_i E_i)$.

{}

\end{document}